\documentclass[preprint,12pt]{elsarticle}
\usepackage{amsmath,amssymb,amsfonts,afterpage,hyperref}
\usepackage[makeroom]{cancel}
\usepackage{mathtools}
\usepackage{amscd,subeqnarray}
\usepackage{comment}
\usepackage{graphics}
\usepackage{wrapfig}
\usepackage{graphicx}
\usepackage{epsfig}
\usepackage{array}
\usepackage{color}
\usepackage{lineno}
\usepackage{ragged2e}
\usepackage{stmaryrd}
\usepackage{soul, ulem}
\usepackage{siunitx}
\usepackage{boldfonts,bm}
\usepackage{bbding}
\usepackage{pifont}

\usepackage{geometry}
\usepackage{pgf}
\usepackage{pgfplotstable}
\usepackage{caption}
\usepackage{subcaption,booktabs}

\usepackage{placeins}
\usepackage{adjustbox}

\newtheorem{theorem}{Theorem}[section]
\newtheorem{lemma}{Lemma}[section]

\newtheorem{remark}{Remark}[section]

\newcommand{\vertiii}[1]{{\left\vert\kern-0.25ex\left\vert\kern-0.25ex\left\vert #1
    \right\vert\kern-0.25ex\right\vert\kern-0.25ex\right\vert}}

\definecolor{mypink}{RGB}{219, 48, 122}

\newcommand{\pp}[2]{\dfrac{\partial {#1}}{\partial {#2}}}

\newcommand{\divv}{\nabla \cdot}

\newcommand{\normz}[1]{\| #1\|} 
\newcommand{\normo}[1]{\| #1\|_1}
\newcommand{\normr}[1]{\| #1\|_{\mathbb{R}^r}}

\newcommand{\triplenorm}[1]{%
 \left\vert\kern-0.9pt\left\vert\kern-0.9pt\left\vert #1
 \right\vert\kern-0.9pt\right\vert\kern-0.9pt\right\vert}

\newcommand{\Om}{\Omega}

\newcommand{\bff}{\mathbf{f}}

\newcommand{\Honez}{H^1_{0}(\Omega)}
\newcommand{\Honedz}{[H^1_{0}(\Omega)]^d}

\newcommand{\Vh}{\boldsymbol{V}_h}
\newcommand{\Wh}{W_h}
\newcommand{\Thetah}{{\Theta}_h}
\newcommand{\Vr}{\boldsymbol{V}_r}
\newcommand{\Wr}{W_r}
\newcommand{\Thetar}{{\Theta}_r}

\newcommand{\Th}{\mathcal{T}_h}

\newcommand{\czero}{c_0}
\newcommand{\cd}{C_d}
\newcommand{\alpham}{\alpha_m}
\newcommand{\alphat}{\alpha_T}
\newcommand{\Kdr}{K_{dr}}
\newcommand{\thetaz}{\theta_0}

\newcommand{\buh}{\bu_h}
\newcommand{\ph}{p_h}
\newcommand{\thetah}{\theta_h}
\newcommand{\bur}{\bu_r}
\newcommand{\pr}{p_r}
\newcommand{\thetar}{\theta_r}

\newcommand{\bsigt}{\tilde{\bsigma}}
\newcommand{\bsigb}{\overline{\bsigma}}

\newcommand{\dt}{\Delta t}

\newcommand{\eu}{e_{\bu}}
\newcommand{\ep}{e_p}
\newcommand{\etheta}{e_\theta}
\newcommand{\half}{\frac{1}{2}}

\newcommand{\mat}[1]{\mathbf{#1}}
\renewcommand{\vec}[1]{\underline{\boldsymbol{\mathbf{#1}}}}

\graphicspath{{images}}

\usetikzlibrary{patterns}
\usetikzlibrary{calc}
\tikzset{axis line style/.style={thin, gray, -stealth}}

\begin{document}

\begin{frontmatter}

\title{Projection-based reduced order modeling of an iterative scheme for linear thermo-poroelasticity}

\author[FB]{Francesco Ballarin}
\author[SL]{Sanghyun Lee}
\author[SY]{Son-Young Yi}

\affiliation[FB]{organization={Department of Mathematics and Physics, Università Cattolica del Sacro Cuore},
            city={25133 Brescia},
            country={Italy}}
\affiliation[SL]{organization={Department of Mathematics, Florida State University},
            city={Tallahassee},
            state={FL 32304},
            country={USA}}
\affiliation[SY]{organization={Department of Mathematical Sciences, The University of Texas at El Paso},
            city={El Paso},
            state={TX 79968},
            country={USA}}

\begin{abstract}
This paper explores an iterative approach to solve linear thermo-poroelasticity problems, with its application as a high-fidelity discretization utilizing finite elements during the training of projection-based reduced order models.
One of the main challenges in addressing coupled multi-physics problems is the complexity and
computational expenses involved. In this study, we introduce a decoupled iterative solution approach, integrated with reduced order modeling, aimed at augmenting the efficiency of the computational algorithm.
The iterative technique we employ builds upon the established fixed-stress splitting scheme that has been extensively investigated for Biot's poroelasticity.
By leveraging solutions derived from this coupled iterative scheme, the reduced order model employs an additional Galerkin projection onto a reduced basis space formed by a small number of modes obtained through proper orthogonal decomposition. The effectiveness of the proposed algorithm is demonstrated through numerical experiments, showcasing its computational prowess.
\end{abstract}


\end{frontmatter}

\section{Introduction}
Thermo-hydro-mechanical (THM) processes refer to the coupled interactions between temperature, fluid flow, and mechanical deformation that occur in a wide range of natural and engineered systems. These processes are relevant in various fields, including environmental science, civil engineering, and material science, to name a few.
Therefore, the ability to understand and predict tightly coupled THM processes in natural and engineering systems has significant impacts, for example, on the environment, public safety, and the economy.

One widely employed mathematical model for describing THM processes is based on Biot's non-isothermal consolidation theory \cite{biot1956thermo}, known as the thermo-poroelasticity model. This model extends Biot's poroelasticity model \cite{biot1941general}, which characterizes the interaction between a deformable porous medium and fluid flow within it under isothermal conditions.
The governing system of partial differential equations (PDEs) of the thermo-poroelasticity model consists of a heat transfer equation, mass conservation equation, and momentum conservation equation. These equations are fully coupled through linear and nonlinear coupling terms.
For simplicity, this paper focuses on a fully-coupled, linearized thermo-poroelasticity model, {where the heat transfer is assumed to be diffusion-dominated and the coupling strength between the THM processes is rather weak.}

Due to the complex nature of the mathematical model for coupled multi-physics thermo-poroelasticity, it is highly challenging to develop an accurate and computationally efficient numerical method.
There are various approaches to handling the coupled nature of the problems.
First, one can solve the coupled system using a  monolithic numerical method, where all subproblems (PDE equations) are solved simultaneously in each time step.
In general, monolithic schemes are unconditionally stable but computationally expensive and may require some special linear or nonlinear solvers, often making them less practical.
Another approach is to use a sequential method, where the subproblems are separated and solved sequentially in each time step.
This approach is much more economical compared to the monolithic method. However, a sequential method may not converge to the right solution, {or effective stabilization is often required to make such methods unconditionally stable~\cite{kim2018unconditionally}.}
To overcome this limitation, one can resort to an iterative  approach, which is a staggered yet tightly coupled scheme.
In iterative coupling schemes, the subproblems are solved in a staggered way while ensuring full
convergence of the solution in each time step.

Iterative coupling schemes for poroelasticity have been extensively studied in recent years~\cite{Park83, FarhatEtAl91,  HuangZienk98, SettariWalters99,wheeler2020ipacs, Bause17, BauseElAl17, BauseElAl21}. In particular, one of the most successful approaches is based on the fixed-stress spitting method~\cite{Wheeler07,kim2011stability,mikelic2013convergence,YiBean17,lee2017iterative}.  In this fixed-stress scheme, one solves the flow problem while fixing the volumetric stress (mean stress), subsequently solving the mechanics problem using the updated pressure solution. The convergence and stability of this fixed-stress splitting method have been studied in the context of various numerical methods, including CG, DG, and mixed finite element methods.
In contrast to poroelasticity, the subject of iterative methods for thermo-poroelasticity has not been heavily studied.
Some recent papers on this topic {include~\cite{brun2020monolithic,kim2018unconditionally,both2019}}.

This paper first describes an iterative scheme for solving the thermal-poroelasticity model based on the classical continuous Galerkin finite element method, which is considered to be the high-fidelity solver or offline solver in the context of the reduced order modeling approach.
We especially utilize the fixed stress iterative method and derive the stability and convergence properties that are close to those of the monolithic approach.

Next, we utilize this fixed stress iterative scheme to train a projection-based reduced order model (ROM) as outlined in \cite{hesthaven2016certified,quarteroni2015reduced}. ROMs generally consist of two distinct phases: the offline (or training) stage, and the online (or evaluation) phase.
In the offline stage, a high-fidelity scheme is employed to generate a dataset containing solutions for the specific problem.
In our context, we utilize the fixed stress iterative method, and this dataset is compressed using proper orthogonal decomposition \cite{berkooz1993proper}, yielding a reduced basis space of significantly lower dimension. During the online phase, the ROM incorporates a Galerkin method that leverages the reduced basis space, in contrast to the original high-fidelity space. If the reduced basis space possesses dimensions substantially smaller than its high-fidelity counterpart, and if the problem operators can be suitably pre-assembled, then the ROM evaluation is anticipated to be significantly more computationally efficient than assessing the high-fidelity scheme.
The application of projection-based ROMs to THM systems is still limited, see in particular \cite{larion2020building,florez2022model}. Furthermore, the use of iterative schemes as high-fidelity solvers in ROMs  \cite{boon2023reduced,NoninoBallarinRozzaMaday2022,NoninoBallarinRozza2021,BallarinRozzaMaday2017,ShahGirfoglioQuintelaRozzaLengominBallarinBarral2021,zancanaro2022segregated,ngan2023reducedorder,kadeethum2022enhancing}
is not as common as monolithic ones, and none of the aforementioned references uses a scheme with the specific properties of the fixed-stress iterative scheme. Therefore, by applying a fixed-stress iterative ROM to a THM problem, this work advances the state-of-the-art and could pave the way for applications of the proposed ROM in geothermal flows.

The rest of the paper is organized as follows.
We describe the governing equations for the thermo-poroelasticity model in Section~\ref{sec: govern}.
We then introduce the high-fidelity discretization based on the fixed-stress iterative scheme in Section~\ref{sec: hf}. In there, we also prove that the solution of the iterative scheme converges to the  monolithic scheme, {assuming some conditions for certain physical parameters}.
In Section~\ref{sec: rom}, we introduce the projection-based ROM, which builds upon the fixed-stress iterative scheme, and show that its solution converges to the one of the ROM trained on the monolithic scheme.
Finally, we present several numerical experiments in Section~\ref{sec: num_exp} to validate and present the capabilities of the proposed algorithm.

\section{Governing Equations}\label{sec: govern}
~Let $\Omega$ be a bounded, connected, and Lipschitz domain in $\mathbb{R}^d$, $d=2, 3$, with the boundary $\partial \Omega$ and let $\mathbb{I} = (0, T]$ with $T >0$.
Then, let  $\bu: \Omega \rightarrow \mathbb{R}^d$ be the vector-valued displacement of the solid, $p: \Omega \rightarrow \mathbb{R}$   the  scalar-valued fluid pressure, and $\theta: \Omega \rightarrow \mathbb{R}$ the  scalar-valued temperature. Then, governing equations for thermo-poroelasticity are derived by coupling momentum balance for mechanics based on linear elasticity, mass balance for the pressure, and energy balance for the temperature as follows:
\begin{subequations}\label{sys: thm}
\begin{alignat}{2}
& - \divv (  \bsigma(\bu)  -  \alpha p \bI - 3 \alpha_T K_{dr} \theta \bI)= \bff  &&\quad  \text{ in } \Omega \times \mathbb{I}, \label{thm-m}\\
& \pp{}{t}(\czero p + \alpha \divv \bu - 3\alpha_m \theta )
- \divv (\bK\nabla p) = g && \quad \text{ in } \Omega \times \mathbb{I}, \label{thm-h} \\
& \dfrac{\partial}{\partial t}
\left(\cd \theta + 3 \alpha_T K_{dr}\thetaz \nabla \cdot \bu - 3 \alpha_m \thetaz p  \right)
- \nabla \cdot  (\bD \nabla \theta)
= \eta
&& \quad \text{ in } \Omega \times \mathbb{I},
\label{thm-t}
\end{alignat}
\end{subequations}
In the momentum balance equation \eqref{thm-m},
$\bsigma(\bu)$ is the standard stress tensor from linear elasticity. It satisfies the constitutive equation $ \bsigma(\bu) := 2\mu  \bepsilon(\bu) +  \lambda ( \divv \bu) \bI$, where $\bepsilon(\bu) := \frac{1}{2} [ \nabla \bu + (\nabla \bu)^T ]$ is the strain tensor, $\bI$ is the $d\times d$ identity tensor, and $\mu, \lambda$ are  the Lam\'{e} constants.  The Lam\'{e} constants are assumed to be  in the range
$\mu \in [\mu_0, \mu_1]$ and
$\lambda \in [0, \infty)$ for some $0 < \mu_0 < \mu_1 < \infty$.
Also, $\bff$ is the body force, $\alpha$ is the {\it Biot-Willis} constant, $K_{dr}:={(d\lambda+2\mu)}/d$ is the drained isothermal bulk modulus, and $\alpha_T$ is the volumetric skeleton thermal dilation coefficient. The total stress tensor is given by $\bsigt(\bu,p,\theta)= \bsigma(\bu) - \alpha p \bI - 3 \alpha_T K_{dr} \theta \bI$.

The second equation \eqref{thm-h} is the mass balance equation for the fluid, assuming the Darcy law for the volumetric fluid flux: $\bq = -\bK \nabla p$.
We ignore the gravity effect and set the fluid viscosity to be one here for a simple presentation of the numerical method. However, including the gravity term and the fluid viscosity in the numerical formulation is straightforward.
Here, $\bK \in \mathbb{R}^{d \times d} $ is the permeability tensor, which is symmetric and uniformly positive-definite and satisfies the following assumption: there exist positive constants $k_{\text{min}}, k_{\text{max}}$  such that for any $\bx \in \Om$,
\begin{equation}\label{Kelliptic_1}
k_{\text{min}} \bxi^T \bxi \leq \bxi^T \bK(\bx) \bxi \leq k_{\text{max}} \bxi^T \bxi, \quad \forall \bxi \in \mathbb{R}^d.
\end{equation}
In addition, $c_0 = 1/M$, where $M$ is Biot's modulus, $3\alpha_m$ is the thermal dilation coefficient, and $g$ is the volumetric fluid source/sink term.

Finally, the energy balance equation, or the heat transfer equation  \eqref{thm-t}, is obtained by assuming local thermal equilibrium between solid and fluid in pores. Therefore, this energy balance equation is expressed in terms of a single temperature variable $\theta$ with the effective total heat conductivity $C_d$.
Additionally, $\eta$ is the volumetric heat source/sink term, and $\theta_0$ is a reference temperature, which is assumed to be nonzero.
In addition, the bulk thermal conductivity tensor $\bD$ is symmetric and uniformly positive-definite and assumed to satisfy the following:
there exist positive constants $d_{\text{min}}, d_{\text{max}}$ such that for any $\bx \in \Om$,
\begin{equation}\label{Kelliptic_2}
d_{\text{min}} \bxi^T \bxi \leq \bxi^T \bD(\bx) \bxi \leq d_{\text{max}} \bxi^T \bxi, \quad \forall \bxi \in \mathbb{R}^d.
\end{equation}
{Note that we adopt a linearized heat transfer equation in this paper, contrary to using a nonlinear equation in other work (e.g. \cite{kim2018unconditionally}). This linearization of the model by introducing the reference temperature and dropping the nonlinear convection term can be justified if we assume small magnitudes of $\alpha_T$ and $\alpha_m$ and that diffusion is a dominant process for heat transfer. Similar linearized models were used in \cite{zimmerman2000coupling,cacace2017flexible}.}

To complete the system of governing equations \eqref{sys: thm}, we have to provide initial conditions and boundary conditions.
As our main focus here is on the iterative coupling scheme, we will consider the pure homogeneous Dirichlet boundary conditions for all three variables,
whereas the initial conditions are given as
\begin{equation}\label{eqn: ic}
\bu( \cdot, 0) = \bu^0, \quad p(\cdot, 0) = p^0, \quad  \theta(\cdot, 0) = \theta^0 \quad \forall x \in \Omega.
\end{equation}

\section{High-fidelity Discretization}\label{sec: hf}
In this section, we present the discretization schemes employed during the reduced order model training in section \ref{sec: rom}, which will be referred to as \textit{high-fidelity schemes} therein.
We discuss first the monolithic scheme (section \ref{sec: coupled}), then introduce the proposed fixed-stress iterative scheme (section \ref{sec: fs}). In section \ref{sec: convergence}, we prove a convergence result for the proposed fixed-stress iterative scheme by utilizing
{continuous Galerkin finite element methods}.
\subsection{Monolithic high-fidelity (M-HF) method}\label{sec: coupled}
~In this section, we present a weak formulation and a fully discrete continuous Galerkin method for the model problem \eqref{sys: thm}.
The standard notation for the $L^2$- and Sobolev spaces and their associated inner products and norms will be used here. In particular, $\left(\cdot, \cdot\right)$ and $\normz{\cdot}$ denote the $L^2$ inner product and $L^2$-norm, respectively, and  $\left(\cdot, \cdot\right)_1$  and $\normo{\cdot}$ denote the (full) $H^1$-inner product and (full) $H^1$-norm, respectively.

To derive a weak formulation, we multiply \eqref{thm-m}, \eqref{thm-h}, and \eqref{thm-t} by $\bv \in \Honedz$, $ w \in \Honez$, and $ s \in \Honez$,  respectively, and integrate by parts, resulting in the following weak problem:
 At  every $t \in (0, T]$, find $(\bu, p, \theta) \in \Honedz \times \Honez  \times \Honez$  such that
\begin{subequations}\label{variational}
\begin{alignat}{2}
2\mu(\varepsilon(\bu), \varepsilon(\bv)) + \lambda ( \divv \bu, \divv \bv)  - \alpha(p, \divv \bv) - 3\alpha_T K_{dr} (\theta , \nabla \cdot \bv)  &= (\bff, \bv), \label{variational-elasticity} \\
\czero ( p_t, w )  + \alpha(\divv \bu_t, w ) - 3\alpha_m (\theta_t, w)
+  ( \bK \nabla p, \nabla w) & = (g, w),
\label{variational-flow}   \\
\cd ( \theta_t, s )  + 3 \alpha_T K_{dr} \thetaz( \nabla \cdot \bu_t, s )  - 3\alpha_m \thetaz (p_t, s)
+  (\bD \nabla \theta, \nabla s)  &= (\eta, s),
\label{variational-thermal}
\end{alignat}
\end{subequations}
 for any $(\bv, w, s) \in \Honedz \times \Honez \times \Honez$.

To present a fully-discrete method, first let $\dt > 0$ be a time step size and  $t^n = n \dt$, where $n = 0, \ldots, N$ with $t^N = T$.
For temporal discretization, we consider the backward Euler time-stepping scheme for simplicity. For spatial discretization, let us consider a shape-regular triangulation $\Th$ of $\Omega$.
Then, we let
$\Vh, \Wh, \Thetah$ be finite-dimensional subspaces of $\Honedz, \Honez, \Honez$, respectively, on the mesh $\Th$.
Then, the fully-discrete monolithic CG method for \eqref{variational} reads as follows:
Given $(\buh^n, \ph^n, \thetah^n) \in \Vh \times \Wh \times \Thetah$, find $(\buh^{n+1}, \ph^{n+1}, \thetah^{n+1})\in \Vh \times \Wh \times \Thetah$ such that
\begin{subequations}\label{sys: cg}
\begin{align}
&2\mu(\varepsilon(\buh^{n+1}), \varepsilon(\bv)) + \lambda ( \divv \buh^{n+1}, \divv \bv)   - \alpha(\ph^{n+1}, \divv \bv) - 3\alpha_T K_{dr} (\thetah^{n+1}, \nabla \cdot \bv) \notag \\
& \qquad \qquad = (\bff^{n+1}, \bv), \label{eqn: cg-elasticity} \\
&\czero \left( \frac{\ph^{n+1} - \ph^{n}}{\dt}, w \right)  + \alpha\left(\frac{\divv  \buh^{n+1} - \divv \buh^{n} }{\dt}, w \right) - 3\alpha_m \left( \frac{\thetah^{n+1}-\thetah^{n}}{\dt}, w \right)
 \notag \\
& \qquad \qquad  +  ( \bK \nabla \ph^{n+1}, \nabla w) = (g^{n+1}, w),
\label{eqn: cg-flow}   \\
&\cd \left( \frac{\thetah^{n+1}-{\thetah^{n}}}{\dt}, s \right)  + 3 \alpha_T K_{dr} \thetaz \left( \frac{\divv  \buh^{n+1} - \divv \buh^{n} }{\dt}, s \right)  - 3\alpha_m \thetaz \left(\frac{\ph^{n+1} - \ph^{n}}{\dt}, s \right)  \notag \\
& \qquad \qquad
+ (\bD \nabla \thetah^{n+1}, \nabla s)  = (\eta^{n+1}, s),
\label{eqn: cg-thermal}
\end{align}
\end{subequations}
for any $(\bv, w, s) \in \Vh \times \Wh \times \Thetah$.

The well-posedness and convergence of the CG method for the same model as \eqref{sys: thm} but with an extra nonlinear coupling term $-\bK \nabla p^n \cdot \nabla \theta^{n+1}$ in the temperature equation were studied in \cite{ZhangRui2022}. {In their paper, they employed the standard continuous piecewise $\mathbb{P}_k$ finite element for all $\Vh, \Wh, \Thetah$ for $k \ge 1$, where $\mathbb{P}_k$ is a space of all polynomials of degree at most $k$. This space, $\Vh \times \Wh \times \Thetah = \mathbb{P}_k\times \mathbb{P}_k \times \mathbb{P}_k$, does not satisfy the inf-sup stability condition. Despite the lack of the inf-sup stability condition, however, they proved that the resulting method is still well-posed and convergent with some constraints on the physical parameters. These constraints are the same as the ones for the convergence of our FS-HF iterative scheme, to be established in \eqref{eqn: assump1}. }


\subsection{Fixed-stress high-fidelity (FS-HF) iterative scheme}\label{sec: fs}
{Solving the monolithic system resulting from the M-HF method \eqref{sys: cg} is very expensive computationally.} Our goal is to develop an iteratively coupled method whose solutions converge to the solution of  the M-HF method. To this end, we will employ a version of the fixed-stress splitting scheme. {Similarly to \cite{kim2018unconditionally,both2019}, }the main idea behind the fixed-stress split is to fix (or lag) the total mean stress  while solving the flow and heat problems. Here, the total mean stress, $\bsigb$, is defined as $\bsigb = \mathrm{tr}(\bsigt)/d$, where $\mathrm{tr}(\cdot)$ is the trace operator on tensors.
By noting the relation between volumetric stress and strain, we can see that
\begin{equation}
   \bsigb=  K_{dr} \nabla \cdot \bu - \alpha p - 3 \alpha_T K_{dr} \theta.
\end{equation}
In light of this relationship, let $\bsigb_{h}^{n+1}$ be the approximate total mean stress at time $t^{n+1}$, defined by
\[
\bsigb_{h}^{n+1} =  K_{dr} \nabla \cdot \buh^{n+1} - \alpha \ph^{n+1} - 3 \alpha_T K_{dr} \thetah^{n+1}.
\]
Now, to  present our fixed-stress iterative scheme, we will denote by
$\phi^{n+1, i}$ the $i$-th iterate of $\phi^{n+1}$ for any variable $\phi = \buh, \ph, \thetah$.
To derive a fixed-stress iterative scheme, we want the condition
\[
\bsigb_{h}^{n+1, i+1} = \bsigb_h^{n+1, i}, \quad i = 0, 1, 2, \ldots
\]
or, equivalently,
\begin{subequations}
\begin{equation}\label{eqn: divu}
\nabla \cdot \buh^{n+1, i+1}   = \frac{\alpha}{K_{dr}} \left( \ph^{n+1, i+1} -\ph^{n+1, i}  \right)
+ 3\alpha_T \left(\theta_h^{n+1, i+1} - \theta_h^{n+1, i} \right)
+   \nabla \cdot \buh^{n+1, i},
\end{equation}
while solving the flow and heat equations. 
Note that the flow and heat equations can be solved simultaneously or sequentially, leading to different versions of the fixed-stress scheme. In this section, we consider a method where we solve these two problems sequentially, solving the flow equation first, then the heat problem, before solving the mechanics problem.
Since we will solve the flow and heat equations sequentially, we will use a modified version of \eqref{eqn: divu} when solving them. Specifically, we will assume that the temperature and pressure variables, besides the total mean stress, are fixed while solving the flow and heat problems, respectively. This means that we use
\begin{equation}\label{eqn: divu2}
\nabla \cdot \buh^{n+1, i+1}   = \frac{\alpha}{K_{dr}} \left( \ph^{n+1, i+1} -\ph^{n+1, i}  \right)
+   \nabla \cdot \buh^{n+1, i},
\end{equation}
when solving the flow problem and
\begin{equation}\label{eqn: divu3}
\nabla \cdot \buh^{n+1, i+1}   =
 3\alpha_T \left(\theta_h^{n+1, i} - \theta_h^{n+1, i} \right)
+   \nabla \cdot \buh^{n+1, i},
\end{equation}
when solving the heat problem.
\end{subequations}

Given an initialization $\buh^{n+1, 0}$, $\ph^{n+1,0}$ and $\thetah^{n+1,0}$, the following is our FS-HF iterative scheme at each time step $n+1$ to generate infinite sequences
$\{ \buh^{n+1, i+1} \}_{i=0}^\infty, \{ \ph^{n+1, i+1} \}_{i=0}^\infty$, and $\{ \thetah^{n+1, i+1} \}_{i=0}^\infty$.
Note that $L > 0$  to appear in the algorithm is a stabilization coefficient.
If we drive a scheme using \eqref{eqn: divu2} and \eqref{eqn: divu}, then $L =  1$.
However, one can consider tuning this coefficient to achieve faster convergence{, see for instance \cite{Storvik2019}}. \\
\noindent
{\bf Step 1.} Given $(\buh^{n+1, i}, \ph^{n+1, i},\thetah^{n+1, i}) \in \Vh \times \Wh \times \Thetah$, solve the flow problem for $\ph^{n+1, i+1} \in \Wh$:
\begin{subequations}\label{sys: fs}
\begin{align}\label{eqn: fs-flow}
& \czero \left( \frac{\ph^{n+1, i+1} - \ph^{n}}{\dt}, w \right) +  ( \bK \nabla \ph^{n+1, i+1}, \nabla w)
+ L \frac{\alpha^2}{K_{dr}} \left(\frac{\ph^{n+1, i+1} -\ph^{n+1, i}}{\dt}, w \right)
 \notag \\
& \quad = (g^{n+1}, w)
- \alpha \left(\divv \left( \frac{\buh^{n+1, i} -  \buh^{n} }{\dt} \right), w \right)
+ 3\alpha_m  \left( \frac{\thetah^{n+1, i}-\thetah^{n}}{\dt}, w \right) \quad \forall w \in \Wh.
\end{align}
{\bf Step 2.} Given $(\buh^{n+1, i}, \ph^{n+1, i},\thetah^{n+1, i}) \in \Vh \times \Wh \times \Thetah$, solve the heat problem for $\thetah^{n+1, i+1} \in \Thetah$:
\begin{align}\label{eqn: fs-thermal}
&\cd   \left( \frac{\thetah^{n+1, i+1}-\thetah^{n}}{\dt}, s \right)
+  ( \bD \nabla \thetah^{n+1, i+1}, \nabla s)
+ 9 L \alphat^2 \Kdr \thetaz \left( \frac{\thetah^{n+1, i+1}- \thetah^{n+1, i}}{\dt}, s \right)
\notag \\
& \quad
  = (\eta^{n+1}, s)
  + 3 \alpham  \thetaz \left(\frac{\ph^{n+1, i} - \ph^{n}}{\dt}, s \right)
  - 3 \alphat \Kdr \thetaz  \left(\divv \left( \frac{\buh^{n+1, i} -  \buh^{n} }{\dt} \right), s \right)
  \quad \forall s \in \Thetah.
\end{align}
{\bf Step 3.} Given $(\buh^{n+1, i}, \ph^{n+1, i+1},\thetah^{n+1, i+1}) \in \Vh \times \Wh \times \Thetah$, solve the mechanics problem for $\buh^{n+1, i+1} \in \Vh$:
\begin{align}\label{eqn: fs-elasticity}
&2\mu(\varepsilon(\buh^{n+1, i+1}), \varepsilon(\bv)) + \lambda ( \divv \buh^{n+1, i+1}, \divv \bv)
 \notag \\
&  = (\bff^{n+1}, \bv) + \alpha(\ph^{n+1, i+1}, \divv \bv)  + 3\alpha_T K_{dr} (\thetah^{n+1, i+1}, \nabla \cdot \bv) \quad \forall \bv \in \Vh.
\end{align}
\end{subequations}

Steps 1, 2, and 3 are repeated by increasing $i$ to $i+1$ until an appropriate stopping criterion is satisfied. The specific choice of the stopping criterion depending on a given problem, as well as the procedure to initialize $\buh^{n+1,0}$, $\ph^{n+1,0}$, and $\thetah^{n+1,0}$, will be discussed in Section \ref{sec: num_exp}.

\subsection{Convergence analysis of the FS-HF scheme}\label{sec: convergence}

In this section, we will prove that the solution of the FS-HF scheme, \eqref{sys: fs}, converges to that of the M-HF scheme, \eqref{sys: cg}.
To facilitate the analysis, we introduce some notations for the errors in the iterates.
Let $\eu^{n+1, i} = \buh^{n+1, i} - \buh^{n+1}$ be the error in the $i$-th iterate of $\buh^{n+1}$ for $ 1 \leq n \leq N$ and $i \ge 1 $.
Likewise, $\ep^{n+1, i} = \ph^{n+1, i} - \ph^{n+1}$ and $\etheta^{n+1, i} = \thetah^{n+1, i} - \thetah^{n+1}$ are the errors in the  $i$-th iterates of $\ph^{n+1}$ an $\thetah^{n+1}$, respectively.
To derive a system of error equations, subtract \eqref{eqn: cg-elasticity}, \eqref{eqn: cg-flow}, and \eqref{eqn: cg-thermal} from  \eqref{eqn: fs-elasticity}, \eqref{eqn: fs-flow}, and \eqref{eqn: fs-thermal}, respectively, and do some algebraic manipulations.
We note that the resulting system involves quantities only at the current time step $n+1$. Therefore, in the following, we will use superscripts indicating only the iteration numbers. Then, the system of error equations reads as follows:
\begin{subequations}\label{sys: err}
\begin{align}
&2\mu(\varepsilon(\eu^{i+1}), \varepsilon(\bv)) + \lambda ( \divv \eu^{i+1}, \divv \bv)
 = \alpha(\ep^{i+1}, \divv \bv)  + 3\alpha_T K_{dr} (\etheta^{i+1}, \nabla \cdot \bv), \label{eqn: err-elasticity} \\
& \czero \left( \frac{\ep^{i+1}}{\dt}, w \right) +  ( \bK \nabla \ep^{i+1}, \nabla w)
+ L \frac{\alpha^2}{K_{dr}} \left(\frac{\ep^{i+1} -\ep^{i}}{\dt}, w \right)
 \notag \\
& \qquad =
 - \alpha \left(\frac{\divv  \eu^{i}}{\dt} , w \right)
+ 3\alpha_m  \left( \frac{\etheta^{i}}{\dt}, w \right), \label{eqn: err-flow} \\
&\cd   \left( \frac{\etheta^{i+1}}{\dt}, s \right)
+  ( \bD \nabla \etheta^{i+1}, \nabla s)
+ 9 L \alphat^2 \Kdr \thetaz \left( \frac{\etheta^{i+1}- \etheta^{i}}{\dt}, s \right)
\notag \\
& \qquad
  =
   3 \alpham  \thetaz \left(\frac{\ep^{i}}{\dt}, s \right)
  - 3 \alphat \Kdr \thetaz  \left(  \frac{\divv \eu^{i} }{\dt}, s \right), \label{eqn: err-thermal}
\end{align}
for any $(\bv, w, s) \in \Vh \times \Wh \times \Thetah$.
\end{subequations}
For the subsequent analysis, we provide two preliminary results here.
\begin{lemma}\label{lem: coer}
The following coercivity condition is satisfied for any $\bv \in \Vh$.
\begin{equation}\label{eqn: coer}
\Kdr \normz{\divv \bv}^2 \leq 2 \mu \normz{\varepsilon(\bv)}^2 + \lambda \normz{\divv \bv}^2.
\end{equation}
\end{lemma}
\textbf{\textit{Proof.}}
It is easy to check that $\normz{\divv \bv}^2 \leq d \normz{\varepsilon(\bv)}^2$ using Young's inequality. Then, \eqref{eqn: coer} follows from the definition of $\Kdr$. \qed
\begin{lemma}\label{lem: eu}
The following inequality holds:
\begin{align}\label{eqn: eu}
& 2\mu \normz{\varepsilon(\eu^{i+1} - \eu^{i})}^2 + \lambda \normz{\divv (\eu^{i+1}- \eu^{i}) }^2 \notag \\
& \quad \leq {2} \left(  \frac{\alpha^2}{K_{dr}} \normz{\ep^{i+1} - \ep^i}^2
+  9 \alphat^2 \Kdr  \normz{\etheta^{i+1} - \etheta^i}^2 \right).
\end{align}
\end{lemma}
\textbf{\textit{Proof.}}
First, subtract the mechanics equation \eqref{eqn: err-elasticity}  at the iteration $i$ from  the iteration $i+1$:
\begin{align}\label{eqn: err-elas1}
&2\mu(\varepsilon(\eu^{i+1}- \eu^{i}), \varepsilon(\bv)) + \lambda ( \divv (\eu^{i+1} - \eu^{i}), \divv \bv)\notag \\
&\quad  = \alpha(\ep^{i+1}- \ep^{i}, \divv \bv)
 + 3\alpha_T K_{dr} (\etheta^{i+1}  - \etheta^{i}, \nabla \cdot \bv).
 \end{align}
Then, take
$\bv = \eu^{i+1} - \eu^{i}$ in \eqref{eqn: err-elas1} and use the coercivity condition \eqref{eqn: coer} on the left-hand side and the Cauchy-Schwarz inequality on the right-hand side to obtain
\begin{align}\label{eqn: err-elas2}
& \Kdr \normz{\divv (\eu^{i+1}- \eu^{i}) }^2 \notag \\
& \leq 2\mu \normz{\varepsilon(\eu^{i+1} - \eu^{i})}^2 + \lambda \normz{\divv (\eu^{i+1}- \eu^{i}) }^2 \notag \\
& \leq \left( \alpha \normz{\ep^{i+1}- \ep^{i}}+ 3\alpha_T K_{dr} \normz{\etheta^{i+1} - \etheta^{i}}\right) \normz{\divv (\eu^{i+1} - \eu^{i})},
\end{align}
which leads to
\begin{equation}\label{eqn: err-elas3}
 \normz{\divv (\eu^{i+1}- \eu^{i}) } \leq \frac{\alpha}{\Kdr} \normz{\ep^{i+1} - \ep^{i}}
 + 3\alpha_T  \normz{\etheta^{i+1}- \etheta^{i}}.
\end{equation}
Using \eqref{eqn: err-elas3} back in \eqref{eqn: err-elas2}, we now get
\begin{align*}
& 2\mu \normz{\varepsilon(\eu^{i+1} - \eu^{i})}^2 + \lambda \normz{\divv (\eu^{i+1}- \eu^{i}) }^2 \notag \\
& \leq \left( \alpha \normz{\ep^{i+1}- \ep^{i}}+ 3\alpha_T K_{dr} \normz{\etheta^{i+1} - \etheta^{i}}\right)
\left(\frac{\alpha}{\Kdr} \normz{\ep^{i+1} - \ep^{i}}
 + 3\alpha_T  \normz{\etheta^{i+1}- \etheta^{i}} \right),
\end{align*}
from which we can get the desired result \eqref{eqn: eu} after multiplying out the terms and using Young's inequality.
\qed

\begin{theorem} \label{thm:main}
Assume that
\begin{equation}\label{eqn: assump1}
\czero > 3 \alpham, \quad \cd > 3 \alpham \thetaz,
\end{equation} and
\begin{equation}\label{eqn: assump2}
 L \ge {2} \delta \quad \text{for some  } \delta \ge \half.
\end{equation}
For a fixed $n + 1$, assume that the FS-HF solution at time step $n + 1$ is initialized to the M-HF solution at time step $n$, i.e.
\[
\buh^{n+1, 0} = \buh^{n}, \quad \ph^{n+1, 0} = \ph^{n}, \quad \thetah^{n+1, 0} = \thetah^{n}.
\]
Then, the FS-HF solution defined in \eqref{sys: fs} converges to the solution of the M-HF method \eqref{sys: cg}:
\[
\normo{\buh^{n+1, i} - \buh^{n+1}} \to 0, \quad
\normz{\ph^{n+1, i} - \ph^{n+1}} \to 0,\quad
\normz{\thetah^{n+1, i} - \thetah^{n+1}} \to 0
\]
as $i \to \infty$.
\end{theorem}
\textbf{\textit{Proof.}}
We first recall that
\begin{equation}\label{eqn: i+1-i}
(e^{i+1}-e^i, e^{i+1}) = \frac{1}{2} \left( \normz{e^{i+1} - e^i}^2 + \normz{e^{i+1}}^{2} - \normz{e^i}^{2} \right)
\end{equation}
holds for any $e^i, e^{i+1}$ in $L^2(\Omega)$. Now, take $\bv = \eu^{i+1}, w = \dt \ \ep^{i+1}, s = \dt \ \etheta^{i+1}/\thetaz$ in \eqref{sys: err}, add the three equations, and apply \eqref{eqn: i+1-i} to obtain
\begin{align*}
& 2\mu \normz{\varepsilon(\eu^{i+1})}^2 + \lambda \normz{\divv \eu^{i+1}}^2 + \czero\normz{\ep^{i+1}}^2
+ \frac{\cd}{\thetaz}\normz{\etheta^{i+1}}^2 \notag \\
& \quad + \dt \left(\normz{ \bK^\half \nabla \ep^{i+1}}^2 + \frac{1}{\thetaz} \normz{ \bD^\half \nabla \etheta^{i+1}}^2 \right)
+ \frac{L}{2} \frac{\alpha^2}{K_{dr}} \left( \normz{\ep^{i+1} - \ep^i}^2
+ \normz{\ep^{i+1}}^{2} - \normz{\ep^i}^{2} \right)
\notag \\
& \qquad
+ \frac{9}{2 } L \alphat^2 \Kdr \left( \normz{\etheta^{i+1} - \etheta^i}^2 + \normz{\etheta^{i+1}}^2 - \normz{\etheta^i}^2
\right) - \alpha (\ep^{i+1}, \divv (\eu^{i+1} - \eu^{i}))  \notag \\
& \qquad \quad
- 3 \alphat \Kdr (\etheta^{i+1}, \divv (\eu^{i+1} - \eu^{i}))
=   3 \alpham ( \etheta^{i}, \ep^{i+1}) + 3 \alpham (\ep^{i}, \etheta^{i+1}).
\end{align*}
The last two terms on the left-hand side of  the above equation can be replaced by
\[
- 2\mu(\varepsilon(\eu^{i+1}), \varepsilon(\eu^{i+1} - \eu^{i})) - \lambda ( \divv \eu^{i+1}, \divv (\eu^{i+1} - \eu^{i})),
\]
which is obtained by taking $ \bv = -(\eu^{i+1} - \eu^{i})$ in \eqref{eqn: err-elasticity}.
Then, applying Young's inequality to the last two terms on the left-hand side and to the two terms on the right-hand side, we obtain
\begin{align*}
& 2\mu \normz{\varepsilon(\eu^{i+1})}^2 + \lambda \normz{\divv \eu^{i+1}}^2 + \czero\normz{\ep^{i+1}}^2
+ \frac{\cd}{\thetaz}\normz{\etheta^{i+1}}^2
+ \dt \left(  \normz{\bK^\half \nabla \ep^{i+1}}^2
+ \frac{1}{\thetaz}\normz{\bD^\half \nabla \etheta^{i+1}}^2 \right) \notag \\
& \ + \frac{L}{2} \frac{\alpha^2}{K_{dr}} \left( \normz{\ep^{i+1} - \ep^i}^2 + \normz{\ep^{i+1}}^{2} - \normz{\ep^i}^{2} \right)
+ \frac{9}{2} L \alphat^2 \Kdr \left( \normz{\etheta^{i+1} - \etheta^i}^2 + \normz{\etheta^{i+1}}^{2} - \normz{\etheta^i}^{2}
\right) \notag \\
& \ \  - \frac{\delta}{2}\left( 2\mu \normz{\varepsilon(\eu^{i+1} - \eu^{i})}^2
+ \lambda \normz{\divv (\eu^{i+1} - \eu^{i})}^2  \right) - \frac{1}{2 \delta}\left(
 2\mu \normz{\varepsilon(\eu^{i+1})}^2
+ \lambda \normz{\divv \eu^{i+1}}^2 \right) \notag \\
& \leq \frac{3 \alpham}{2}\left(\normz{\ep^{i+1} }^2 + \normz{\etheta^{i+1} }^2
+ \normz{\ep^{i} }^2 + \normz{\etheta^{i} }^2\right),
\end{align*}
where $\delta>0$ is a positive constant.
{After collecting terms with the superscript $i+1$ on the left-hand side and the terms with $i$ on the right-hand side, we bound the left-hand side below. Specifcally, we use the assumptions in \eqref{eqn: assump1} for the coefficients of $\normz{\ep^{i+1}}^2$ and $\normz{\etheta^{i+1}}^2$, and the uniform ellipticity \eqref{Kelliptic_1} and \eqref{Kelliptic_2}, and Poincar\'{e} inequality for  the flow and heat flux terms, and \eqref{eqn: eu} for the terms involving $\eu^{i+1} - \eu^{i}$}.
Then, we can arrive at the following inequality:
\begin{align*}
& \left( 1 - \frac{1}{2\delta}\right) \left( 2\mu \normz{\varepsilon(\eu^{i+1})}^2 + \lambda \normz{\divv \eu^{i+1}}^2 \right)  + \left( \frac{3\alpham}{2} + \frac{L}{2}\frac{\alpha^2}{K_{dr}} + \frac{k_{min}}{C_\Omega} \dt \right)  \normz{\ep^{i+1}}^2 \notag \\
& \quad +\left(\frac{3\alpham}{2} + \frac{L}{2} 9 \alphat^2 \Kdr + \frac{d_{min}}{C_\Omega} \frac{\dt}{\thetaz} \right) \normz{\etheta^{i+1}}^2
  + \left( \frac{L}{2} -{\delta}\right) \frac{\alpha^2}{K_{dr}} \left( \normz{\ep^{i+1} - \ep^i}^2  \right)
  \notag \\
& \qquad + \left( \frac{L}{2} -{\delta} \right)9 \alphat^2 \Kdr\left( \normz{\etheta^{i+1} - \etheta^i}^2
\right) \notag \\
& \leq
\left(\frac{3 \alpham}{2} + \frac{L}{2}\frac{\alpha^2}{K_{dr}}  \right) \normz{\ep^i}^2
+ \left(\frac{3 \alpham}{2} + \frac{9}{2} L \alphat^2 \Kdr \right)  \normz{\etheta^i}^2,
\end{align*}
where $C_\Omega$ is a constant from the Poincar\'{e} inequality.
Now, due to the assumptions in \eqref{eqn: assump2},
both  $\left( 1 - \frac{1}{2\delta}\right)$ and $\left( \frac{L}{2} -{\delta}\right) $ are  nonnegative, hence we have
\begin{align*}
&  \left( 3\alpham + L\frac{\alpha^2}{K_{dr}} + \frac{2 k_{min}}{C_\Omega}\dt \right)  \normz{\ep^{i+1}}^2
+\left(3\alpham + 9 L \alphat^2 \Kdr + \frac{2 d_{min}}{C_\Omega} \frac{\dt}{\thetaz} \right) \normz{\etheta^{i+1}}^2 \notag \\
& \leq
\left( 3\alpham + L \frac{\alpha^2}{K_{dr}}  \right) \normz{\ep^i}^2
+ \left( 3\alpham + 9 L  \alphat^2 \Kdr \right)  \normz{\etheta^i}^2.
\end{align*}
Now, let
\begin{equation}\label{eqn: tildeC}
\tilde{C} = \text{min}\left \{\frac{2 k_{min} \dt}{C_\Omega(3\alpham + L\frac{\alpha^2}{K_{dr}})}, \,
 \frac{2 d_{min} \dt}{C_\Omega \thetaz (3\alpham + 9 L \alphat^2 \Kdr )} \right \}.
\end{equation}
Then, we have
\begin{align}\label{eqn: err-flowheat}
& ( 1 + \tilde{C}) \left[ \left( 3\alpham + L\frac{\alpha^2}{K_{dr}} \right) \normz{\ep^{i+1}}^2+
(3\alpham +  9 L \alphat^2 \Kdr)\normz{\etheta^{i+1}}^2 \right ] \notag \\
& \leq \left( 3\alpham + L \frac{\alpha^2}{K_{dr}}  \right) \normz{\ep^i}^2
+ \left( 3\alpham +  9 L \alphat^2 \Kdr \right)  \normz{\etheta^i}^2.
\end{align}
If we define a norm  $\triplenorm{\cdot}$ in $\Wh \times \Thetah$ by
\[
\triplenorm{(w, s)}^2 = \left( 3\alpham + L \frac{\alpha^2}{K_{dr}}  \right) \normz{w}^2
+ \left( 3\alpham + 9 L \alphat^2 \Kdr \right)  \normz{s}^2 \quad \forall (w, s) \in \Wh \times \Thetah,
\]
then, \eqref{eqn: err-flowheat} can be rewritten as
\begin{equation}\label{eqn: contract}
\triplenorm{(\ep^{i+1},\etheta^{i+1})}^2
\leq \frac{1}{1 + \tilde{C}}\triplenorm{(\ep^{i},\etheta^{i})}^2.
\end{equation}
Therefore, $\triplenorm{(\ep^{i},\etheta^{i})} \to 0$  as $ i \to \infty$, which also implies that
$\normz{\ep^{i}} \to 0$
and $\normz{\etheta^{i}} \to 0$.
In order to prove the convergence of $\normo{\eu^{i}} \to 0$, first obtain
\begin{equation}\label{eqn: eu2}
 2\mu \normz{\varepsilon(\eu^{i+1})}^2 + \lambda \normz{\divv (\eu^{i+1}) }^2
\leq {2} \left(  \frac{\alpha^2}{K_{dr}} \normz{\ep^{i+1}}^2
+  9 \alphat^2 \Kdr  \normz{\etheta^{i+1}}^2 \right)
\end{equation}
in the same manner as in the proof of Lemma~\ref{lem: eu}, and
recall Korn's inequality, stating that there is a constant $C_{korn} >0$ such that
\[
C_{korn} \normo{\bv}^2 \leq \normz{\varepsilon(\bv)}^2
\quad \forall \bv \in \Honedz.
\]
Using these two results, we have
\[
\normo{\eu^{i+1}}^2 \leq  {\frac{1}{\mu C_{korn}}} \left(  \frac{\alpha^2}{K_{dr}} \normz{\ep^{i+1}}^2
+  9 \alphat^2 \Kdr  \normz{\etheta^{i+1}}^2 \right).
\]
As the quantity on the right-hand side goes to $0$ as $i \to \infty$, $\normo{\eu^{i}}$ also goes to $0$ as $i \to \infty$.

\qed

{
\begin{remark}
Similar constraints to the ones in (16) were needed to prove the well-posedness and convergence of numerical methods for thermo-poroelasticity in \cite{brun2019wellposed, brun2020monolithic, zhao2020locking}. On the other hand, Kim \cite{kim2018unconditionally} proved unconditional stability for an extended fixed-stress split for a nonlinear thermo-poroelasticity model.
\end{remark}
}


\section{Projection-Based Reduced Order Models}\label{sec: rom}

In this section, we discuss projection-based reduced order models (ROMs) based on compression of the time history of the solutions provided by one of the high-fidelity methods (either M-HF or FS-HF method) from Section~\ref{sec: hf}. Such compression comprises the offline stage (or training) of the ROM and will result in a reduced basis of small dimension $r$. The reduced basis will be obtained by means of the proper orthogonal decomposition (POD) algorithm, see Section~\ref{sec: POD}. The online stage (or evaluation) of the ROM will then consist of a Galerkin method applied on either the fully-coupled scheme (Section~\ref{sec: coupled ROM}) or the fixed-stress iterative scheme (Section~\ref{sec: fs ROM}). Convergence of the fixed-stress iterative ROM scheme to the solution of the fully-coupled ROM scheme will be proven in section \ref{sec: convergence ROM}.

\subsection{Proper orthogonal decomposition}\label{sec: POD}
Let $\phi = \buh, \ph, \thetah$ denote one of the components of the solution field, the displacement, pressure, and temperature, respectively.
Then, $\left\{\phi^{n}\right\}_{n = 0}^{N}$ denote the sequence obtained by collecting such a component at different times, $t^0, \hdots, t^N$. The proper orthogonal decomposition \cite{berkooz1993proper,quarteroni2015reduced,hesthaven2016certified} is a data analysis method that aims to obtain an equivalent representation of the sequence $\left\{\phi^{n}\right\}_{n = 0}^{N}$ as a linear combination of at most $N + 1$ orthogonal modes $\varphi_{0}, \hdots, \varphi_{N}$ and define information content indices $\nu_{0}, \hdots, \nu_{N}$ associated to each corresponding mode.  Assuming the latter to be sorted as $\nu_{0} \geq \nu_{1} \geq \hdots \geq \nu_{N}$, POD is typically used as a data \textit{compression} tool to obtain an approximate representation in the linear space $\Phi_r = \mathrm{span}(\varphi_{0}, \hdots, \varphi_{r - 1})$, truncating the combination to the first $r$ modes for some $r \leq N$ and neglecting the modes $\varphi_{r}, \hdots, \varphi_{N}$ associated to indices $\nu_{r}, \hdots, \nu_{N}$ with small information content.

In the context of ROMs, the primary interest of the application of POD is in the generation of the space $\Phi_r$, to be called  the \textit{reduced basis space}. To this end, we will focus on the so-called \textit{method of snapshots}, even though alternative presentations are possible, for instance by means of the singular value decomposition or the principal component analysis \cite{berkooz1993proper,quarteroni2015reduced,hesthaven2016certified}.
First of all, the method of snapshots requires solving the following eigenvalue problem
\begin{equation}
\label{eq:POD_eigenvalues}
\mat{C}^{\phi} \vec{v} = \nu \ \vec{v}, \text{  where } \left[\mat{C}^{\phi}\right]_{nm} = (\phi^n, \phi^m)_1 \text{ for every } n, m = 0, \hdots, N.
\end{equation}
Since the matrix $\mat{C}^{\phi} \in \mathbb{R}^{(N + 1) \times (N + 1)}$ is symmetric positive definite by construction, its eigenvalues $\nu_{0}, \hdots, \nu_{N}$ are real and positive and, without loss of generality, can be assumed to be sorted in decreasing order. Denoting by $\vec{v}_n \in \mathbb{R}^{N + 1}$ the eigenvector associated with the eigenvalue $\nu_n$, the $n$-th POD mode is then computed as
\begin{equation}\label{eq: POD mode}
\varphi_n = \frac{1}{\sqrt{\nu_n}} \sum_{\beta = 0}^{N} \left[\vec{v}_n\right]_{\beta} \phi^{\beta}.
\end{equation}
Indeed, the obtained POD modes are $(\cdot, \cdot)_1$-orthonormal. To see this, first note that $\vec{v}_n^T \vec{v}_m = \delta_{nm}$, where $\delta_{nm}$ is the Kronecker delta function. Then, we obtain
\begin{align}
\left(\varphi_n, \varphi_m\right)_1 &= \frac{1}{\sqrt{\nu_n \nu_m}} \left(\sum_{\beta = 0}^{N} \left[\vec{v}_n\right]_{\beta} \phi^{\beta}, \sum_{\gamma = 0}^{N} \left[\vec{v}_m\right]_{\gamma} \phi^{\gamma}\right)_1
= \frac{1}{\sqrt{\nu_n \nu_m}} \sum_{\beta = 0}^{N} \sum_{\gamma = 0}^{N} \left[\vec{v}_n\right]_{\beta} \left[\vec{v}_m\right]_{\gamma}
\left[\mat{C}^{\phi}\right]_{\beta, \gamma}\notag\\
&= \frac{1}{\sqrt{\nu_n \nu_m}} \vec{v}_n^T \mat{C}^{\phi} \vec{v}_m
= \frac{\nu_m}{\sqrt{\nu_n \nu_m}} \vec{v}_n^T \vec{v}_m = \delta_{nm}.\label{eq: POD orthonormal}
\end{align}

Finally, we notice that POD modes $\varphi_{0}, \hdots, \varphi_{N}$ satisfy homogeneous Dirichlet boundary conditions if the input sequence $\phi^{0}, \hdots, \phi^{N}$ is zero on $\partial\Omega$. This easily follows by taking the trace of \eqref{eq: POD mode} on $\partial\Omega$. Therefore, any element in $\Phi_r$ is by construction equal to zero on $\partial\Omega$.

\subsection{Monolithic reduced order model (M-ROM)}\label{sec: coupled ROM}
In this section, we describe a projection-based ROM, which builds upon the fully-coupled monolithic (M-HF) scheme introduced in Section~\ref{sec: coupled} as a high-fidelity discretization. During the offline stage, the high-fidelity discretization is queried to obtain three sequences, $\left\{\buh^{n}\right\}_{n = 0}^{N}$, $\left\{\ph^{n}\right\}_{n = 0}^{N}$, and $\left\{\thetah^{n}\right\}_{n = 0}^{N}$,  representing the time history of the displacement, pressure, and  temperature fields, respectively.
Upon selecting a  reduced basis size $r \leq N$ and applying POD to each sequence $\left\{\buh^{n}\right\}_{n = 0}^{N}$, $\left\{\ph^{n}\right\}_{n = 0}^{N}$, and $\left\{\thetah^{n}\right\}_{n = 0}^{N}$ as discussed in Section~\ref{sec: POD}, we obtain the following reduced basis spaces
$$
\Vr = \mathrm{span}\left\{\varphi_1^{\bu}, \hdots, \varphi_r^{\bu}\right\}, \,\,
\Wr  = \mathrm{span}\left\{\varphi_1^{p}, \hdots, \varphi_r^{p}\right\}, \,\,
\Thetar = \mathrm{span}\left\{\varphi_1^{\theta}, \hdots, \varphi_r^{\theta}\right\},
$$
respectively.

During the online stage, the ROM is a Galerkin method for problem \eqref{variational} on the reduced basis space $\Vr \times \Wr \times \Thetar$. Therefore, the fully-coupled monolithic ROM (M-ROM) scheme for \eqref{variational} reads as follows:
Given $(\bur^n, \pr^n, \thetar^n) \in \Vr \times \Wr \times \Thetar$, find $(\bur^{n+1}, \pr^{n+1}, \thetar^{n+1}) \in \Vr \times \Wr \times \Thetar$ such that
\begin{subequations}\label{sys: cg ROM}
\begin{align}
&2\mu(\varepsilon(\bur^{n+1}), \varepsilon(\bv)) + \lambda ( \divv \bur^{n+1}, \divv \bv)   - \alpha(\pr^{n+1}, \divv \bv) - 3\alpha_T K_{dr} (\thetar^{n+1}, \nabla \cdot \bv) \notag \\
& \qquad \qquad = (\bff^{n+1}, \bv), \label{eqn: cg-elasticity-ROM} \\
&\czero \left( \frac{\pr^{n+1} - \pr^{n}}{\dt}, w \right)  + \alpha\left(\frac{\divv  \bur^{n+1} - \divv \bur^{n} }{\dt}, w \right) - 3\alpha_m \left( \frac{\thetar^{n+1}-\thetar^{n}}{\dt}, w \right)
 \notag \\
& \qquad \qquad  +  ( \bK \nabla \pr^{n+1}, \nabla w) = (g^{n+1}, w),
\label{eqn: cg-flow-ROM}   \\
&\cd \left( \frac{\thetar^{n+1}-{\thetar^{n}}}{\dt}, s \right)  + 3 \alpha_T K_{dr} \thetaz \left( \frac{\divv  \bur^{n+1} - \divv \bur^{n} }{\dt}, s \right)  - 3\alpha_m \thetaz \left(\frac{\pr^{n+1} - \pr^{n}}{\dt}, s \right)  \notag \\
& \qquad \qquad
+ (\bD \nabla \thetar^{n+1}, \nabla s)  = (\eta^{n+1}, s),
\label{eqn: cg-thermal-ROM}
\end{align}
\end{subequations}
for any $(\bv, w, s) \in \Vr \times \Wr \times \Thetar$. Even though \eqref{sys: cg} and \eqref{sys: cg ROM} are both obtained by means of a Galerkin method, the fundamental difference is that the reduced basis space $\Vr \times \Wr \times \Thetar$ employed in \eqref{sys: cg ROM} is of small dimension $3r$, owing to the compression carried out by POD.

In order to highlight further differences between \eqref{sys: cg} and \eqref{sys: cg ROM}, let $\vec{u}_r^{n+1} \in \mathbb{R}^r$ be the vector whose components are the M-ROM degrees of freedom for $\bur^{n+1} \in \Vr$, i.e.
\begin{equation*}
\bur^{n+1} = \sum_{\beta = 0}^r [\vec{u}_r^{n+1}]_\beta \varphi_\beta^{\bu}.
\end{equation*}
Similarly, let $\vec{p}_r^{n+1} \in \mathbb{R}^r$ and $\vec{\theta}_r^{n+1} \in \mathbb{R}^r$ be the vector whose components are the M-ROM degrees of freedom for $\pr^{n+1} \in \Wr$ and $\thetar^{n+1} \in \Thetar$, respectively. The system \eqref{sys: cg ROM} can thus be written in the following block-matrix form
\begin{equation}\label{eq: block cg ROM}
\begin{bmatrix}
\mat{A}_r^{\bu \bu} & \mat{A}_r^{\bu p} & \mat{A}_r^{\bu \theta}\\
\mat{M}_r^{p \bu} & \mat{M}_r^{p p} + \mat{A}_r^{p p} & \mat{M}_r^{p \theta}\\
\mat{M}_r^{\theta \bu} & \mat{M}_r^{\theta p} & \mat{M}_r^{\theta \theta} + \mat{A}_r^{\theta \theta},\\
\end{bmatrix}
\begin{bmatrix}
\vec{u}_r^{n+1}\\
\vec{p}_r^{n+1}\\
\vec{\theta}_r^{n+1}
\end{bmatrix} =
\begin{bmatrix}
\textbf{0} & \textbf{0} & \textbf{0}\\
\mat{M}_r^{p \bu} & \mat{M}_r^{p p} & \mat{M}_r^{p \theta}\\
\mat{M}_r^{\theta \bu} & \mat{M}_r^{\theta p} & \mat{M}_r^{\theta \theta}\\
\end{bmatrix}
\begin{bmatrix}
\vec{u}_r^{n}\\
\vec{p}_r^{n}\\
\vec{\theta}_r^{n}
\end{bmatrix}
+ \begin{bmatrix}
\vec{f}_r^{n+1}\\
\vec{g}_r^{n+1}\\
\vec{\eta}_r^{n+1}
\end{bmatrix}
\end{equation}
where the expressions of matrices and vectors appearing in the block formulation are summarized in Table \ref{table:rom_operators}.
Indeed, matrices in Table \ref{table:rom_operators} can be assembled once and for all at the end of the offline stage since such computations involve integration over the finite element mesh $\Th$. Similarly, vectors in Table \ref{table:rom_operators} can be pre-assembled at the end of the offline stage for every $n = 1, \hdots, N$ and stored. During the online stage, iterating in time through \eqref{eq: block cg ROM} can finally be carried out at a vastly decreased computational cost since i) the assembly of the block system only requires loading $r \times r$ matrices and vectors of dimension $r$ in Table \ref{table:rom_operators} from storage without necessitating any operation involving the finite element $\Th$, and ii) the resulting linear system is of small size $3r \times 3r$.

\begin{table}
\centering
\small
\begin{tabular}{|ll|l|}
 \hline
 \multicolumn{2}{|c|}{\textbf{ROM matrices}} & \textbf{ROM vectors} \\
 \hline
\multicolumn{2}{|l|}{$[\mat{A}_r^{\bu \bu}]_{\beta\gamma} = 2\mu(\varepsilon(\varphi_\gamma^{\bu}), \varepsilon(\varphi_\beta^{\bu})) + \lambda ( \divv \varphi_\gamma^{\bu}, \divv \varphi_\beta^{\bu})$,}&
$[\vec{f}_r^{n+1}]_\beta = (\bff^{n+1}, \varphi_\beta^{\bu})$,\\
$[\mat{A}_r^{\theta \theta}]_{\beta\gamma} = (\bD \nabla \varphi_\gamma^{\theta} , \nabla \varphi_\beta^{\theta})$,&
$[\mat{A}_r^{p p}]_{\beta\gamma} = ( \bK \nabla \varphi_\gamma^{p}, \nabla \varphi_\beta^{p})$,&
$[\vec{g}_r^{n+1}]_\beta = (g^{n+1}, \varphi_\beta^{p})$,\\
$[\mat{A}_r^{\bu \theta}]_{\beta\gamma} = - 3\alpha_T K_{dr} (\varphi_\gamma^{\theta}, \nabla \cdot \varphi_\beta^{\bu})$,&
$[\mat{A}_r^{\bu p}]_{\beta\gamma} = - \alpha(\varphi_\gamma^{p}, \divv \varphi_\beta^{\bu})$,&
$[\vec{\eta}_r^{n+1}]_\beta = (\eta^{n+1}, \varphi_\beta^{\theta}).$\\
$[\mat{M}_r^{\theta \theta}]_{\beta\gamma} =
\cd \left( \frac{\varphi_\gamma^{\theta}}{\dt}, \varphi_\beta^{\theta} \right)$,&
$[\mat{M}_r^{p p}]_{\beta\gamma} = \czero \left( \frac{\varphi_\gamma^{p}}{\dt}, \varphi_\beta^{p} \right)$,&\\
$[\mat{M}_r^{\theta \bu}]_{\beta\gamma} =
 3 \alpha_T K_{dr} \thetaz \left( \frac{\divv  \varphi_\gamma^{\bu} }{\dt}, \varphi_\beta^{\theta} \right)$,&
$[\mat{M}_r^{\theta p}]_{\beta\gamma} =
- 3\alpha_m \thetaz \left(\frac{\varphi_\gamma^{p}}{\dt}, \varphi_\beta^{\theta} \right)$,&\\
$[\mat{M}_r^{p \bu}]_{\beta\gamma} = \alpha\left(\frac{\divv  \varphi_\gamma^{\bu}}{\dt}, \varphi_\beta^{p} \right)$,&
$[\mat{M}_r^{p \theta}]_{\beta\gamma} = - 3\alpha_m \left( \frac{\varphi_\gamma^{\theta}}{\dt}, \varphi_\beta^{p} \right)$,&\\
$[\mat{S}_r^{\theta \theta}]_{\beta\gamma} =
9 L \alphat^2 \Kdr \thetaz \left( \frac{\varphi_\gamma^{\theta}}{\dt}, \varphi_\beta^{\theta} \right)$,&
$[\mat{S}_r^{p p}]_{\beta\gamma} = L \frac{\alpha^2}{K_{dr}} \left(\frac{\varphi_\gamma^{p}}{\dt}, \varphi_\beta^{p} \right)$.&\\
 \hline
\end{tabular}
\caption{Definitions of ROM matrices and vectors in \eqref{eq: block cg ROM}. Here $\beta, \gamma = 1, \hdots, r$.}
\label{table:rom_operators}
\end{table}

\subsection{Fixed-stress reduced order model (FS-ROM)}\label{sec: fs ROM}
In this section, we further describe a projection-based ROM built upon the fixed-stress high-fidelity (FS-HF) iterative scheme described in Section~\ref{sec: fs}. During the offline stage, the FS-HF scheme is queried to obtain the time evolution $\left\{\buh^{n,\infty}\right\}_{n = 0}^{N}$, $\left\{\ph^{n, \infty}\right\}_{n = 0}^{N}$ and $\left\{\thetah^{n, \infty}\right\}_{n = 0}^{N}$ of displacement, pressure and temperature fields, respectively, where the superscript  ${n,\infty}$ denotes the converged solutions at time step $n$. We notice that in practice the number of iterations will not be $\infty$ upon defining suitable stopping criteria in Section~\ref{sec: num}.
Those sequences are then compressed by means of POD to obtain the reduced basis spaces $\Vr^{FS}, \Wr^{FS}$ and $\Thetar^{FS}$ upon proceeding as in Section~\ref{sec: coupled ROM}. We notice that the obtained reduced basis spaces $\Vr^{FS}, \Wr^{FS}$ and $\Thetar^{FS}$ are in principle different from $\Vr, \Wr$ and $\Thetar$ obtained in Section~\ref{sec: coupled ROM} since the former is obtained by applying POD to solutions of the FS-HF scheme, while the latter is the result of a compression of solutions computed by the M-HF scheme. For the sake of a simpler notation, however, we will drop the suffix $FS$ from the reduced basis spaces $\Vr^{FS}, \Wr^{FS}$, and $\Thetar^{FS}$, still understanding that reduced basis spaces employed in this section are generated from the FS-HF scheme, rather than the M-HF one.

During the online stage, given an initialization $(\bur^{n+1, 0}, \pr^{n+1,0}, \thetar^{n+1,0}) \in \Vr \times \Wr \times \Thetar$, the ROM uses a Galerkin projection of a fixed-stress iterative  scheme to generate infinite sequences
$\{ \bur^{n+1, i+1} \}_{i=0}^\infty \subset \Vr, \{ \pr^{n+1, i+1} \}_{i=0}^\infty \subset \Wr$, and $\{ \thetar^{n+1, i+1} \}_{i=0}^\infty \subset \Thetar$. \\

\noindent
{\bf Step 1-ROM.} Given $(\bur^{n+1, i}, \pr^{n+1, i},\thetar^{n+1, i}) \in \Vr \times \Wr \times \Thetar$, solve the flow problem for $\pr^{n+1, i+1} \in \Wr$:
\begin{subequations}\label{sys: fs ROM}
\begin{align}\label{eqn: fs-flow ROM}
& \czero \left( \frac{\pr^{n+1, i+1} - \pr^{n}}{\dt}, w \right) +  ( \bK \nabla \pr^{n+1, i+1}, \nabla w)
+ L \frac{\alpha^2}{K_{dr}} \left(\frac{\pr^{n+1, i+1} -\pr^{n+1, i}}{\dt}, w \right)
 \notag \\
& \quad = (g^{n+1}, w)
- \alpha \left(\divv \left( \frac{\bur^{n+1, i} -  \bur^{n} }{\dt} \right), w \right)
+ 3\alpha_m  \left( \frac{\thetar^{n+1, i}-\thetar^{n}}{\dt}, w \right) \quad \forall w \in \Wr.
\end{align}
{\bf Step 2-ROM.} Given $(\bur^{n+1, i}, \pr^{n+1, i},\thetar^{n+1, i}) \in \Vr \times \Wr \times \Thetar$, solve the heat problem for $\thetar^{n+1, i+1} \in \Thetar$:
\begin{align}\label{eqn: fs-thermal ROM}
&\cd   \left( \frac{\thetar^{n+1, i+1}-\thetar^{n}}{\dt}, s \right)
+  ( \bD \nabla \thetar^{n+1, i+1}, \nabla s)
+ 9 L \alphat^2 \Kdr \thetaz \left( \frac{\thetar^{n+1, i+1}- \thetar^{n+1, i}}{\dt}, s \right)
\notag \\
&  = (\eta^{n+1}, s)
  + 3 \alpham  \thetaz \left(\frac{\pr^{n+1, i} - \pr^{n}}{\dt}, s \right)
  - 3 \alphat \Kdr \thetaz  \left(\divv \left( \frac{\bur^{n+1, i} -  \bur^{n} }{\dt} \right), s \right) \quad \forall s \in \Thetar.
\end{align}
{\bf Step 3-ROM.} Given $(\bur^{n+1, i}, \pr^{n+1, i+1},\thetar^{n+1, i+1}) \in \Vr \times \Wr \times \Thetar$, solve the mechanics problem for $\bur^{n+1, i+1} \in \Vr$:
\begin{align}\label{eqn: fs-elasticity ROM}
&2\mu(\varepsilon(\bur^{n+1, i+1}), \varepsilon(\bv)) + \lambda ( \divv \bur^{n+1, i+1}, \divv \bv)
 \notag \\
& \quad = (\bff^{n+1}, \bv) + \alpha(\pr^{n+1, i+1}, \divv \bv)  + 3\alpha_T K_{dr} (\thetar^{n+1, i+1}, \nabla \cdot \bv) \quad \forall \bv \in \Vr.
\end{align}\\
\end{subequations}
Steps 1-ROM, 2-ROM, 3-ROM are repeated by increasing $i$ to $i +1$ until appropriate stopping criteria are satisfied. The specific choice of the stopping criteria, as well as the procedure to initialize $\bur^{n+1, 0}$, $\pr^{n+1,0}$ and $\thetar^{n+1,0}$, will be discussed in Section~\ref{sec: num}.

The FS-ROM scheme can be equivalently reformulated in matrix form as follows. It generates the sequences
$\{ \vec{\bu}_r^{n+1, i+1} \}_{i=0}^\infty \subset \mathbb{R}^r, \{ \vec{p}_r^{n+1, i+1} \}_{i=0}^\infty \subset \mathbb{R}^r$, and $\{ \vec{\theta}_r^{n+1, i+1} \}_{i=0}^\infty \subset \mathbb{R}^r$ containing their ROM degrees of freedom.\\

\noindent
{\bf Step i-ROM.} Given $(\vec{\bu}_r^{n+1, i}, \vec{p}_r^{n+1, i},\vec{\theta}_r^{n+1, i}) \in \mathbb{R}^r \times \mathbb{R}^r \times \mathbb{R}^r$, solve the flow problem for $\vec{p}_r^{n+1, i+1} \in \mathbb{R}^r$:
\begin{subequations}\label{sys: fs ROMeq}
\begin{align}\label{eqn: fs-flow ROMeq}
&\left(\mat{M}_r^{p p} + \mat{A}_r^{p p} + \mat{S}_r^{p p}\right) \vec{p}_r^{n+1, i + 1} = \vec{g}_r^{n+1}
- \mat{M}_r^{p p} \vec{p}_r^{n} - \mat{S}_r^{p p} \vec{p}_r^{n+1, i} \notag\\
&\quad  - \mat{M}_r^{p \bu} \left(\vec{u}_r^{n+1,i} - \vec{u}_r^{n}\right) - \mat{M}_r^{p \theta} \left(\vec{\theta}_r^{n+1, i} - \vec{\theta}_r^{n}\right).
\end{align}
{\bf Step ii-ROM.} Given $(\vec{\bu}_r^{n+1, i}, \vec{p}_r^{n+1, i},\vec{\theta}_r^{n+1, i}) \in \mathbb{R}^r \times \mathbb{R}^r \times \mathbb{R}^r$, solve the heat problem for $\vec{\theta}_r^{n+1, i+1} \in \mathbb{R}^r$:
\begin{align}\label{eqn: fs-thermal ROMeq}
&\left(\mat{M}_r^{\theta \theta} + \mat{A}_r^{\theta \theta} + \mat{S}_r^{\theta \theta}\right) \vec{\theta}_r^{n+1,i+1}
= \vec{\eta}_r^{n+1} - \mat{M}_r^{\theta \theta} \vec{\theta}_r^{n} - \mat{S}_r^{\theta \theta} \vec{\theta}_r^{n+1,i}\notag\\
&\quad- \mat{M}_r^{\theta \bu}\left(\vec{u}_r^{n+1, i} - \vec{u}_r^{n}\right)  - \mat{M}_r^{\theta p} \left(\vec{p}_r^{n+1, i} - \vec{p}_r^{n}\right).
\end{align}
{\bf Step iii-ROM.}  Given $(\vec{\bu}_r^{n+1, i}, \vec{p}_r^{n+1, i+1},\vec{\theta}_r^{n+1, i+1}) \in \mathbb{R}^r \times \mathbb{R}^r \times \mathbb{R}^r$, solve the mechanics problem for $\vec{\bu}_r^{n+1, i+1} \in \mathbb{R}^r$:
\begin{align}\label{eqn: fs-elasticity ROMeq}
&\mat{A}_r^{\bu \bu} \vec{u}_r^{n+1,i + 1} = \vec{f}_r^{n+1} - \mat{A}_r^{\bu p}\vec{p}_r^{n+1,i+1} - \mat{A}_r^{\bu \theta}\vec{\theta}_r^{n+1,i+1}.
\end{align}
\end{subequations}

To conclude, we note that \eqref{sys: fs ROMeq} requires solving three linear systems of dimension $r \times r$. Matrices and vectors appearing in \eqref{sys: fs ROMeq} can be precomputed as in Table \ref{table:rom_operators}, where the final row of the table contains two further matrices associated with the stabilization terms that were not required by M-ROM.

{\begin{remark}
Due to the evaluation of the ROM by means of a Galerkin projection, inf-sup stability does not necessarily get preserved by the ROM, even if the HF scheme was inf-sup stable. This issue is studied especially for saddle point problems arising from Stokes and Navier-Stokes equations, see \cite{RozzaVeroy2007,BallarinManzoniQuarteroniRozza2015}. The typical recipe to ensure inf-sup stability is to perform enrichment of the reduced order spaces, as proposed in \cite{RozzaVeroy2007} for the Stokes problem. However, this procedure can sometimes be avoided in practice, as discussed recently in  ~\cite{PacciariniRozza2014,AliBallarinRozza2020,StabileRozza2018}. The work in \cite{AliBallarinRozza2020} shows how performing a Galerkin projection of a stabilized Navier-Stokes problem, which incorporates finite element stabilization terms, can effectively replace the enrichment of the reduced order spaces. Furthermore, \cite{StabileRozza2018} shows that ROM based on iterative coupling techniques for Navier-Stokes problems do not need enrichment. Based on this experience, and the fact that (i) in \cite{BallarinManzoniQuarteroniRozza2015}, ROM solutions derived from spaces which lack inf-sup stability are very visibly incorrect, with relative errors up to $100\%$, (ii) as in \cite{StabileRozza2018}, our goal is to discuss a scheme with iterative coupling, and (iii) as in \cite{AliBallarinRozza2020}, the iterative scheme contains some stabilization, in the next numerical results we will not perform any enrichment. This will only empirically confirm that the resulting ROM is inf-sup stable, while further theoretical justifications are of relevant mathematical interest, but beyond the scope of this work.
\end{remark}}

\subsection{Convergence analysis of the FS-ROM}\label{sec: convergence ROM}
In this section, we retrace the convergence analysis carried out in Section~\ref{sec: convergence} in order to adapt it to the FS-ROM introduced in Section~\ref{sec: fs ROM}.

\begin{theorem} \label{thm:main ROM}
Assume that \eqref{eqn: assump1} and \eqref{eqn: assump2} hold. {Furthermore, assume that the M-ROM online stage employs a Galerkin projection onto the reduced basis spaces associated to FS-ROM}. For a fixed $n + 1$, assume that the FS-ROM solution at time step $n + 1$ is initialized to the M-ROM solution at time step $n$, i.e.
\[
\bur^{n+1, 0} = \bur^{n}, \quad \pr^{n+1, 0} = \pr^{n}, \quad \thetar^{n+1, 0} = \thetar^{n}.
\]

Then, the FS-ROM solution defined in \eqref{sys: fs ROM} converges to the solution of the M-ROM method \eqref{sys: cg ROM}:
\[
\normo{\bur^{n+1, i} - \bur^{n+1}} \to 0, \quad
\normz{\pr^{n+1, i} - \pr^{n+1}} \to 0,\quad
\normz{\thetar^{n+1, i} - \thetar^{n+1}} \to 0
\]
as $i \to \infty$.
\end{theorem}
\textbf{\textit{Proof.}}
Let $\eu^{i} = \bur^{n+1, i} - \bur^{n+1} \in \Vr$, $\ep^{i} = \pr^{n+1, i} - \pr^{n+1} \in \Wr$ and $\etheta^{i} = \thetar^{n+1, i} - \thetar^{n+1} \in \Vr$.
Upon subtracting \eqref{eqn: cg-elasticity-ROM}, \eqref{eqn: cg-flow-ROM}, and \eqref{eqn: cg-thermal-ROM} from  \eqref{eqn: fs-elasticity ROM}, \eqref{eqn: fs-flow ROM}, and \eqref{eqn: fs-thermal ROM} we obtain again \eqref{sys: err}, with the only exception that it must be now intended for any $(\bv, w, s) \in \Vr \times \Wr \times \Thetar$.
Also, Lemma~\ref{lem: coer} implies that the coercivity condition \eqref{eqn: coer} also holds for any $\bv \in \Vr$, since $\Vr \subset \Vh$. Even with the modified definitions of $\eu^{i}, \ep^{i}$ and $\etheta^{i}$ defined above, the proof of Lemma \ref{lem: eu} and the proof of the convergence in Theorem \ref{thm:main} follow in the same manner as the original proofs.
\qed

\begin{remark}
{We notice that we cannot provide bounds to
\begin{equation}\label{eq:fsROM_to_cg}
\normo{\bur^{n+1, i} - \buh^{n+1}}, \quad
\normz{\pr^{n+1, i} - \ph^{n+1}},\quad
\normz{\thetar^{n+1, i} - \thetah^{n+1}}
\end{equation}
with the same technique used in the proofs of Lemma \ref{lem: eu} and Theorem \ref{thm:main}.}
Indeed, if we define our errors as $\eu^{i} = \bur^{n+1, i} - \buh^{n+1} \in \Vh$, $\ep^{i} = \pr^{n+1, i} - \ph^{n+1} \in \Wh$ and $\etheta^{i} = \thetar^{n+1, i} - \thetah^{n+1} \in \Vh$, then we can obtain the same error equation, \eqref{sys: err}, as before but for any $(\bv, w, s) \in \Vr \times \Wr \times \Thetar$. However, since $(\eu^{i}, \ep^{i}, \etheta^{i}) \in \Vh \times \Wh \times \Thetar$, we cannot use these errors as test functions. As a result, the proofs of Lemma \ref{lem: eu} and Theorem \ref{thm:main} would not follow through.
Instead, we will empirically {quantify} \eqref{eq:fsROM_to_cg} by means of the numerical experiments in the next section.
\end{remark}

\section{Numerical Experiments}
\label{sec: num_exp}
In this section, we provide numerical experiments to validate and demonstrate the computational efficiency of the proposed algorithm. All the following computations are conducted by utilizing \texttt{FEniCSx}\footnote{\texttt{https://github.com/FEniCS/dolfinx}} for HF computations, and \texttt{RBniCSx} \footnote{\texttt{https://github.com/RBniCS/RBniCSx}} for ROM computations. \texttt{FEniCSx} and \texttt{RBniCSx} rely on \texttt{PETSc}\footnote{\texttt{https://petsc.org}} for sparse and dense linear algebra, respectively. \texttt{RBniCSx} further queries \texttt{LAPACK}\footnote{\texttt{https://www.netlib.org/lapack/}} for the solution of the dense eigenproblem \eqref{eq:POD_eigenvalues}.

\subsection{Example 1. Validation of ROM on test cases with smooth analytical solution}
\label{sec: num}
In this example, we validate our numerical methods against four different test cases with the given analytical solution, namely
\begin{subequations}\label{eq:analytical solution}
\begin{align}
&\bu(x, y; t) = \left[ \sin(\pi x t)\cos(\pi  y t), \cos(\pi x t) \sin(\pi y t) \right]^T x y (1 - x) (1 - y),\\
&p(x, y; t) = \cos(t + x -y)  x y (1 - x) (1 - y),\\
&\theta(x, y; t) = \sin(t + x - y)  x y (1 - x) (1 - y),
\end{align}
\end{subequations}
where $(x,y) \in \Omega = (0, 1)^2$ and $t \in \mathbb{I} = (0, T]$ for some $T > 0$ to be specified in each test case. We note that the factor $x y (1 - x) (1 - y)$ in each solution variable guarantees that the solution satisfies homogeneous Dirichlet boundary conditions, in agreement with the assumptions in the previous sections.
The right-hand side functions $\bff,g,$ and $\eta$ in \eqref{sys: thm} are obtained from the expressions in \eqref{eq:analytical solution}.
Furthermore, the physical parameters are chosen to be
\begin{align*}
&\lambda = 10^2, \  \mu = 10^2,  \ c_0 = 1, \  \alpha = 1,
\ K_{dr} = \mu +\lambda = 2 \cdot 10^2, \  \bK = 10^{-5},\\
&C_d = 1, \ \alphat=10^{-3}, \ \thetaz = 1, \ \alpham = 10^{-5}, \text{ and } \bD = 10^{-5}.
\end{align*}

To compute the initial conditions at $n=0$, the M-HF and FS-HF schemes employ the linear interpolation of \eqref{eq:analytical solution}, while M-ROM and FS-ROM employ the $L^2$-projection of \eqref{eq:analytical solution} on the reduced basis spaces. The stabilization parameter $L=1$ is employed in the FS-HF and FS-ROM schemes. Initialization of the iterations at time $n + 1$ is based on the assignment of the converged solution at time $n$.
Iterations for FS-HF will continue until the following conditions are satisfied
\begin{align*}
&\normo{\buh^{n+1, i+1} - \buh^{n+1, i}} \leq \epsilon \normo{\buh^{n+1, i+1}},\\
&\normo{\ph^{n+1, i+1} - \ph^{n+1, i}} \leq \epsilon \normo{\ph^{n+1, i+1}},\\
&\normo{\thetah^{n+1, i+1} - \thetah^{n+1, i}} \leq \epsilon \normo{\thetah^{n+1, i+1}},
\end{align*}
where $\epsilon = 10^{-10}$. Similarly, the FS-ROM iterations will continue until
\begin{align*}
&\normr{\vec{\bu}_r^{n+1, i+1} - \vec{\bu}_r^{n+1, i}} \leq \epsilon \normr{\vec{\bu}_r^{n+1, i+1}},\\
&\normr{\vec{p}_r^{n+1, i+1} - \vec{p}_r^{n+1, i}} \leq \epsilon \normr{\vec{p}_r^{n+1, i+1}},\\
&\normr{\vec{\theta}_r^{n+1, i+1} - \vec{\theta}_r^{n+1, i}} \leq \epsilon \normr{\vec{\theta}_r^{n+1, i+1}},
\end{align*}
{where $\normr{\cdot}$ denotes the Euclidean norm;}
we note that such conditions still represent a stopping criteria on the relative $\normo{\cdot}$-norm of the increment because, due to \eqref{eq: POD orthonormal}, the matrix representing the $(\cdot, \cdot)_1$ inner product on $\Vr, \Wr$ and $\Thetar$ is the identity.
Furthermore, we set the maximum iteration number as $20$ throughout this example.

\subsubsection{Example 1A. Validation of {HF} and {ROM} solvers by a mesh  refinement analysis}

\begin{figure}
    \centering
    \begin{subfigure}[b]{0.48\textwidth}
        \centering
        \includegraphics[width=\textwidth, trim={0 0.35cm 0 0},clip]{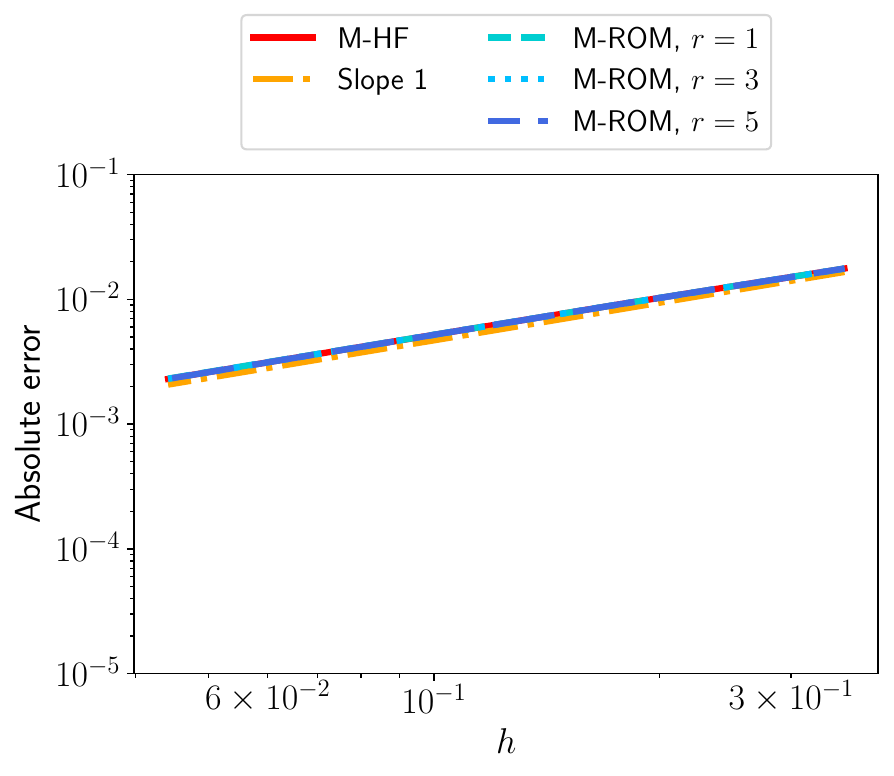}
        \caption{$\normo{\bu_\delta - \bu}$ for $\delta=h,r$.}
        \label{fig:test1_convergence_analysis_wrt_exact_u_H1_cg}
    \end{subfigure}
    \hfill
    \begin{subfigure}[b]{0.48\textwidth}
        \centering
        \includegraphics[width=\textwidth, trim={0 0.35cm 0 0},clip]{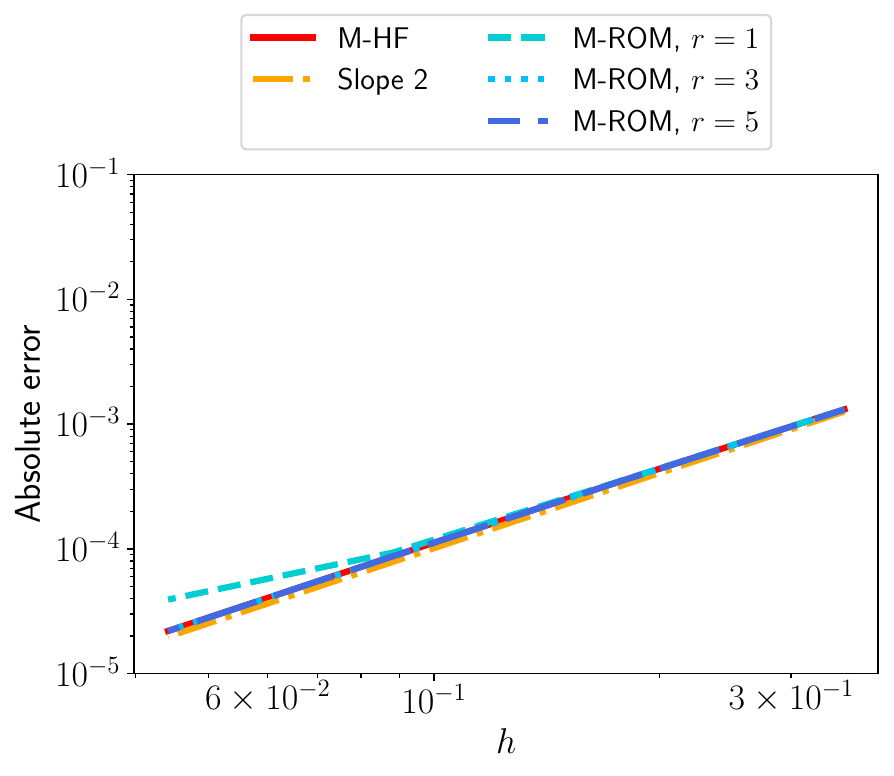}
        \caption{$\normz{\bu_\delta - \bu}$ for $\delta=h,r$.}
        \label{fig:test1_convergence_analysis_wrt_exact_u_L2_cg}
    \end{subfigure}\\
    \centering
    \begin{subfigure}[b]{0.48\textwidth}
        \centering
        \adjincludegraphics[width=\textwidth,Clip={0} {0.15cm} {0} {.2\height}]{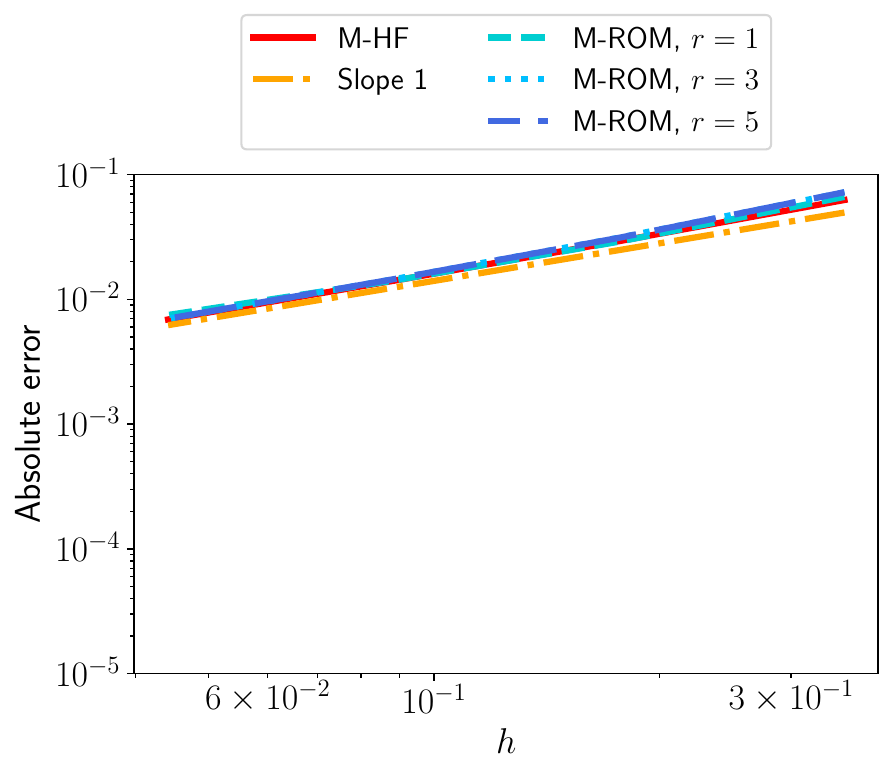}
        \caption{$\normo{p_\delta - p}$ for $\delta=h,r$.}
        \label{fig:test1_convergence_analysis_wrt_exact_p_H1_cg}
    \end{subfigure}
    \hfill
    \begin{subfigure}[b]{0.48\textwidth}
        \centering
        \adjincludegraphics[width=\textwidth,Clip={0} {0.15cm} {0} {.2\height}]{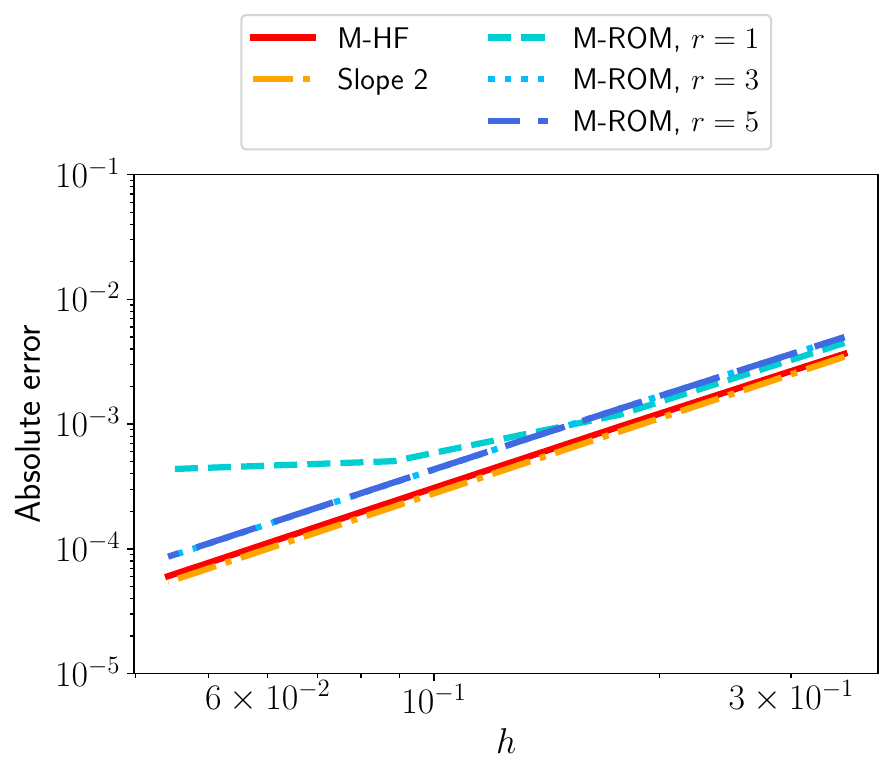}
        \caption{$\normz{p_\delta -p}$ for $\delta=h,r$.}
        \label{fig:test1_convergence_analysis_wrt_exact_p_L2_cg}
    \end{subfigure}\\
    \centering
    \begin{subfigure}[b]{0.48\textwidth}
        \centering
        \adjincludegraphics[width=\textwidth,Clip={0} {0.15cm} {0} {.2\height}]{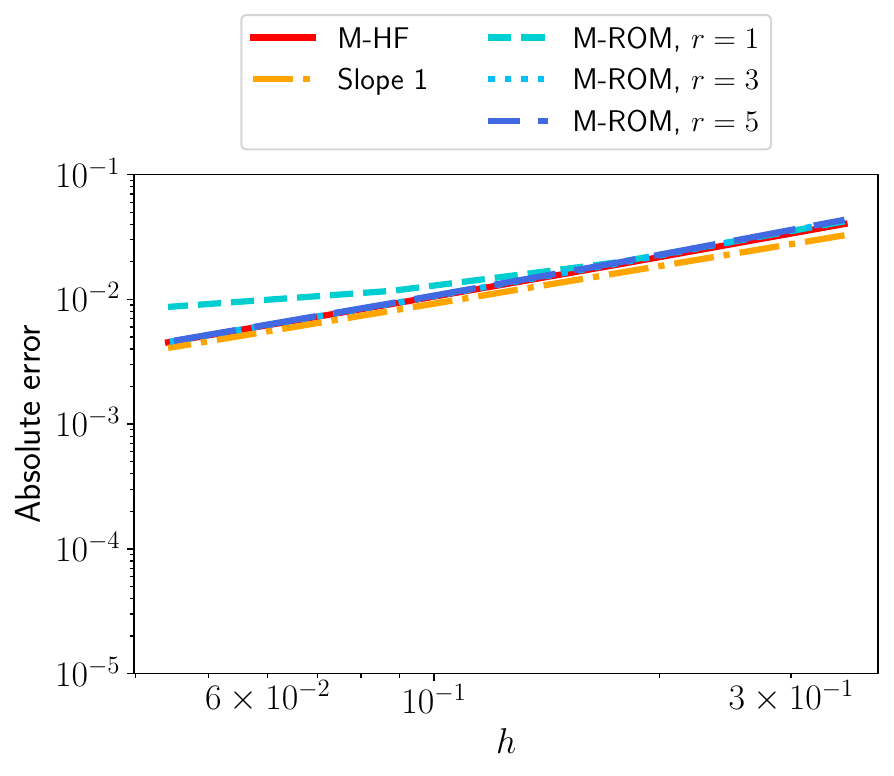}
        \caption{$\normo{\theta_\delta - \theta}$ for $\delta=h,r$.}
        \label{fig:test1_convergence_analysis_wrt_exact_theta_H1_cg}
    \end{subfigure}
    \hfill
    \begin{subfigure}[b]{0.48\textwidth}
        \centering
        \adjincludegraphics[width=\textwidth,Clip={0} {0.15cm} {0} {.2\height}]{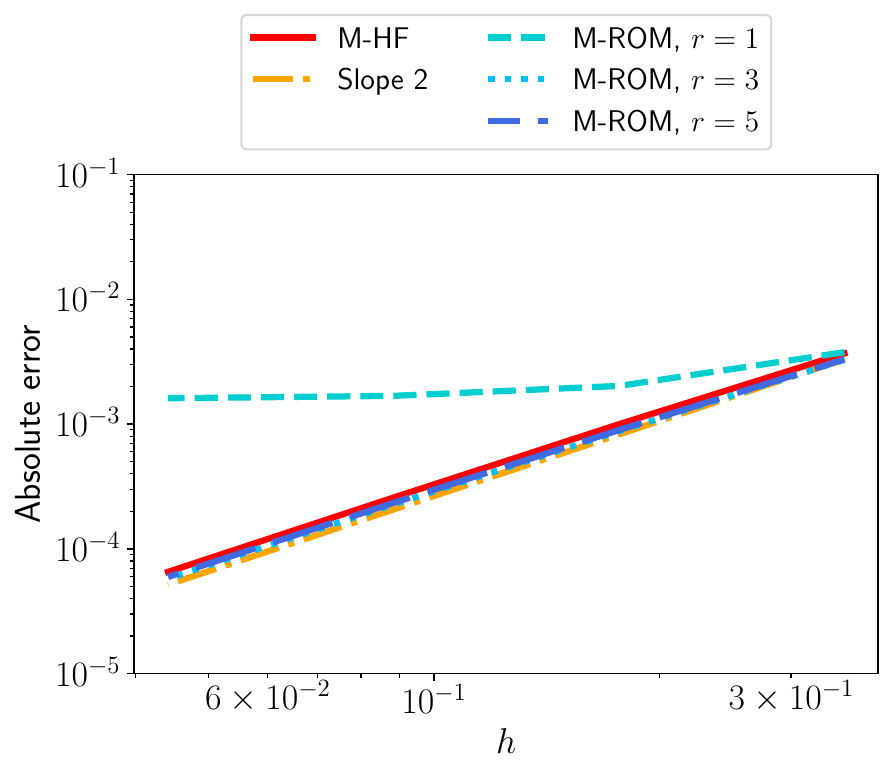}
        \caption{$\normz{\theta_\delta - \theta}$ for $\delta=h,r$.}
        \label{fig:test1_convergence_analysis_wrt_exact_theta_L2_cg}
    \end{subfigure}
    \caption{Example 1A: convergence of the errors for the monolithic schemes (M-HF and M-ROM).}
\label{fig:test1_convergence_analysis_wrt_exact_cg}
\end{figure}

\begin{figure}
    \centering
    \begin{subfigure}[b]{0.48\textwidth}
        \centering
        \includegraphics[width=\textwidth, trim={0 0.35cm 0 0},clip]{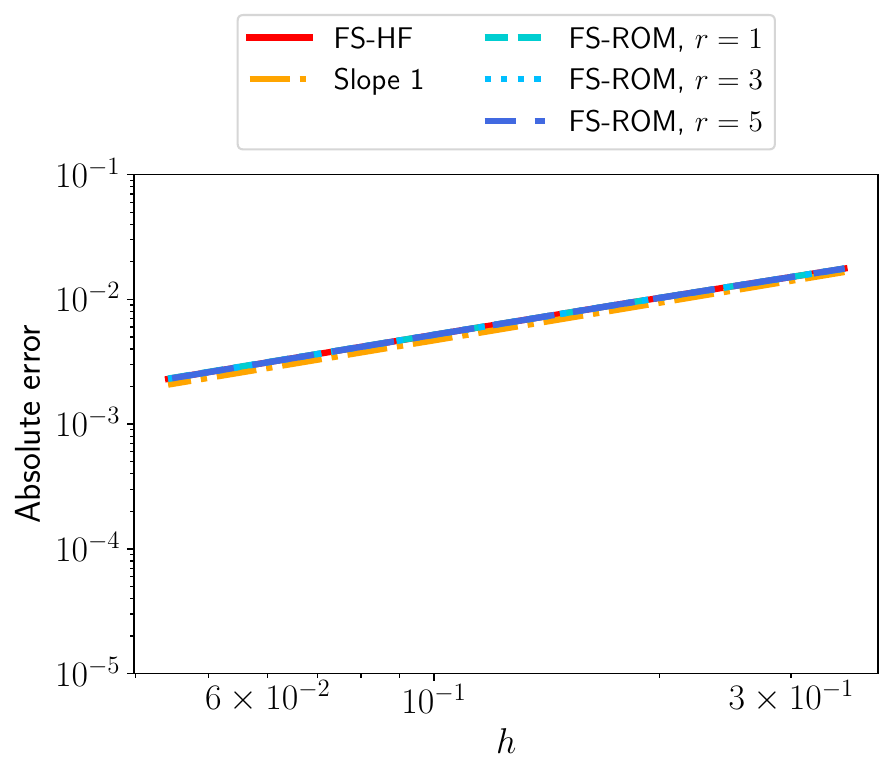}
        \caption{$\normo{\bu_\delta - \bu}$ for $\delta=h,r$.}
        \label{fig:test1_convergence_analysis_wrt_exact_u_H1_fs}
    \end{subfigure}
    \hfill
    \begin{subfigure}[b]{0.48\textwidth}
        \centering
        \includegraphics[width=\textwidth, trim={0 0.35cm 0 0},clip]{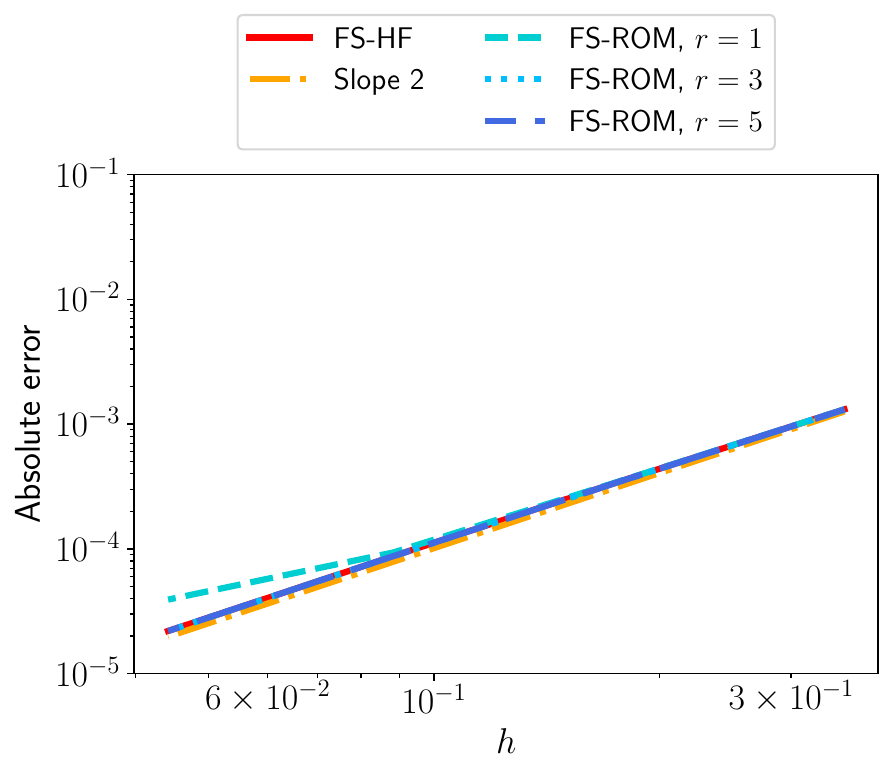}
        \caption{$\normz{\bu_\delta - \bu}$ for $\delta=h,r$.}
        \label{fig:test1_convergence_analysis_wrt_exact_u_L2_fs}
    \end{subfigure}\\
    \centering
    \begin{subfigure}[b]{0.48\textwidth}
        \centering
        \adjincludegraphics[width=\textwidth,Clip={0} {0.15cm} {0} {.2\height}]{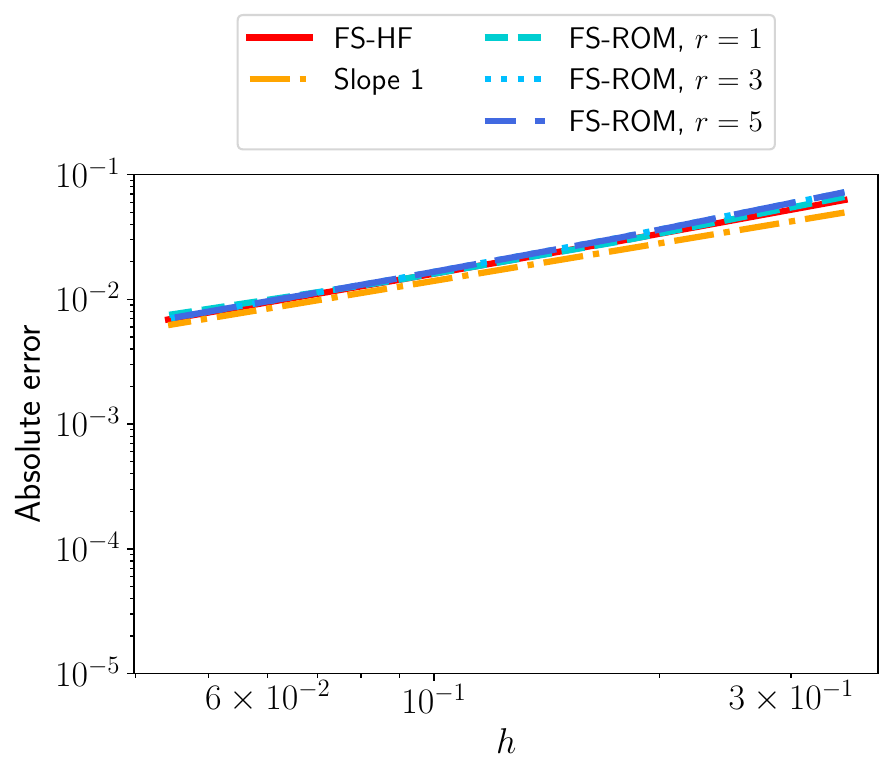}
        \caption{$\normo{p_\delta - p}$ for $\delta=h,r$.}
        \label{fig:test1_convergence_analysis_wrt_exact_p_H1_fs}
    \end{subfigure}
    \hfill
    \begin{subfigure}[b]{0.48\textwidth}
        \centering
        \adjincludegraphics[width=\textwidth,Clip={0} {0.15cm} {0} {.2\height}]{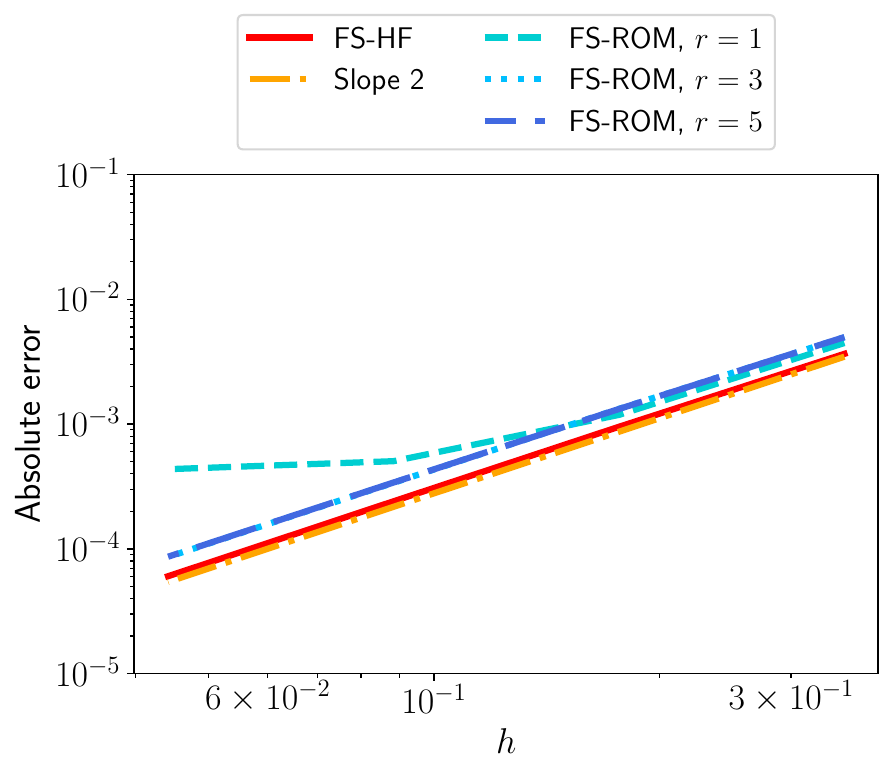}
        \caption{$\normz{p_\delta - p}$ for $\delta=h,r$.}
        \label{fig:test1_convergence_analysis_wrt_exact_p_L2_fs}
    \end{subfigure}\\
    \centering
    \begin{subfigure}[b]{0.48\textwidth}
        \centering
        \adjincludegraphics[width=\textwidth,Clip={0} {0.15cm} {0} {.2\height}]{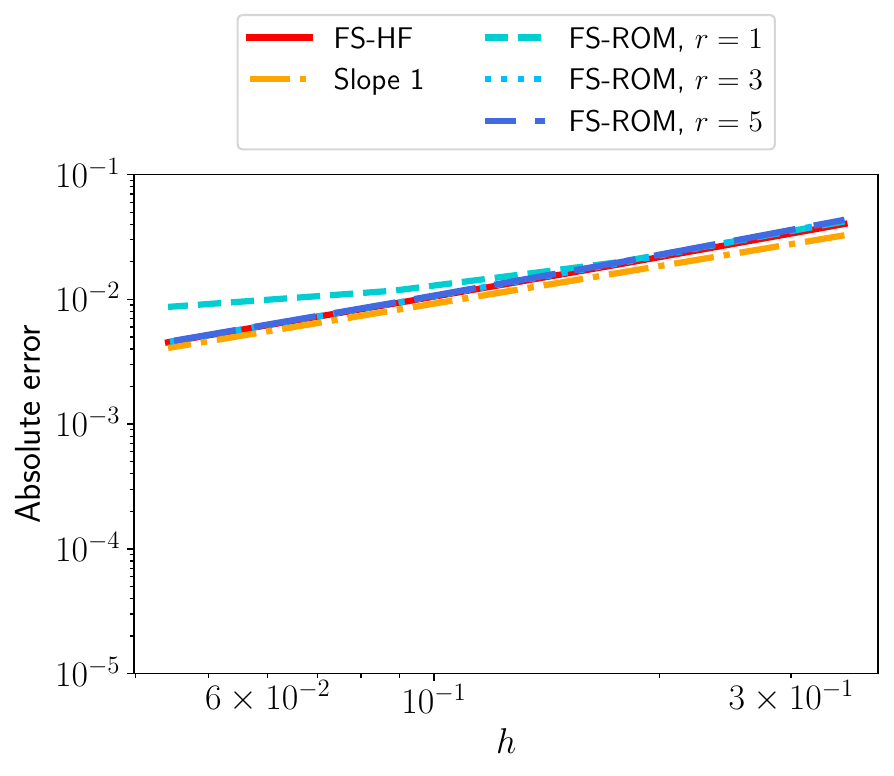}
        \caption{$\normo{\theta_\delta - \theta}$ for $\delta=h,r$.}
        \label{fig:test1_convergence_analysis_wrt_exact_theta_H1_fs}
    \end{subfigure}
    \hfill
    \begin{subfigure}[b]{0.48\textwidth}
        \centering
        \adjincludegraphics[width=\textwidth,Clip={0} {0.15cm} {0} {.2\height}]{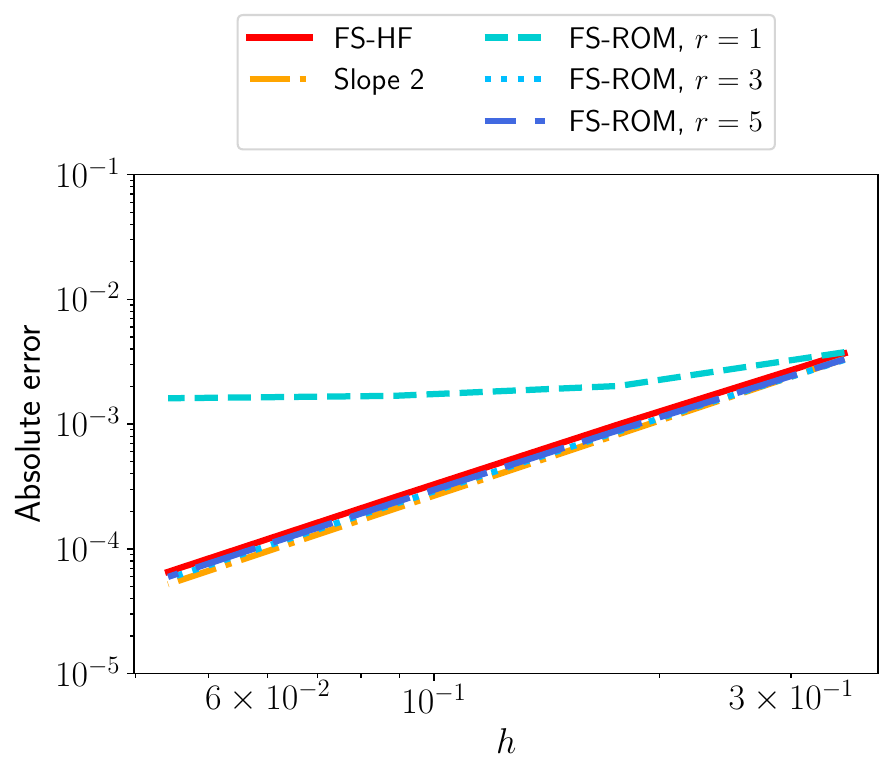}
        \caption{$\normz{\theta_\delta - \theta}$ for $\delta=h,r$.}
        \label{fig:test1_convergence_analysis_wrt_exact_theta_L2_fs}
    \end{subfigure}
    \caption{Example 1A: convergence of the errors for the fixed-stress iterative schemes (FS-HF and FS-ROM).}
    \label{fig:test1_convergence_analysis_wrt_exact_fs}
\end{figure}

As the first test case in Example 1, we validate the four distinct solvers discussed in the previous sections:
M-HF (Section \ref{sec: coupled}),
FS-HF (Section \ref{sec: fs}), M-ROM (Section \ref{sec: coupled ROM}), and FS-ROM (Section \ref{sec: fs ROM}).
In the domain $\Omega$ of the unit square {and for a final time $T=1$}, we employ an initial spatial discretization of $h=0.25$, {linear finite elements}, and a temporal discretization of $\dt = 0.0025$.
Subsequently, with each cycle, we reduce $h$ by half and divide $\dt$ by four. The convergence rate is then computed by employing corresponding norms to calculate the errors (maximum over $(0, T]$) against the provided analytical solutions in \eqref{eq:analytical solution}.

To validate the ROM algorithm, for each cycle, we initially run the HF solver to generate sequences $\left\{\buh^{n}\right\}_{n = 0}^{N}$, $\left\{\ph^{n}\right\}_{n = 0}^{N}$, $\left\{\thetah^{n}\right\}_{n = 0}^{N}$ in time. Then, we employ the POD algorithm as discussed in Section \ref{sec: POD} to compress these sequences, focusing on retaining the first five modes i.e., $r \leq 5$.
The ROM's performance is evaluated across five distinct scenarios, specifically for values of $r$ equal to 1, 2, 3, 4, and 5.

Results of the convergence analyses are plotted in Figure \ref{fig:test1_convergence_analysis_wrt_exact_cg} for the monolithic schemes (M-HF and M-ROM) and in Figure \ref{fig:test1_convergence_analysis_wrt_exact_fs} for the fixed-stress iterative schemes (FS-HF and FS-ROM).
Specifically,
the left column illustrates results in the $\normo{\cdot}$-norm, whereas the panels on the right show results in the $\normz{\cdot}$-norm. The convergence analysis for the displacement $\bu$, pressure $p$, and temperature $\theta$ are arranged in three rows within the figure. The first row pertains to the displacement $\bu$, the second row to the pressure $p$, and the third row to the temperature $\theta$.
For the M-HF scheme we observe a convergence rate of 2 for errors in the $\normz{\cdot}$-norm and of 1 for errors in the $\normo{\cdot}$-norm. Furthermore, since the solution of FS-HF converges to that of M-HF, as presented in Theorem \ref{thm:main}, FS-HF is expected to converge at the same rate as M-HF. Indeed, numerical results show that FS-HF achieves the same convergence rates.
Furthermore, with the exception of $r = 1$, the ROM scheme also achieves identical convergence rates to the HF scheme used in its training.

\begin{figure}[!h]
    \centering
    \begin{subfigure}[b]{0.48\textwidth}
        \centering
        \includegraphics[width=\textwidth, trim={0 0.35cm 0 0},clip]{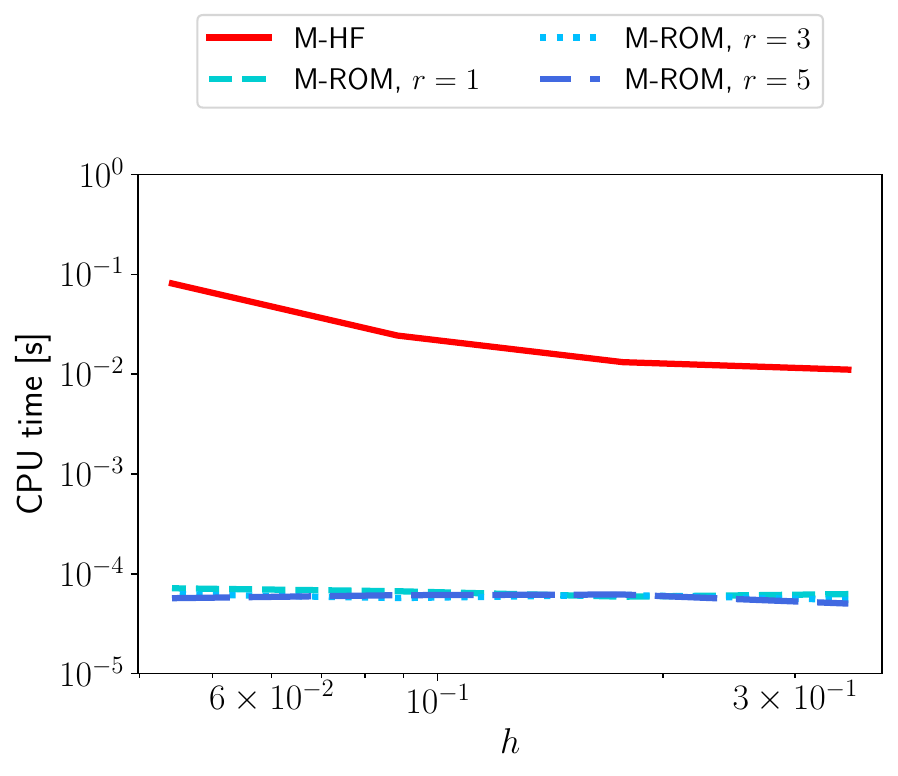}
        \caption{Monolithic schemes (M-HF and M-ROM).}
        \label{fig:test1_CPU_time_cg}
    \end{subfigure}
    \hfill
    \begin{subfigure}[b]{0.48\textwidth}
        \centering
        \includegraphics[width=\textwidth, trim={0 0.35cm 0 0},clip]{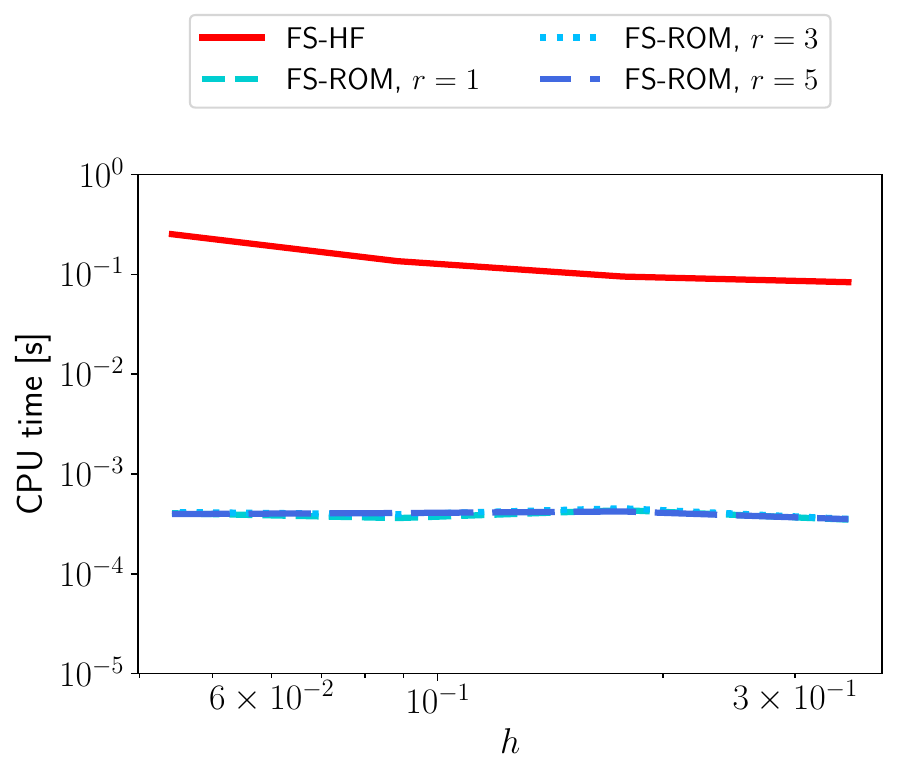}
        \caption{Iterative schemes (FS-HF and FS-ROM).}
        \label{fig:test1_CPU_time_fs}
    \end{subfigure}
    \caption{Example 1A: average CPU time per time step.}
    \label{fig:test1_CPU_time}
\end{figure}

\begin{figure}[!h]
    \centering
    \includegraphics[width=0.48\textwidth, trim={0 0.35cm 0 0},clip]{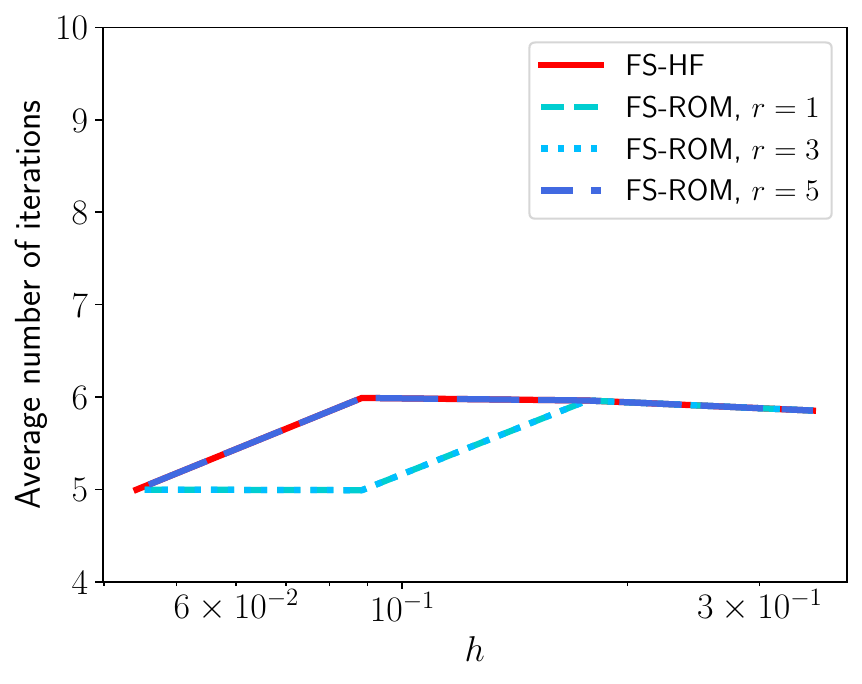}
    \caption{Example 1A: average number of fixed-stress iterations.}
    \label{fig:test1_iterations}
\end{figure}

Figure \ref{fig:test1_CPU_time} highlights one of the main advantages of the ROM schemes, which is the reduction of computational cost.
This advantage becomes evident through the significant reduction in the number of degrees of freedom achieved by the ROM schemes in comparison to the HF schemes. {By accompanying this reduction with the precomputation of matrices and vectors in Table \ref{table:rom_operators} we ensure that the leading cost in the evaluation of the ROM scales only with the ROM dimension $r$, and is not affected by the HF mesh size}. Consequently, ROMs exhibit a remarkable speedup of two orders of magnitude when compared to the underlying HF discretization for both monolithic and fixed-stress cases.

We further remark that, while decreasing the mesh size $h$ demands a higher CPU time for solving with the HF method, the results in Figure ~\ref{fig:test1_CPU_time} show that the ROM CPU time is not affected by mesh refinements. Indeed, the evaluation of the ROM only requires solutions of $r \times r$ linear systems for the FS-ROM or $3 r \times 3 r$ linear systems for the M-ROM regardless of the value of $h$. Moreover, increasing $r$ from 1 to 5 does not seem to increase the required CPU time either, since the solution of linear systems of dimension at most $15 \times 15$ is computationally inexpensive.

{On the other hand, a comparison between Figure \ref{fig:test1_CPU_time_cg} and Figure \ref{fig:test1_CPU_time_fs} highlights that the M-HF solver places lower computational demands than the FS-HF solver. Indeed, the ratio between the CPU time of the M-HF solver and the FS-HF one ranges from 0.13 (largest value of $h$) to 0.32 (smallest value of $h$). We anticipate that the M-HF solver's computational cost will eventually surpass that of the FS-HF solver in an asymptotic sense, but we refrain from further reducing the value of $h$ due to the computational costs.
A similar comparison illustrates that the M-ROM solver is an order of magnitude less expensive than the FS-ROM solver. Indeed, the ratio between the M-ROM CPU time and the FS-ROM one is around 0.15 for every value of $h$ and $r$. This is again due to the fact that, in the current realization, solving linear systems of dimension up to $15 \times 15$ requires a CPU time that is almost independent of the dimension of the linear system. The M-ROM solver needs to solve only one small linear system, while the FS-ROM requires solving 3 linear systems per iteration.}

Finally, Figure \ref{fig:test1_iterations} demonstrates that the average number of iterations necessary by FS-ROM for any given value of $r$ falls within the range of 5 to 6. Moreover, except for the case when $r = 1$, the FS-ROM exhibits an equivalent number of iterations as the FS-HF solver.

\FloatBarrier

\subsubsection{Example 1B. Comparison of HF and ROM solvers in time}

In the second test case in Example 1, we focus on the fixed stress iteration scheme and compare the errors between HF and ROM solvers.
Here, we use a mesh size of $h=1/16$, {linear finite elements}, and a time step size of $\dt = 0.001$.  We set the final time to be $T = 1$, which is larger than the one used in Example 1A.
We then train both M-ROM and FS-ROM schemes, compute at most $10$ POD modes, and evaluate the ROMs on the same time interval $(0, T]$ with the same time step size $\dt$. The larger number of POD modes selected in this test case compared to Example 1A is due to the fact that the final time is now ten times larger than the one employed in Example 1A.

\begin{figure}
    \centering
    \begin{subfigure}[b]{0.48\textwidth}
        \centering
        \includegraphics[width=\textwidth, trim={0 0.35cm 0 0},clip]{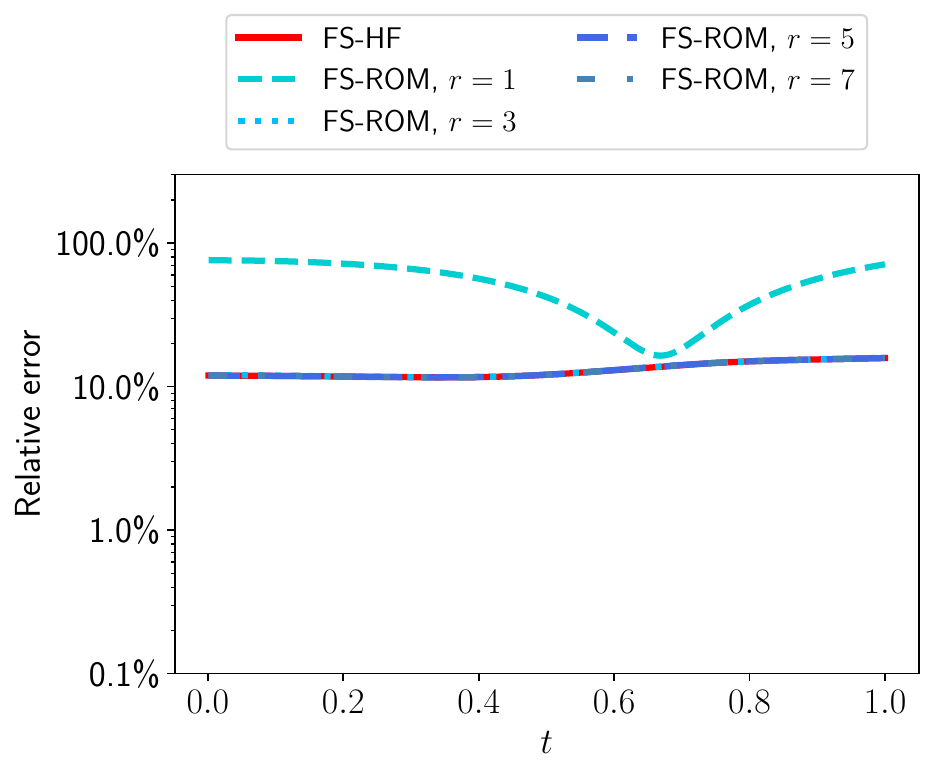}
        \caption{$\normo{\bu_\delta - \bu} / \normo{\bu}$.}
        \label{fig:test2_relative_error_wrt_exact_u_H1_fs}
    \end{subfigure}
    \hfill
    \begin{subfigure}[b]{0.48\textwidth}
        \centering
        \includegraphics[width=\textwidth, trim={0 0.35cm 0 0},clip]{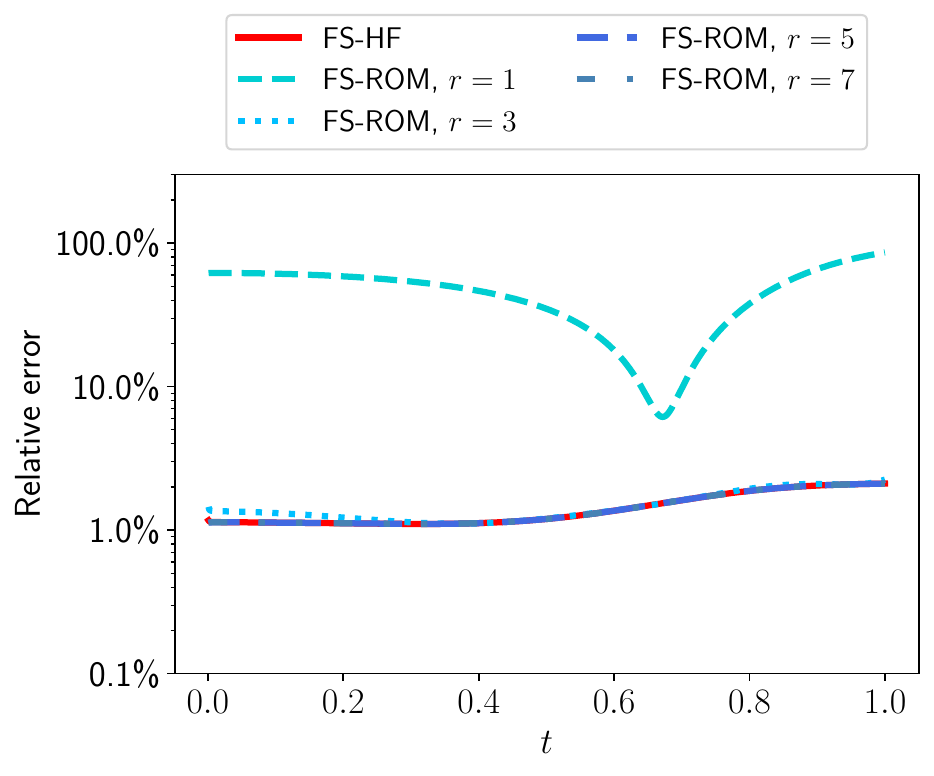}
        \caption{$\normz{\bu_\delta - \bu} / \normz{\bu}$.}
        \label{fig:test2_relative_error_wrt_exact_u_L2_fs}
    \end{subfigure}\\
    \centering
    \begin{subfigure}[b]{0.48\textwidth}
        \centering
        \adjincludegraphics[width=\textwidth,Clip={0} {0.15cm} {0} {.2\height}]{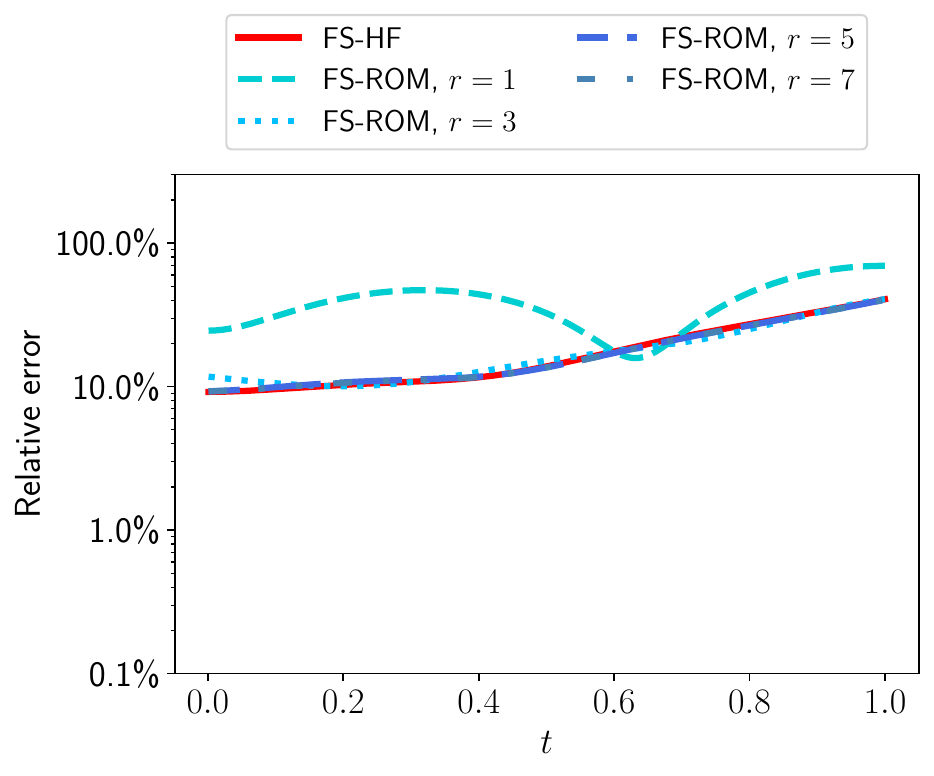}
        \caption{$\normo{p_\delta - p} / \normo{p}$.}
        \label{fig:test2_relative_error_wrt_exact_p_H1_fs}
    \end{subfigure}
    \hfill
    \begin{subfigure}[b]{0.48\textwidth}
        \centering
        \adjincludegraphics[width=\textwidth,Clip={0} {0.15cm} {0} {.2\height}]{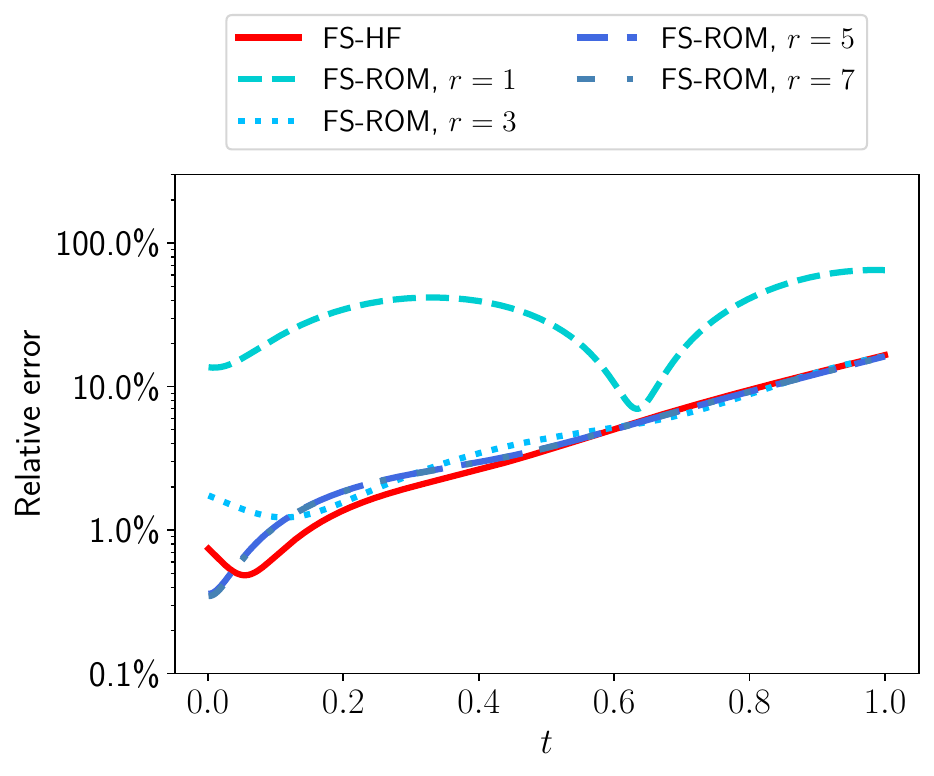}
        \caption{$\normz{p_\delta - p} / \normz{p}$.}
        \label{fig:test2_relative_error_wrt_exact_p_L2_fs}
    \end{subfigure}\\
    \centering
    \begin{subfigure}[b]{0.48\textwidth}
        \centering
        \adjincludegraphics[width=\textwidth,Clip={0} {0.15cm} {0} {.2\height}]{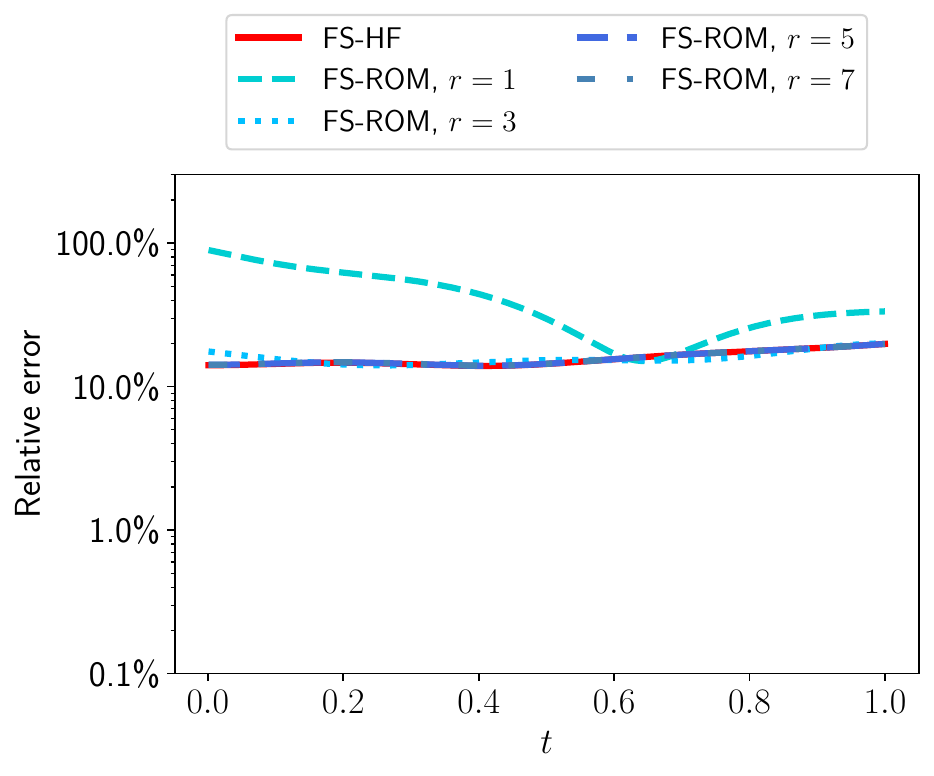}
        \caption{$\normo{\theta_\delta - \theta} / \normo{\theta}$.}
        \label{fig:test2_relative_error_wrt_exact_theta_H1_fs}
    \end{subfigure}
    \hfill
    \begin{subfigure}[b]{0.48\textwidth}
        \centering
        \adjincludegraphics[width=\textwidth,Clip={0} {0.15cm} {0} {.2\height}]{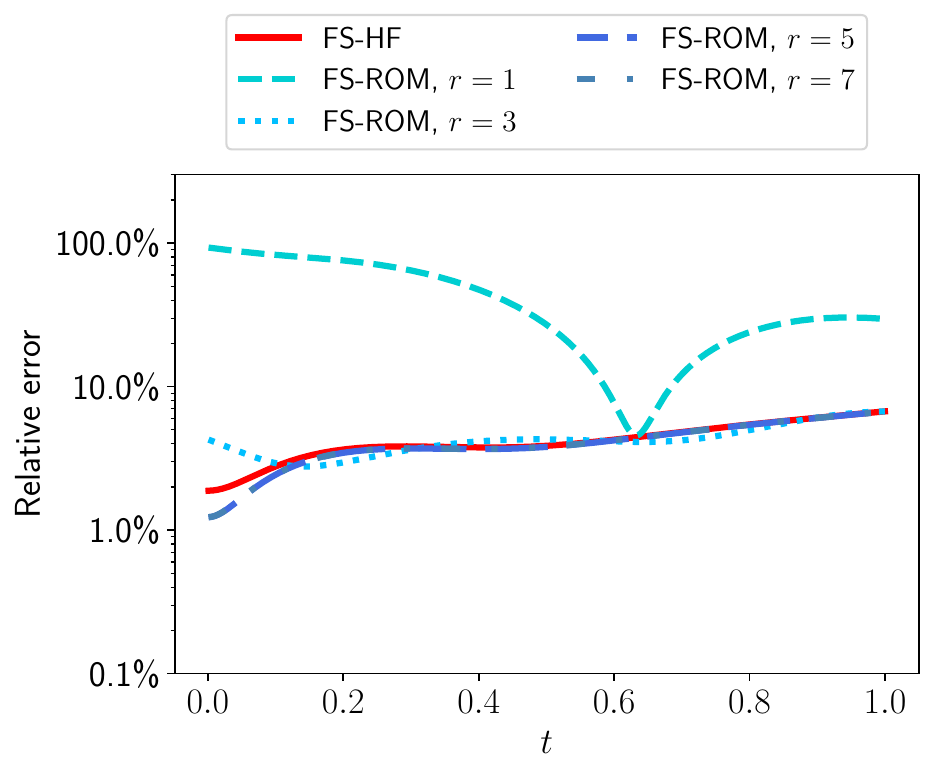}
        \caption{$\normz{\theta_\delta - \theta} / \normz{\theta}$.}
        \label{fig:test2_relative_error_wrt_exact_theta_L2_fs}
    \end{subfigure}
    \caption{Example 1B: relative errors ($\delta=h,r$) with respect to the analytical solution for the fixed-stress iterative schemes (FS-HF and FS-ROM).}
    \label{fig:test2_relative_error_wrt_exact_fs}
\end{figure}

First, Figure \ref{fig:test2_relative_error_wrt_exact_fs} illustrates the relative errors with respect to the analytical solutions in \eqref{eq:analytical solution} for the fixed-stress schemes.
In particular,  the left column presents the results in the $\normo{\cdot}$-norm, whereas the right column shows the results in the $\normz{\cdot}$-norm. The first row pertains to the displacement $\bu$, the second row to the pressure $p$, and the third row to the temperature $\theta$.
We observe that the ROM for $r = 1$ is inaccurate, with relative errors up to $100\%$ over the time interval $(0, T]$, regardless of the solution component $\bu, p$ or $\theta$. For $r \geq 3$, the ROM solutions are as accurate as the HF results, with only minor discrepancies appearing at the beginning of the time interval $(0, T]$ in Figures \ref{fig:test2_relative_error_wrt_exact_p_L2_fs} and \ref{fig:test2_relative_error_wrt_exact_theta_L2_fs}.
Furthermore, the error rises as time progresses.

\begin{figure}
    \centering
    \begin{subfigure}[b]{0.48\textwidth}
        \centering
        \includegraphics[width=\textwidth, trim={0 0.35cm 0 0},clip]{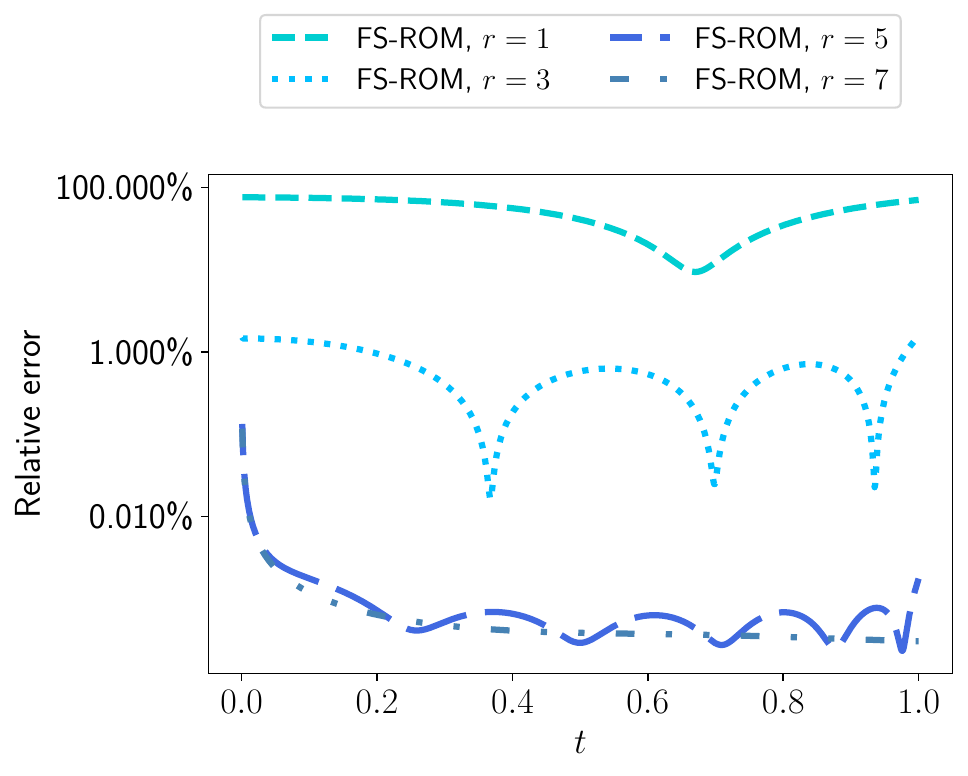}
        \caption{$\normo{\bu_r - \bu_h} / \normo{\bu_h}$.}
        \label{fig:test2_relative_error_wrt_FE_u_H1_fs}
    \end{subfigure}
    \hfill
    \begin{subfigure}[b]{0.48\textwidth}
        \centering
        \includegraphics[width=\textwidth, trim={0 0.35cm 0 0},clip]{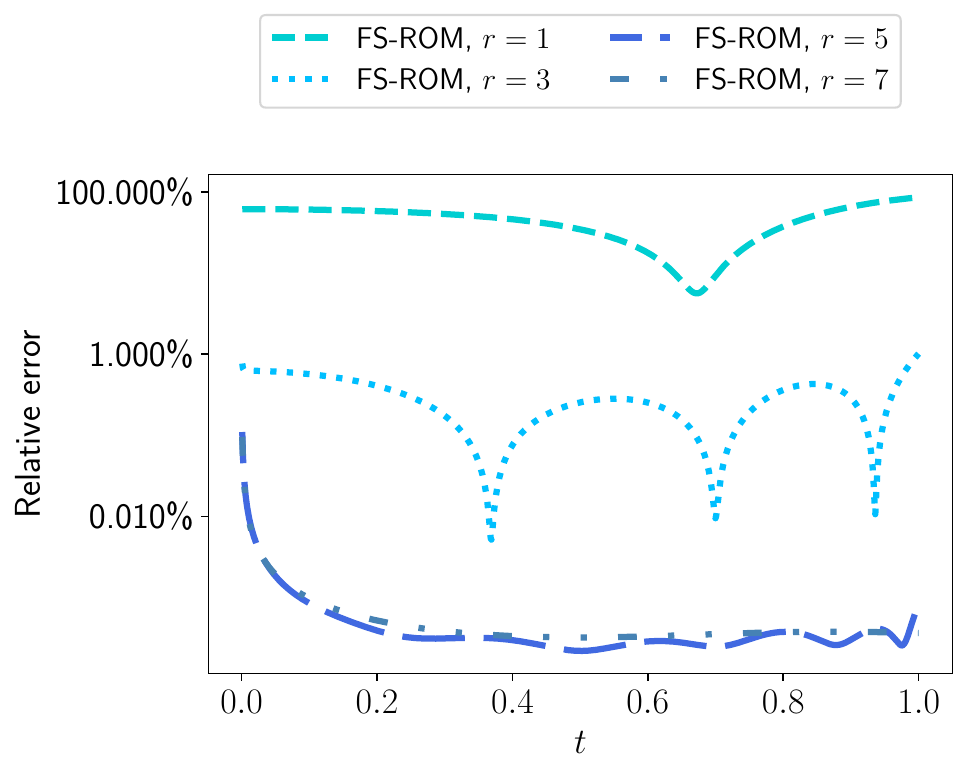}
        \caption{$\normz{\bu_r - \bu_h} / \normz{\bu_h}$.}
        \label{fig:test2_relative_error_wrt_FE_u_L2_fs}
    \end{subfigure}\\
    \centering
    \begin{subfigure}[b]{0.48\textwidth}
        \centering
        \adjincludegraphics[width=\textwidth,Clip={0} {0.15cm} {0} {.2\height}]{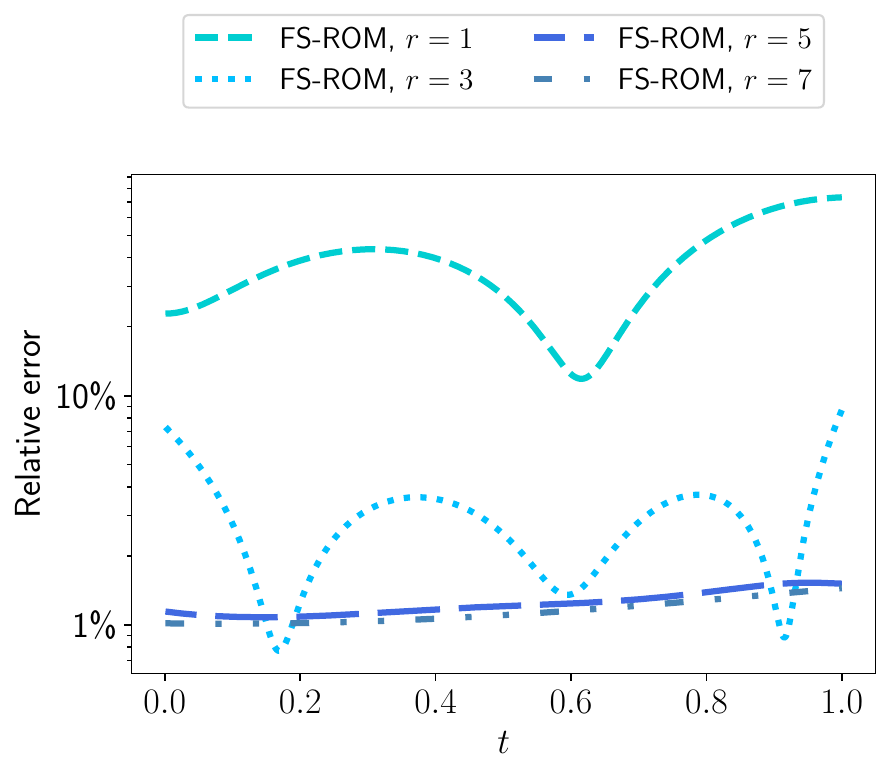}
        \caption{$\normo{p_r - p_h} / \normo{p_h}$.}
        \label{fig:test2_relative_error_wrt_FE_p_H1_fs}
    \end{subfigure}
    \hfill
    \begin{subfigure}[b]{0.48\textwidth}
        \centering
        \adjincludegraphics[width=\textwidth,Clip={0} {0.15cm} {0} {.2\height}]{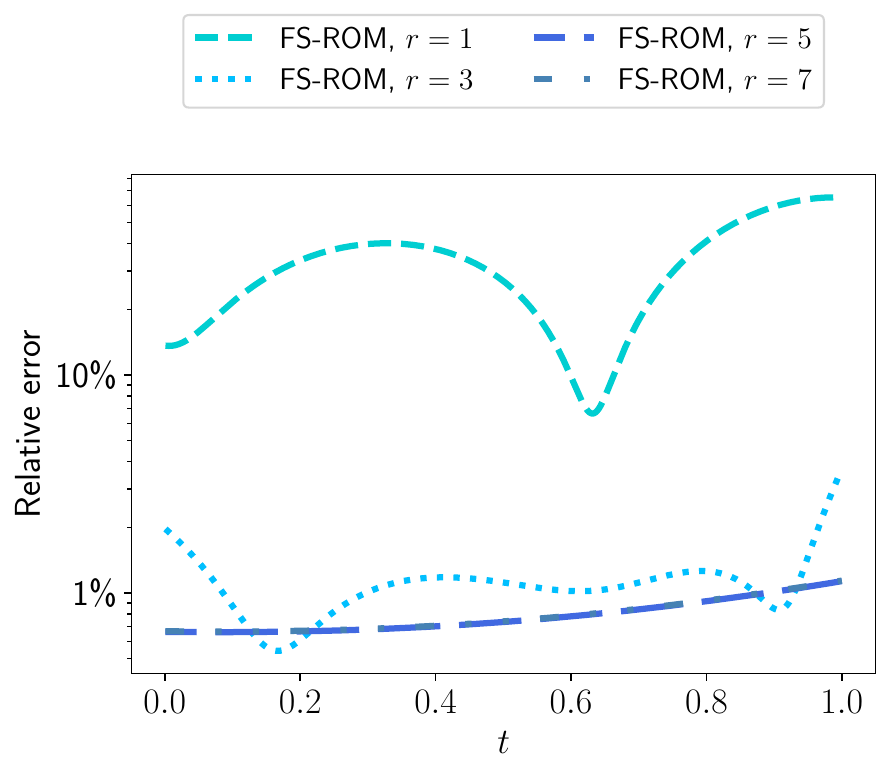}
        \caption{$\normz{p_r - p_h} / \normz{p_h}$.}
        \label{fig:test2_relative_error_wrt_FE_p_L2_fs}
    \end{subfigure}\\
    \centering
    \begin{subfigure}[b]{0.48\textwidth}
        \centering
        \adjincludegraphics[width=\textwidth,Clip={0} {0.15cm} {0} {.2\height}]{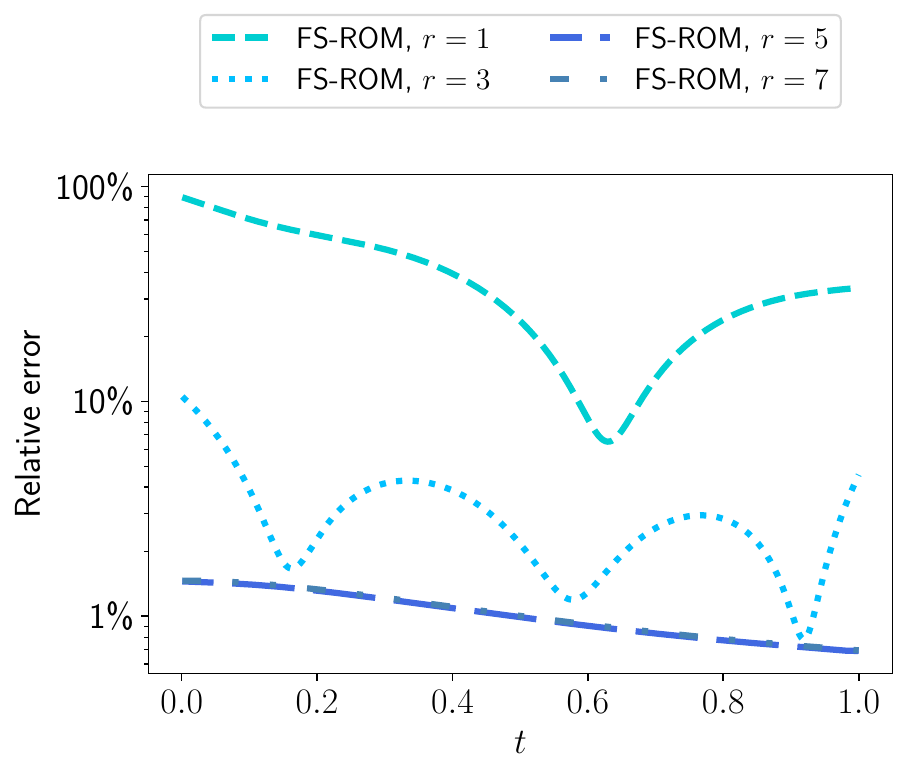}
        \caption{$\normo{\theta_r - \theta_h} / \normo{\theta_h}$.}
        \label{fig:test2_relative_error_wrt_FE_theta_H1_fs}
    \end{subfigure}
    \hfill
    \begin{subfigure}[b]{0.48\textwidth}
        \centering
        \adjincludegraphics[width=\textwidth,Clip={0} {0.15cm} {0} {.2\height}]{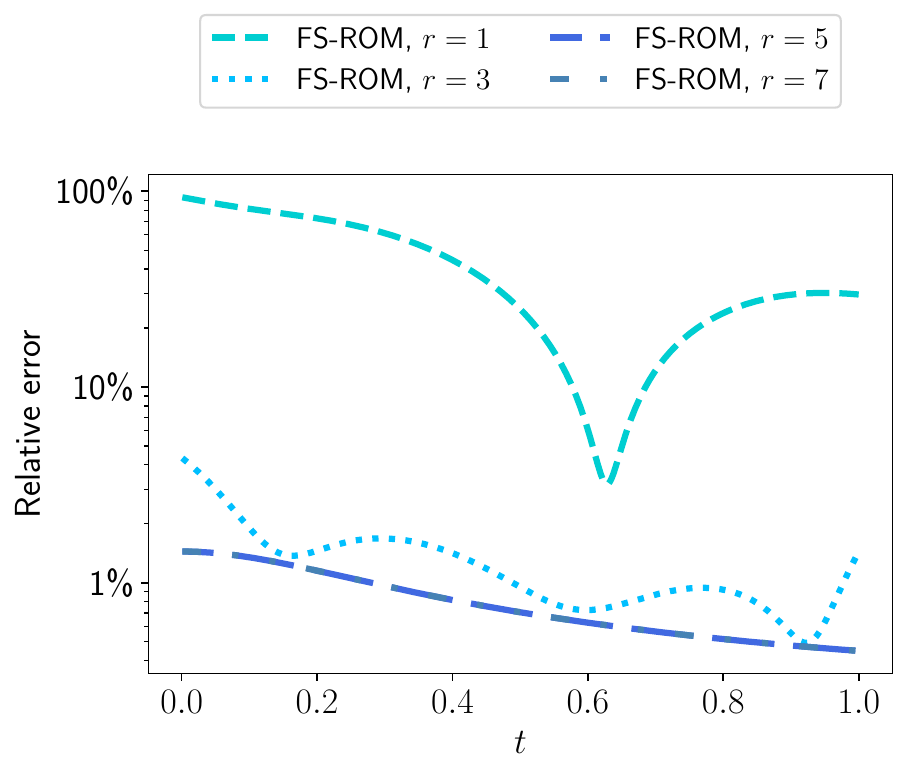}
        \caption{$\normz{\theta_r - \theta_h} / \normz{\theta_h}$.}
        \label{fig:test2_relative_error_wrt_FE_theta_L2_fs}
    \end{subfigure}\\
    \caption{Example 1B: relative errors between the FS-ROM and the FS-HF scheme.}
    \label{fig:test2_relative_error_wrt_FE_fs}
\end{figure}

\begin{figure}
    \centering
    \begin{subfigure}[b]{0.45\textwidth}
        \centering
        \includegraphics[width=\textwidth, trim={0 0.35cm 0 0},clip]{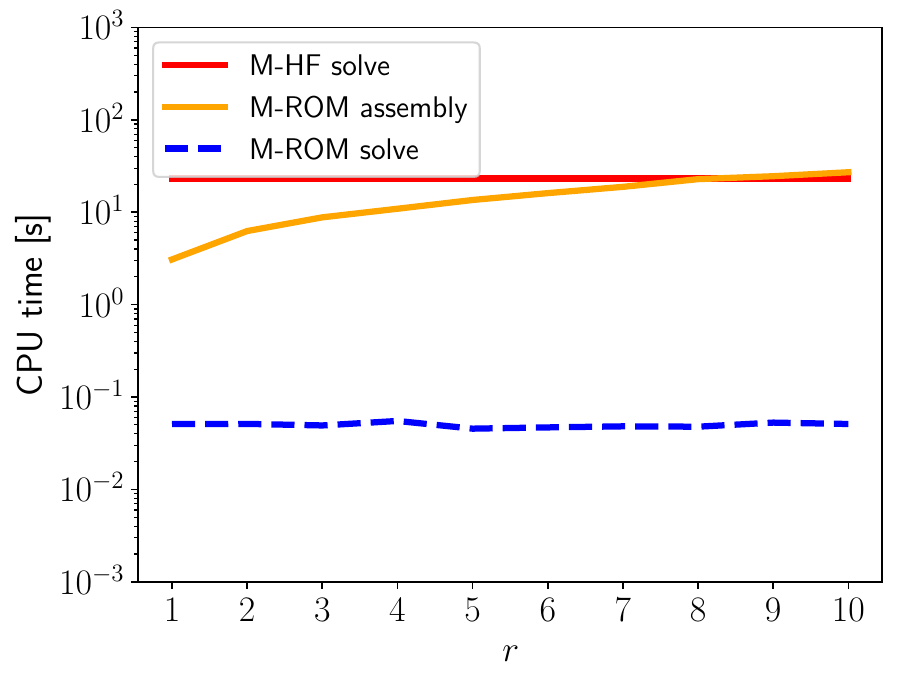}
        \caption{Monolithic schemes (M-HF vs. M-ROM).}
        \label{fig:test2_CPU_time_cg}
    \end{subfigure}
    \hfill
    \begin{subfigure}[b]{0.45\textwidth}
        \centering
        \includegraphics[width=\textwidth, trim={0 0.35cm 0 0},clip]{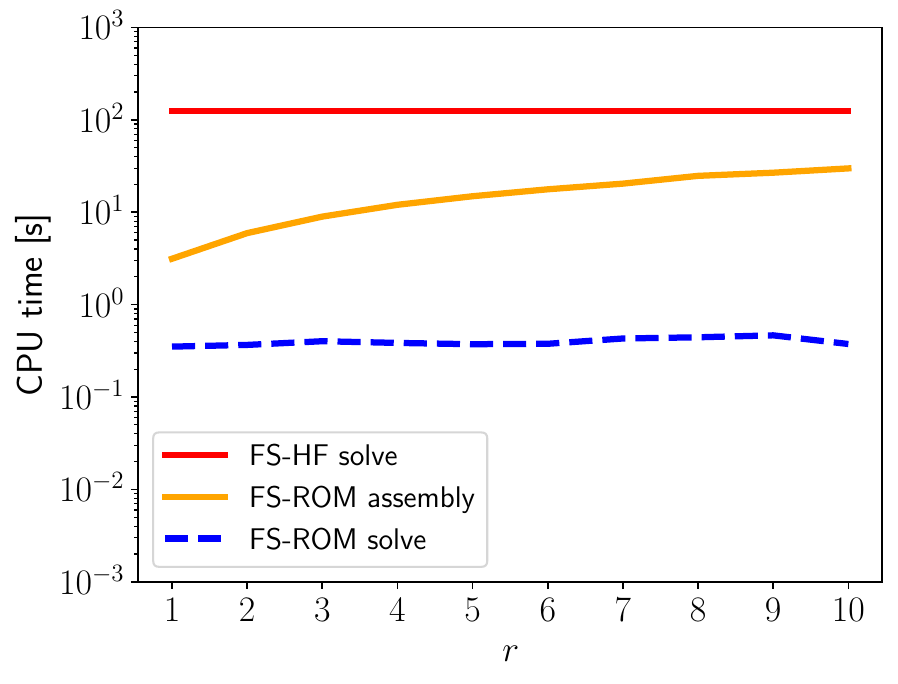}
        \caption{Iterative schemes (FS-HF vs. FS-ROM).}
        \label{fig:test2_CPU_time_fs}
    \end{subfigure}
    \caption{Example 1B: comparisons of total CPU times.}
    \label{fig:test2_CPU_time}
\end{figure}

\begin{figure}
    \centering
    \begin{subfigure}[b]{0.45\textwidth}
        \centering
        \includegraphics[width=\textwidth, trim={0 0.35cm 0 0},clip]{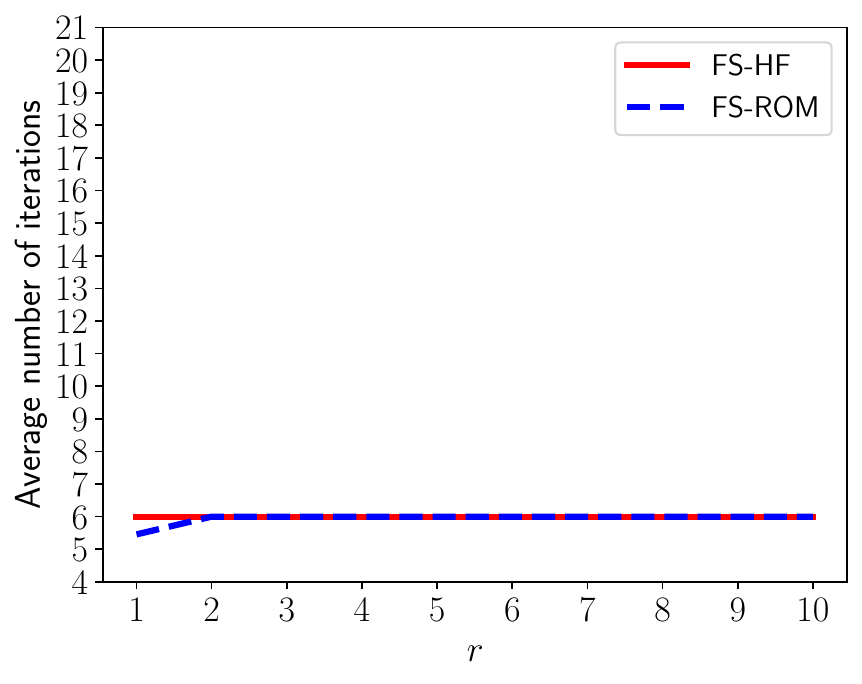}
        \caption{Average number of FS-ROM iterations.}
        \label{fig:test2_iterations}
    \end{subfigure}
    \hfill
    \begin{subfigure}[b]{0.475\textwidth}
        \centering
        \includegraphics[width=\textwidth, trim={0 0.35cm 0 0},clip]{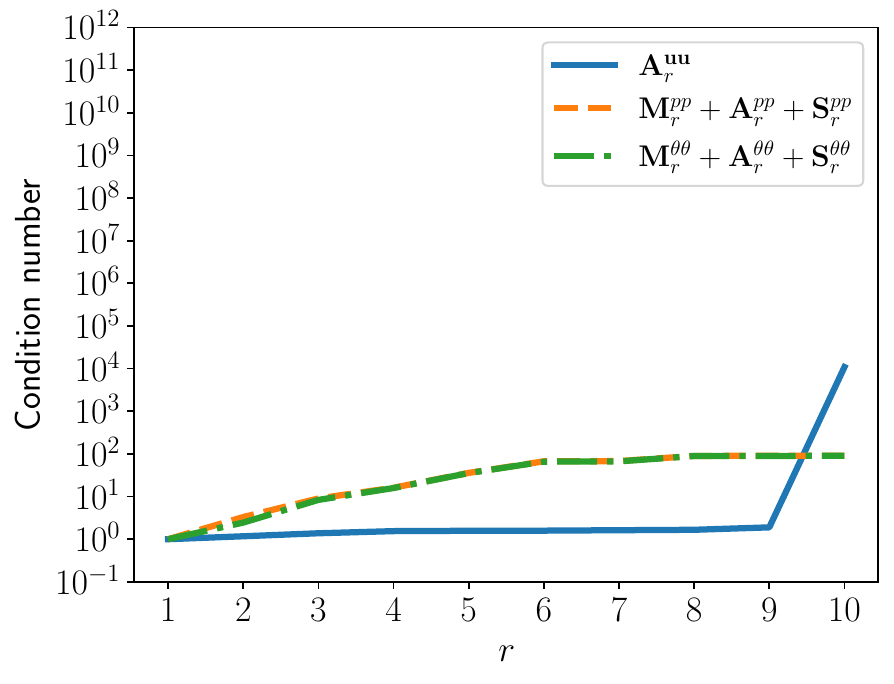}
        \caption{Condition numbers of FS-ROM iteration matrices.}
        \label{fig:test2_condition_numbers}
    \end{subfigure}\\
    \begin{subfigure}[b]{0.45\textwidth}
        \centering
        \includegraphics[width=\textwidth,trim={0 0.35cm 0 0},clip]{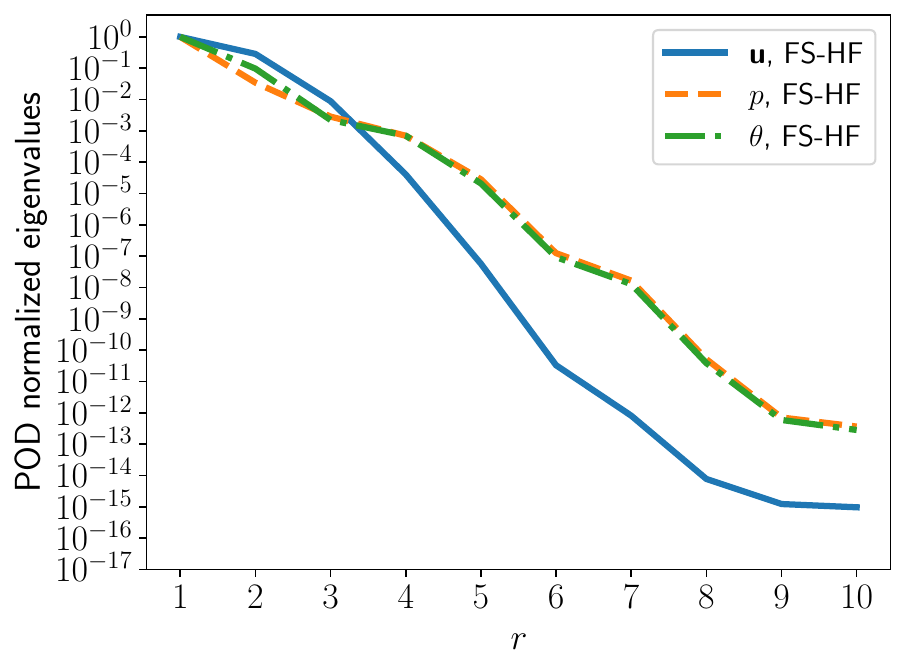}
        \caption{POD normalized eigenvalues $\nu_k / \nu_0$.}
        \label{fig:test2_POD_eigenvalues}
    \end{subfigure}
    \caption{Example 1B: (a) average number of fixed-stress iterations, (b) condition numbers of FS-ROM iteration matrices, and (c) POD normalized eigenvalues.}
    \label{fig:test2_iterations_condition_numbers_POD_eigenvalues}
\end{figure}

Next,  Figure \ref{fig:test2_relative_error_wrt_FE_fs} demonstrates the performance of the FS-ROM scheme by comparing its solutions to those of the FS-HF scheme. Figures \ref{fig:test2_relative_error_wrt_FE_u_H1_fs}-\ref{fig:test2_relative_error_wrt_FE_u_L2_fs} show that $\bur$ converges very fast to $\buh$  in the $\normo{\cdot}$ and $\normz{\cdot}$-norm as
 $r$ increases.
 In particular, we observe an error reduction of five orders of magnitude when increasing the value of $r$ from $r = 1$ to $r = 5$.
Similarly, Figures \ref{fig:test2_relative_error_wrt_FE_p_H1_fs}-\ref{fig:test2_relative_error_wrt_FE_p_L2_fs} show a decrease of more than one order of magnitude for $p$ within the first five POD modes, and Figures \ref{fig:test2_relative_error_wrt_FE_theta_H1_fs}-\ref{fig:test2_relative_error_wrt_FE_theta_L2_fs} a decrease of almost two orders of magnitude for $\theta$.

Furthermore, Figure \ref{fig:test2_CPU_time} reports the comparisons of the total CPU time for both HF and ROM solvers. Analogous to Example 1A, the ROM approach ensures a speedup of over two orders of magnitude, regardless of the specific value of $r$.
Figure \ref{fig:test2_CPU_time} also depicts the CPU time necessary for precomputation of the ROM operators summarized in Table \ref{table:rom_operators}. This observation leads us to conclude that the precomputation of the ROM operators at the end of the offline stage is essential to achieve the desired speedup. This becomes particularly crucial for larger values of $r$ as attempting to carry out ROM assembly during an online ROM solve would result in an online evaluation that is nearly as computationally intensive as the HF solve itself.

Lastly, we present the fixed-stress iteration counts for the FS-ROM in Figure \ref{fig:test2_iterations}. Similar to Example 1A, the FS-ROM exhibits an identical number of iterations when compared to the FS-HF solver.
In addition, Figure \ref{fig:test2_condition_numbers} presents the condition numbers for each matrix in Step i-ROM, Step ii-ROM, and Step iii-ROM. We observe that the condition numbers increase mildly while $r$ increases. Furthermore,  Figure \ref{fig:test2_POD_eigenvalues} presents the normalized eigenvalues of the POD compression, where $\nu_0$ is the largest eigenvalue; the fast decay of the eigenvalues justifies limiting the ROM dimension to $r = 10$ at most.

\FloatBarrier

\subsubsection{Example 1C. ROM on a larger time interval than the HF solver}

\begin{figure}
    \centering
    \begin{subfigure}[b]{0.48\textwidth}
        \centering 
        \includegraphics[width=\textwidth, trim={0 0.35cm 0 0},clip] {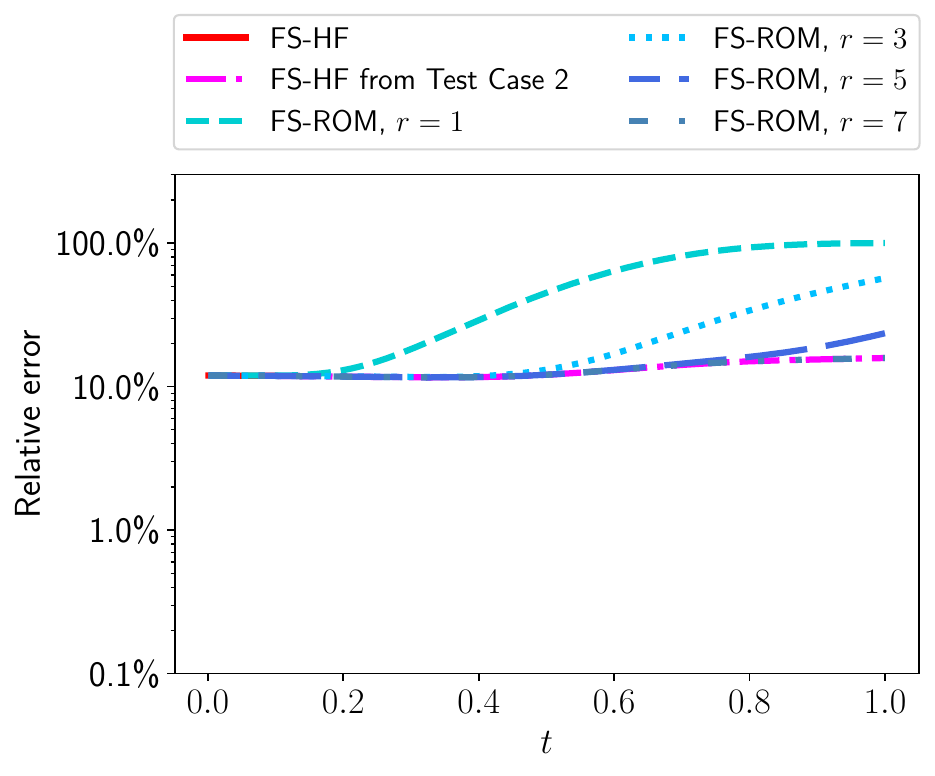}
        \caption{$\normo{\bu_\delta - \bu} / \normo{\bu}$.}
        \label{fig:test3_relative_error_wrt_exact_u_H1_fs}
    \end{subfigure}
    \hfill
    \begin{subfigure}[b]{0.48\textwidth}
        \centering
        \includegraphics[width=\textwidth, trim={0 0.35cm 0 0},clip]{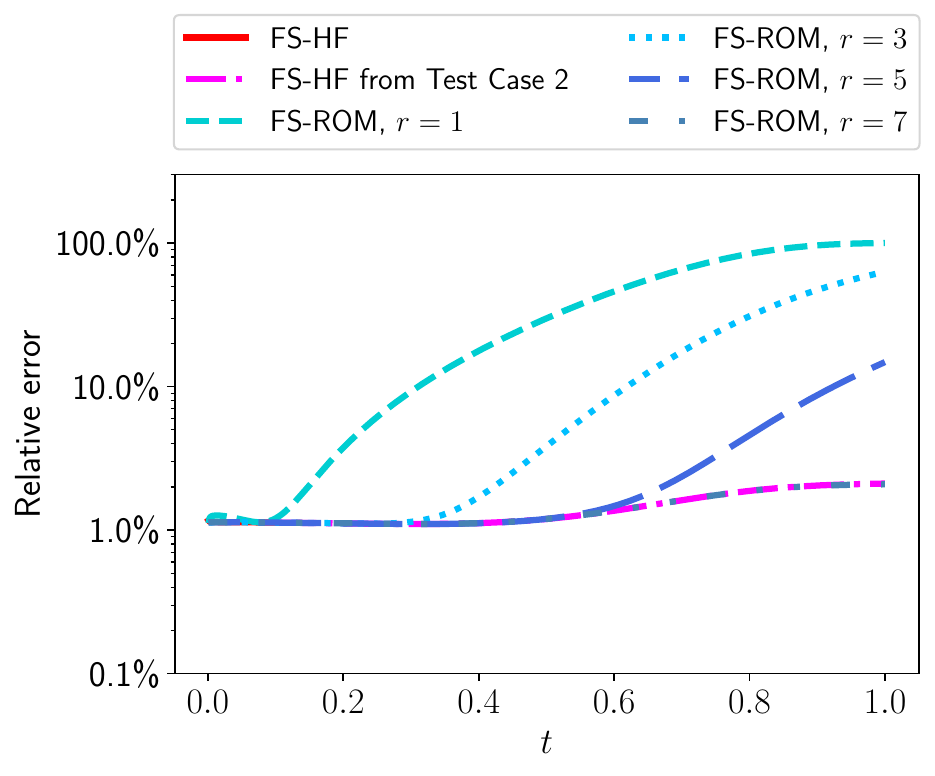}
        \caption{$\normz{\bu_\delta - \bu} / \normz{\bu}$.}
        \label{fig:test3_relative_error_wrt_exact_u_L2_fs}
    \end{subfigure}\\
    \centering
    \begin{subfigure}[b]{0.48\textwidth}
        \centering
        \adjincludegraphics[width=\textwidth,Clip={0} {0.15cm} {0} {.2\height}]{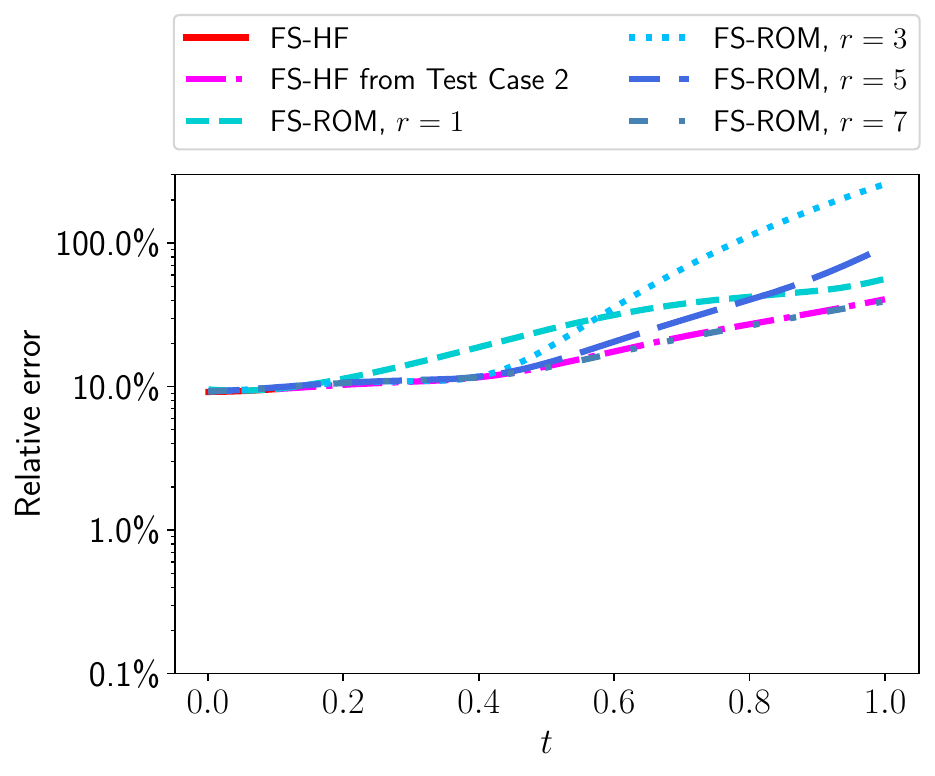}
        \caption{$\normo{p_\delta - p} / \normo{p}$}
        \label{fig:test3_relative_error_wrt_exact_p_H1_fs}
    \end{subfigure}
    \hfill
    \begin{subfigure}[b]{0.48\textwidth}
        \centering
        \adjincludegraphics[width=\textwidth,Clip={0} {0.15cm} {0} {.2\height}]{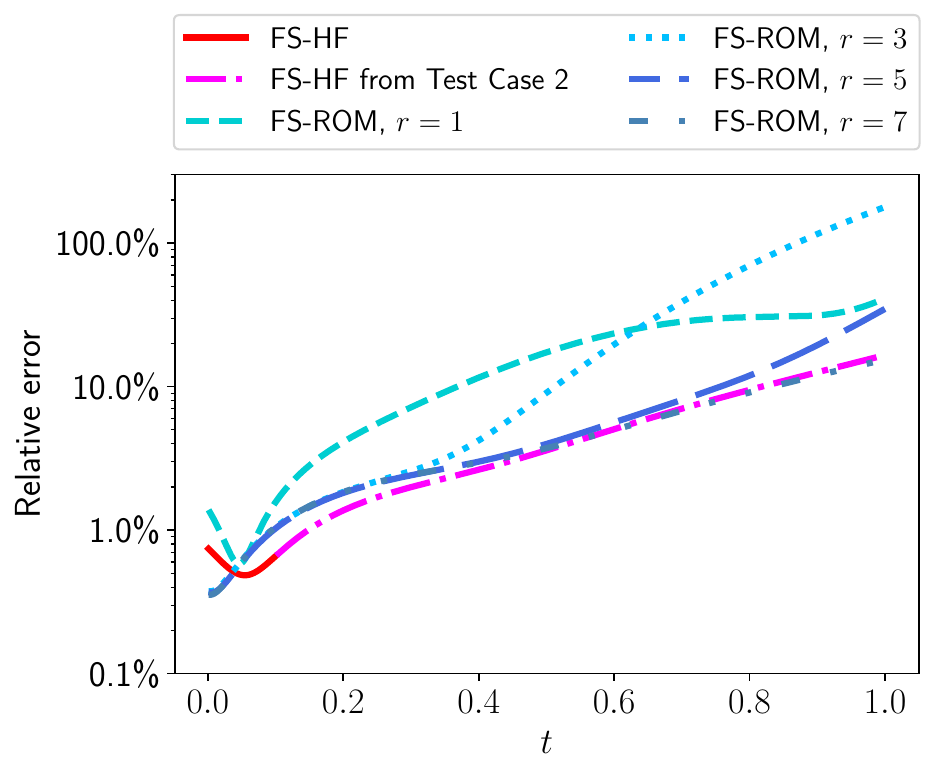}
        \caption{$\normz{p_\delta - p} / \normz{p}$.}
        \label{fig:test3_relative_error_wrt_exact_p_L2_fs}
    \end{subfigure}\\
    \centering
    \begin{subfigure}[b]{0.48\textwidth}
        \centering
        \adjincludegraphics[width=\textwidth,Clip={0} {0.15cm} {0} {.2\height}]{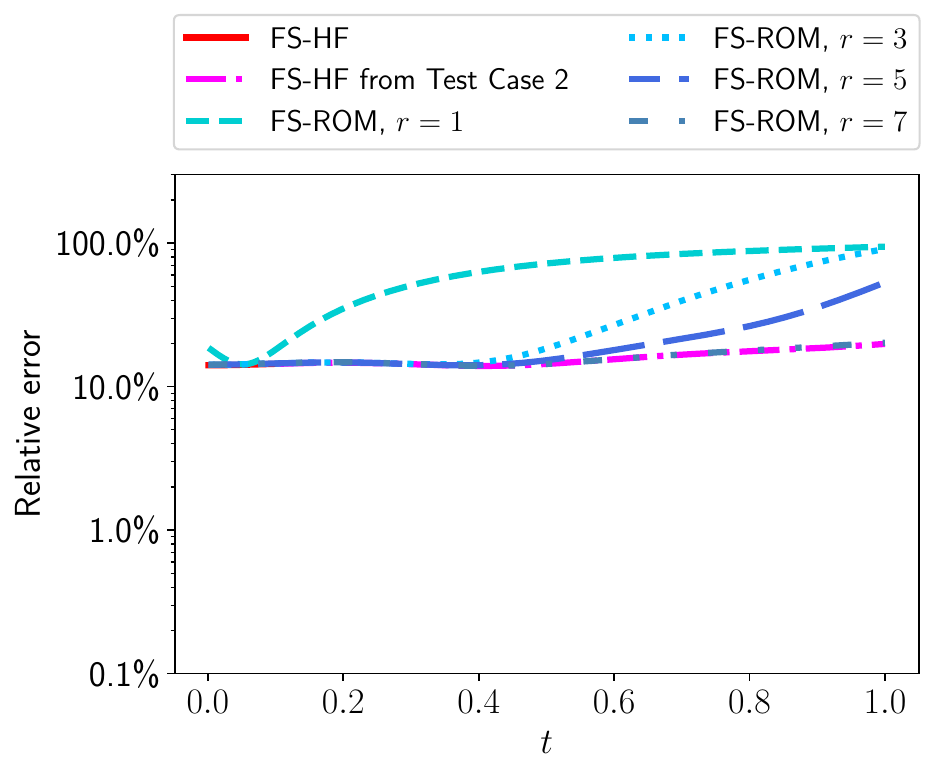}
        \caption{$\normo{\theta_\delta - \theta} / \normo{\theta}$.}
        \label{fig:test3_relative_error_wrt_exact_theta_H1_fs}
    \end{subfigure}
    \hfill
    \begin{subfigure}[b]{0.48\textwidth}
        \centering
        \adjincludegraphics[width=\textwidth,Clip={0} {0.15cm} {0} {.2\height}]{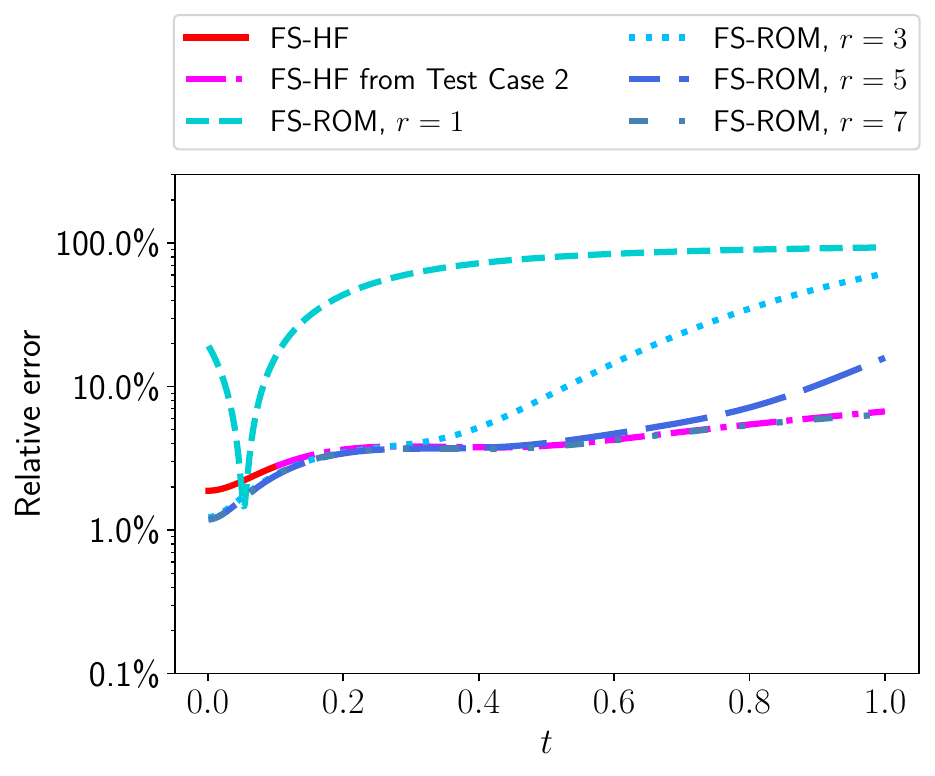}
        \caption{$\normz{\theta_\delta - \theta} / \normz{\theta}$.}
        \label{fig:test3_relative_error_wrt_exact_theta_L2_fs}
    \end{subfigure}\\
    \caption{Example 1C: relative errors for $\delta=h,r$ with respect to the analytical solution for the fixed-stress iterative schemes (FS-HF and FS-ROM).}
    \label{fig:test3_relative_error_wrt_exact_fs}
\end{figure}

\begin{figure}[!h]
    \centering
    \begin{subfigure}[b]{0.45\textwidth}
        \centering
        \includegraphics[width=\textwidth,trim={0 0.35cm 0 0},clip]{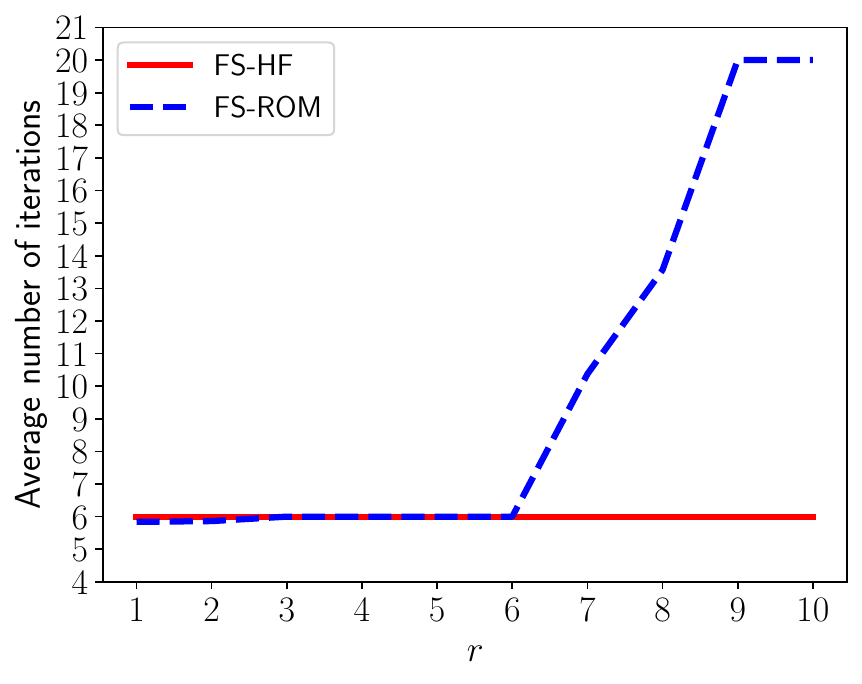}
        \caption{Average number of FS-ROM iterations.}
        \label{fig:test3_iterations}
    \end{subfigure}
    \hfill
    \begin{subfigure}[b]{0.475\textwidth}
        \centering
        \includegraphics[width=\textwidth,trim={0 0.35cm 0 0},clip]{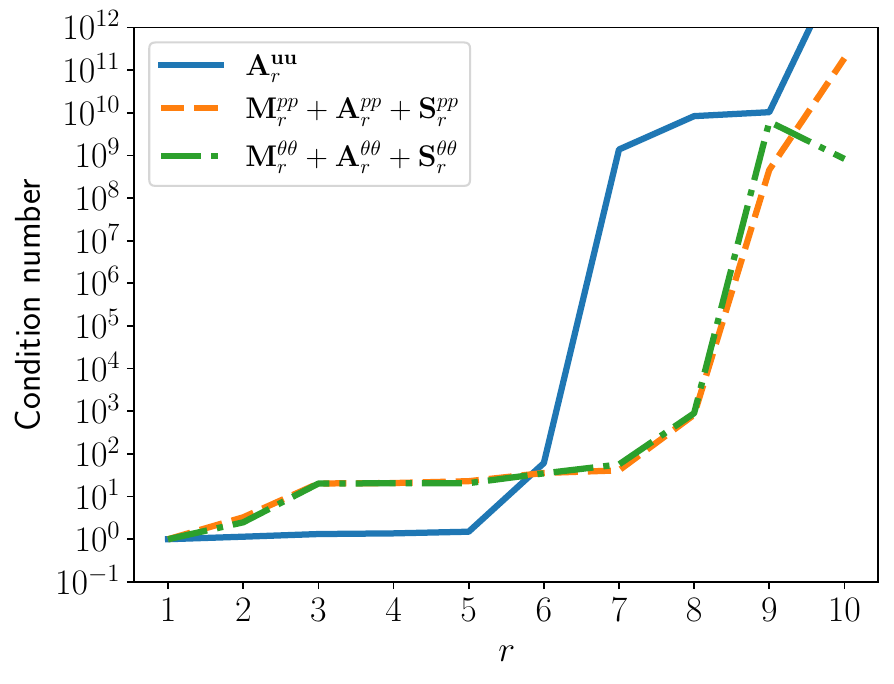}
        \caption{Condition numbers of FS-ROM iteration matrices.}
        \label{fig:test3_condition_numbers}
    \end{subfigure}\\
    \begin{subfigure}[b]{0.45\textwidth}
        \centering
        \includegraphics[width=\textwidth,trim={0 0.35cm 0 0},clip]{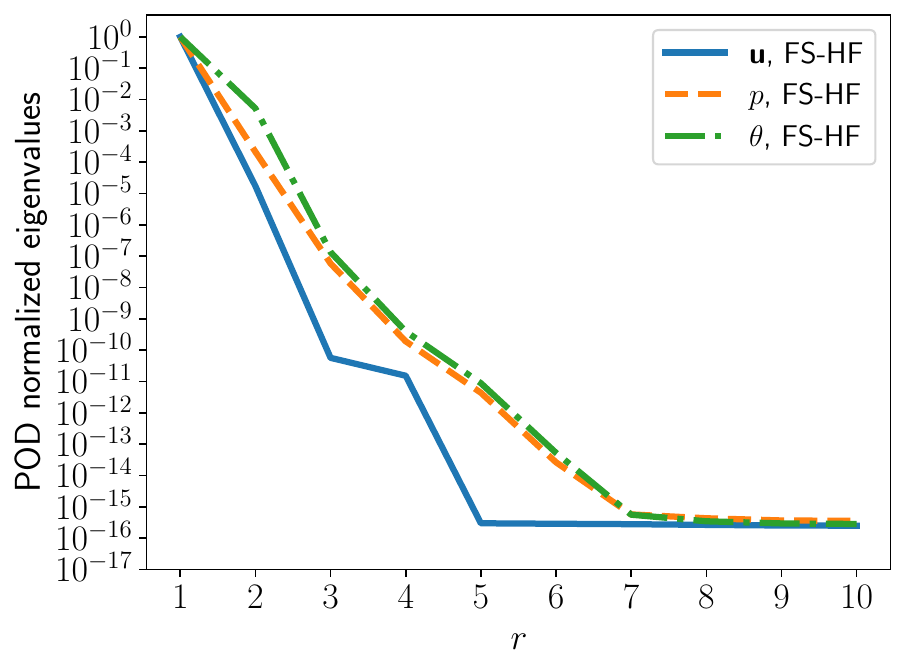}
        \caption{POD normalized eigenvalues $\nu_k / \nu_0$.}
        \label{fig:test3_POD_eigenvalues}
    \end{subfigure}
    \hfill
    \begin{subfigure}[b]{0.45\textwidth}
        \centering
     \includegraphics[width=\textwidth,trim={0 0.35cm 0 0},clip]{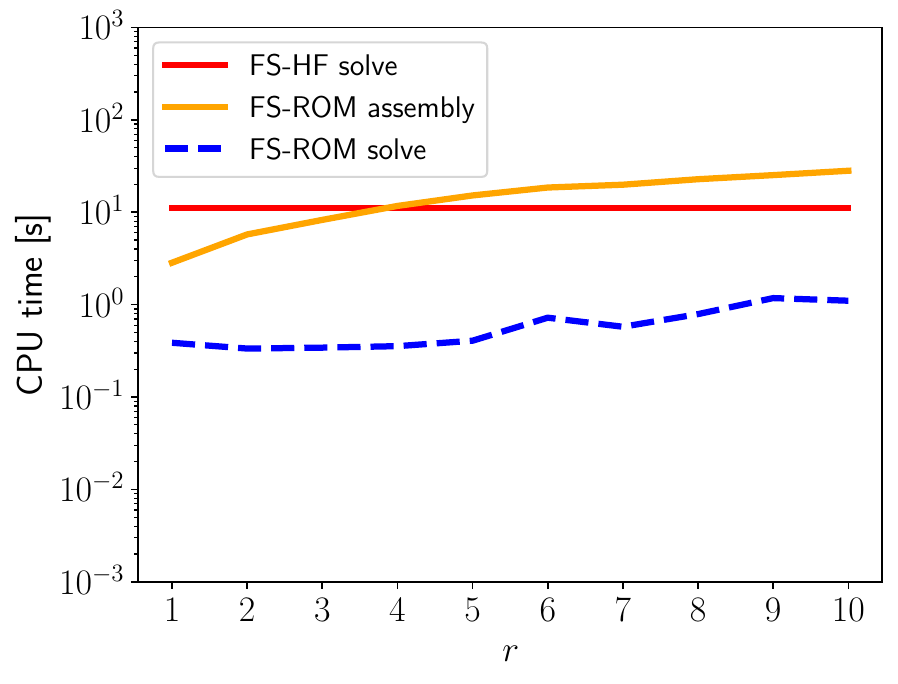}
        \caption{Total CPU time.}
        \label{fig:test3_CPU_time_fs}
        \end{subfigure}
    \caption{Example 1C: (a) average number of fixed-stress iterations, (b) condition numbers of FS-ROM iteration matrices, (c) POD normalized eigenvalues, and (d) total CPU time.}
    \label{fig:test3_iterations_condition_numbers_POD_eigenvalues}
\end{figure}

Employing ROM in Examples 1A and 1B serves well for validation purposes; however, its practical utility remains limited. This is because the ROM, {when trained on the time interval $(0, 1]$ and evaluated on the same time interval}, merely replicates pre-existing solutions within the same time span. In the context of the third test case, we delve into the ROM's extrapolation potential. Here, we train the model on the time span $(0, 0.1]$ and assess its performance on a more extensive interval, $(0, 1]$.

First, Figure \ref{fig:test3_relative_error_wrt_exact_fs} illustrates the relative errors with respect to the analytical solutions in \eqref{eq:analytical solution} for the FS-HF and FS-ROM schemes in the time interval $(0, 1]$. As in Figure~\ref{fig:test2_relative_error_wrt_FE_fs},  the left column presents the results in the $\normo{\cdot}$-norm, whereas the right column shows the results in the $\normz{\cdot}$-norm. The first row pertains to the displacement $\bu$, the second row to the pressure $p$, and the third row to the temperature $\theta$.
We note that the errors for the FS-HF solver are now limited to the interval $(0, 0.1]$ but the errors
displayed in the rest of the time interval $[0.1, 1]$ is only reported for comparison purposes from the previous Example 1B.

Comparing the results in Figure \ref{fig:test3_relative_error_wrt_exact_fs} to those obtained in Figure \ref{fig:test2_relative_error_wrt_exact_fs} for Example 1B, one can notice that a ROM with $r = 5$ is not as accurate in Example 1B due to the errors from extrapolation. Still, increasing the reduced basis size to $r = 7$ allows us to obtain a solution that is as accurate as the solution of the FS-HF scheme.

Although this test case demonstrates a successful extrapolation in accuracy, such an achievement introduces additional challenges for the ROM.
Figure \ref{fig:test3_iterations} indicates that the average number of fixed-stress iterations from the FS-ROM is greater than that of the FS-HF solver for $r \geq 7$. In addition, for $r \geq 9$, the FS-ROM scheme surpasses the maximum allowable limit of 20 iterations at each time step. It is important to emphasize that this behavior does not contravene Theorem \ref{thm:main ROM}; rather, it stems from limitations in arithmetic precision.
To support this statement, we plot
 the condition numbers of the matrices appearing on the left-hand side of Step i-ROM, Step ii-ROM, and Step iii-ROM (from the algorithm \eqref{sys: fs ROMeq}) in Figure \ref{fig:test3_condition_numbers}.
 We notice that the condition numbers of the matrices from Step iii-ROM is above $10^9$ for $r \geq 7$ and the condition numbers of all matrices from Steps i-ROM, ii-ROM, and iii-ROM are above $10^9$ for $r \geq 9$. From Figure \ref{fig:test3_POD_eigenvalues} we understand that this behavior is due to having added to the reduced basis some trailing POD modes associated with eigenvalues which, in machine precision, are zero. Since such eigenvalues appear on the denominator of \eqref{eq: POD mode}, the amplification of numerical errors causes the trailing POD modes not to satisfy \eqref{eq: POD orthonormal} anymore and may result in linearly dependent POD modes.
This behavior can be attributed to the use of a very small time interval during the training process. To observe this, compare Figure \ref{fig:test3_condition_numbers} with Figure \ref{fig:test2_condition_numbers}, as well as Figure \ref{fig:test3_POD_eigenvalues} with Figure \ref{fig:test2_POD_eigenvalues}.
The increase in the number of fixed stress iterations also results in a mild increase of the CPU time required by the ROM, as shown in Figure \ref{fig:test3_CPU_time_fs} for $r \geq 7$.

\FloatBarrier

\subsubsection{Example 1D. ROM using a larger time step than the HF solver}

\begin{figure}
    \centering
    \begin{subfigure}[b]{0.48\textwidth}
        \centering
        \includegraphics[width=\textwidth, trim={0 0.35cm 0 0},clip]{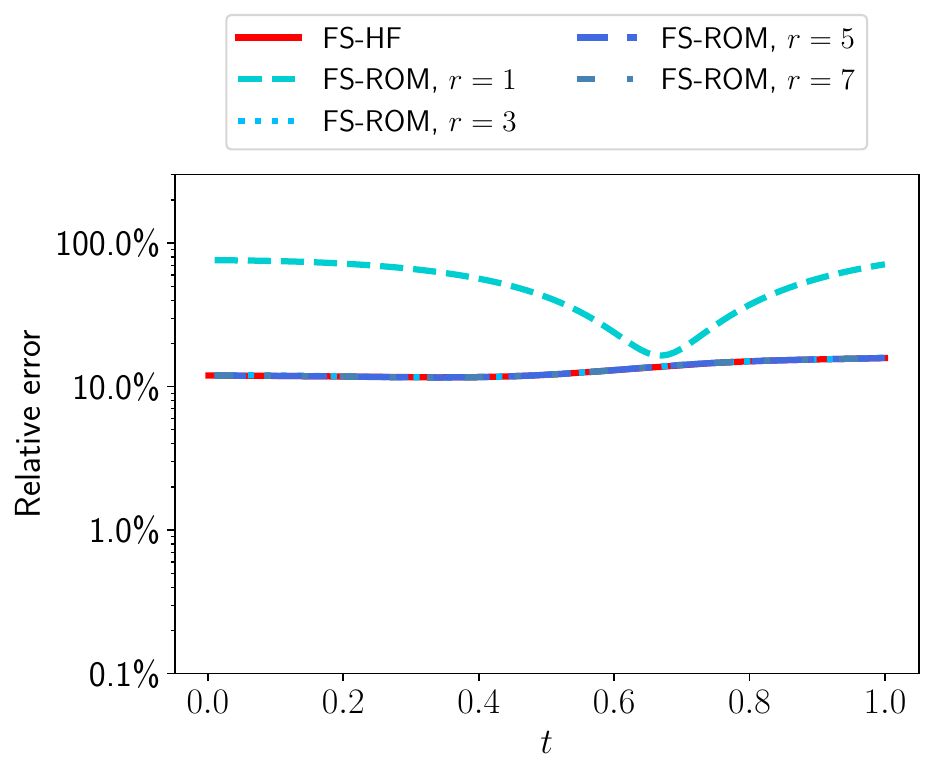}
        \caption{$\normo{\bu_\delta - \bu} / \normo{\bu}$.}
        \label{fig:test4_relative_error_wrt_exact_u_H1_fs}
    \end{subfigure}
    \hfill
    \begin{subfigure}[b]{0.48\textwidth}
        \centering
        \includegraphics[width=\textwidth, trim={0 0.35cm 0 0},clip]{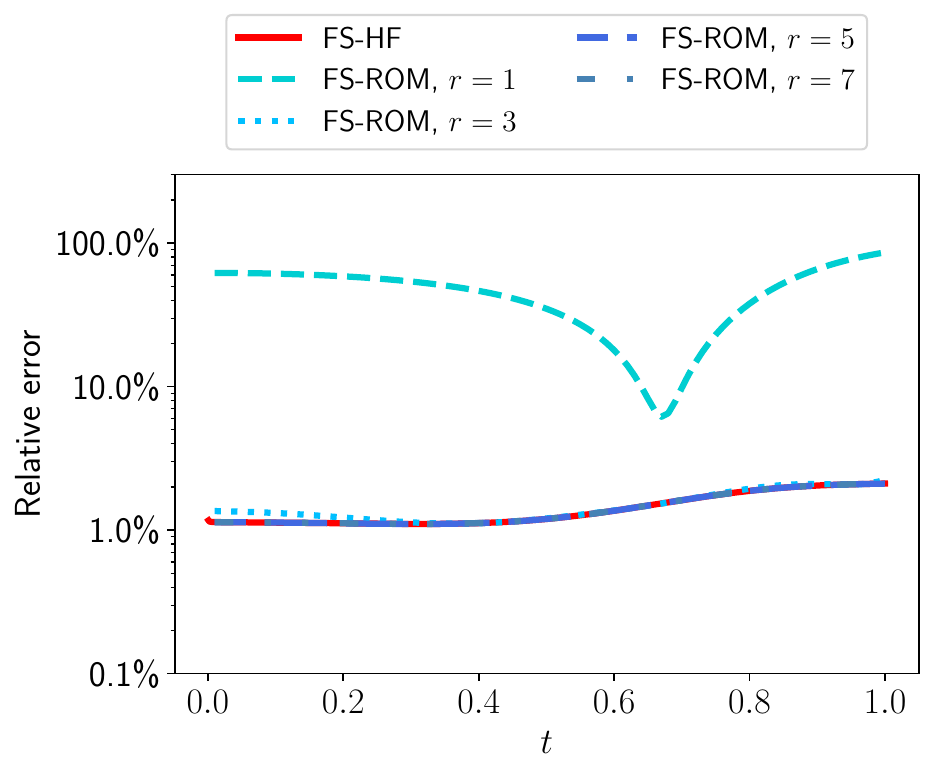}
        \caption{$\normz{\bu_\delta - \bu} / \normz{\bu}$.}
        \label{fig:test4_relative_error_wrt_exact_u_L2_fs}
    \end{subfigure}\\
    \centering
    \begin{subfigure}[b]{0.48\textwidth}
        \centering
        \adjincludegraphics[width=\textwidth,Clip={0} {0.15cm} {0} {.2\height}]{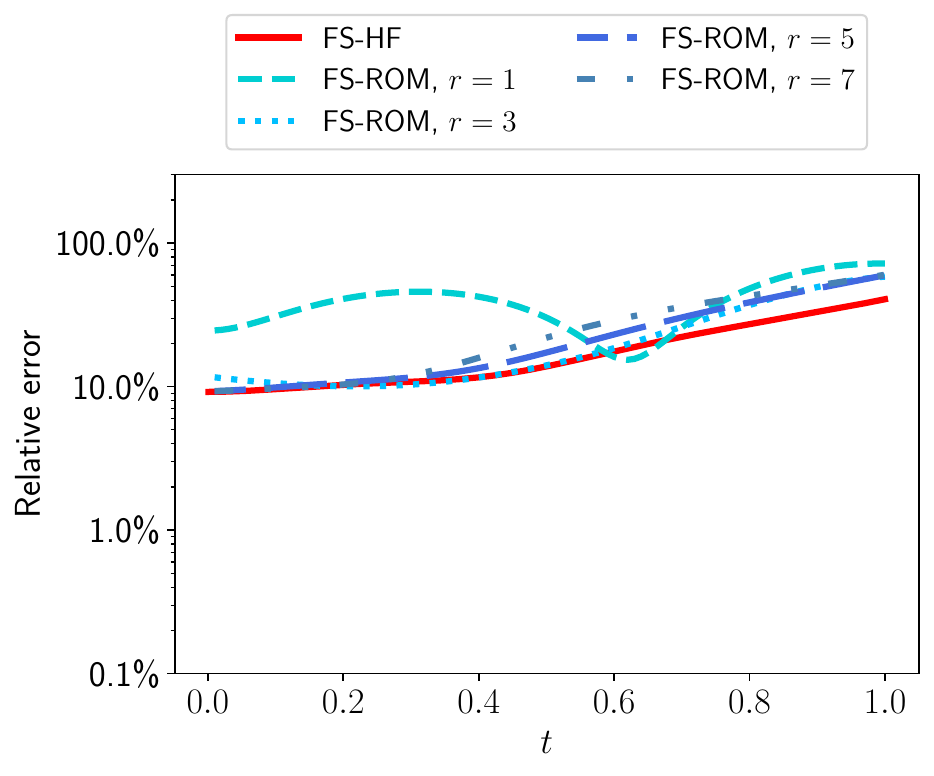}
        \caption{$\normo{p_\delta - p} / \normo{p}$.}
        \label{fig:test4_relative_error_wrt_exact_p_H1_fs}
    \end{subfigure}
    \hfill
    \begin{subfigure}[b]{0.48\textwidth}
        \centering
        \adjincludegraphics[width=\textwidth,Clip={0} {0.15cm} {0} {.2\height}]{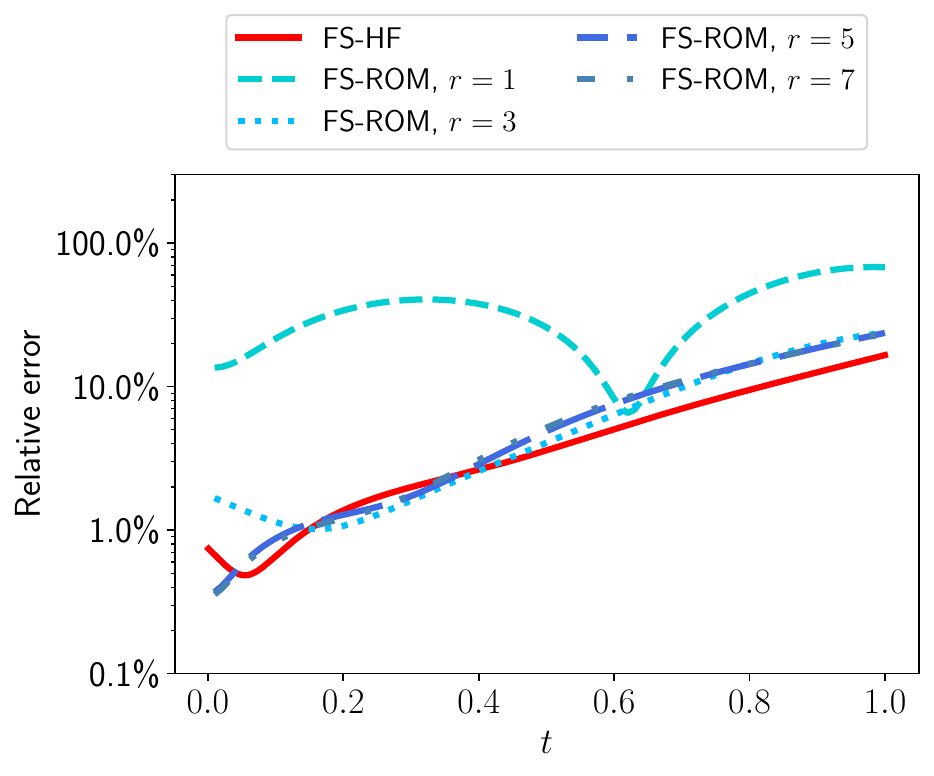}
        \caption{$\normz{p_\delta - p} / \normz{p}$.}
        \label{fig:test4_relative_error_wrt_exact_p_L2_fs}
    \end{subfigure}\\
    \centering
    \begin{subfigure}[b]{0.48\textwidth}
        \centering
        \adjincludegraphics[width=\textwidth,Clip={0} {0.15cm} {0} {.2\height}]{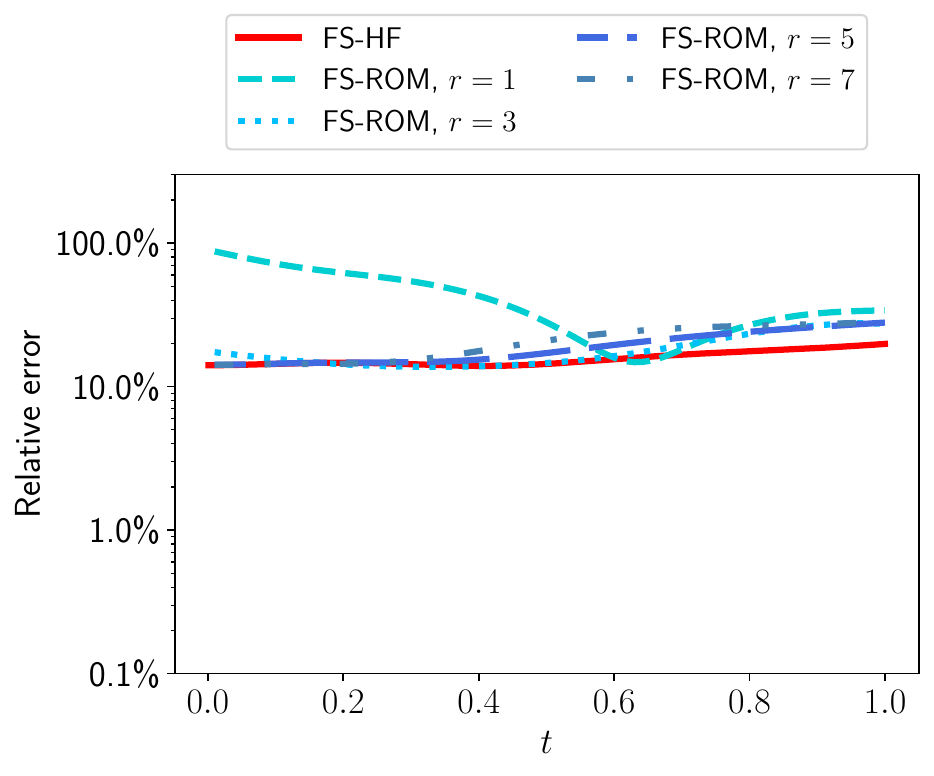}
        \caption{$\normo{\theta_\delta - \theta} / \normo{\theta}$}
        \label{fig:test4_relative_error_wrt_exact_theta_H1_fs}
    \end{subfigure}
    \hfill
    \begin{subfigure}[b]{0.48\textwidth}
        \centering
        \adjincludegraphics[width=\textwidth,Clip={0} {0.15cm} {0} {.2\height}]{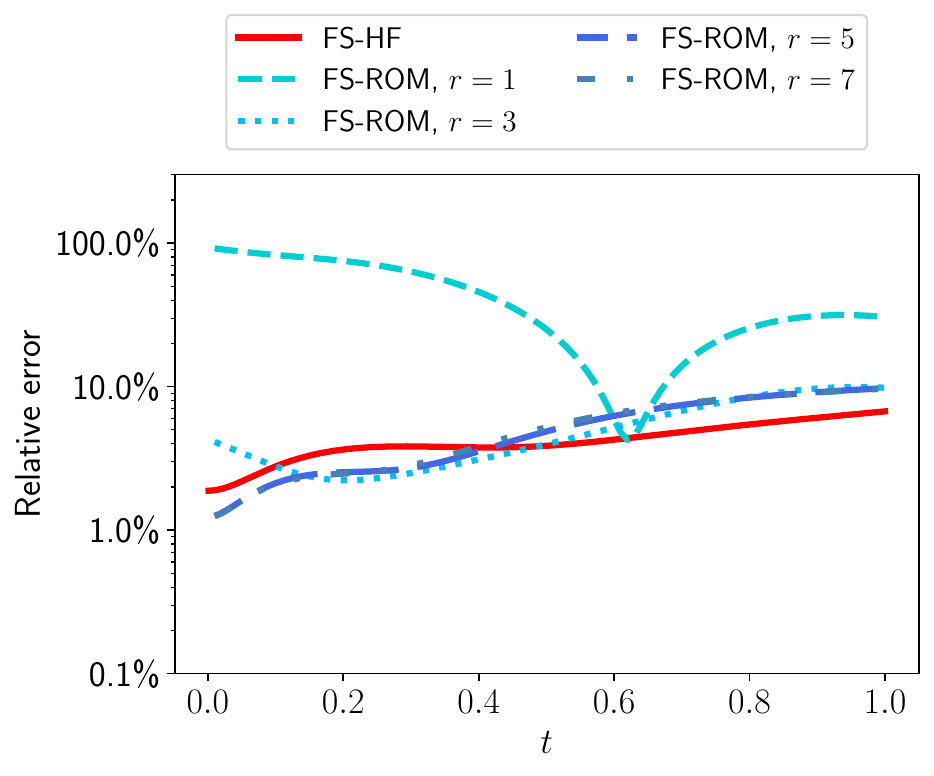}
        \caption{$\normz{\theta_\delta - \theta} / \normz{\theta}$.}
        \label{fig:test4_relative_error_wrt_exact_theta_L2_fs}
    \end{subfigure}
    \caption{Example 1D: relative errors  ($\delta = h,r$) with respect to the analytical solution for the fixed-stress iterative schemes (FS-HF and FS-ROM).}
    \label{fig:test4_relative_error_wrt_exact_fs}
\end{figure}

\begin{figure}
    \centering
    \begin{subfigure}[b]{0.48\textwidth}
        \centering
        \includegraphics[width=\textwidth]{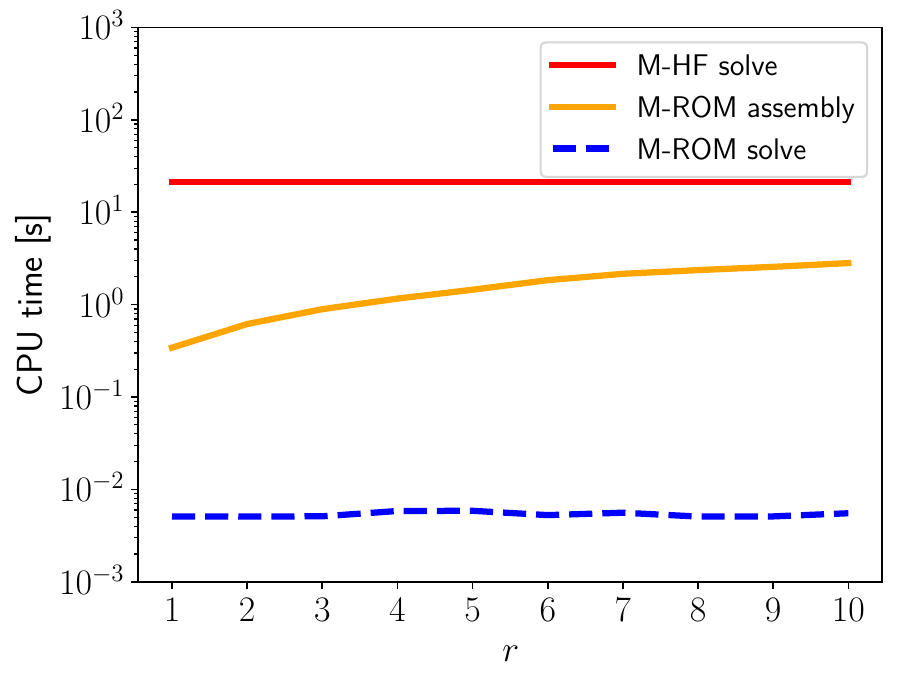}
        \caption{Monolithic schemes (M-HF and M-ROM).}
        \label{fig:test4_CPU_time_cg}
    \end{subfigure}
    \hfill
    \begin{subfigure}[b]{0.48\textwidth}
        \centering
        \includegraphics[width=\textwidth]{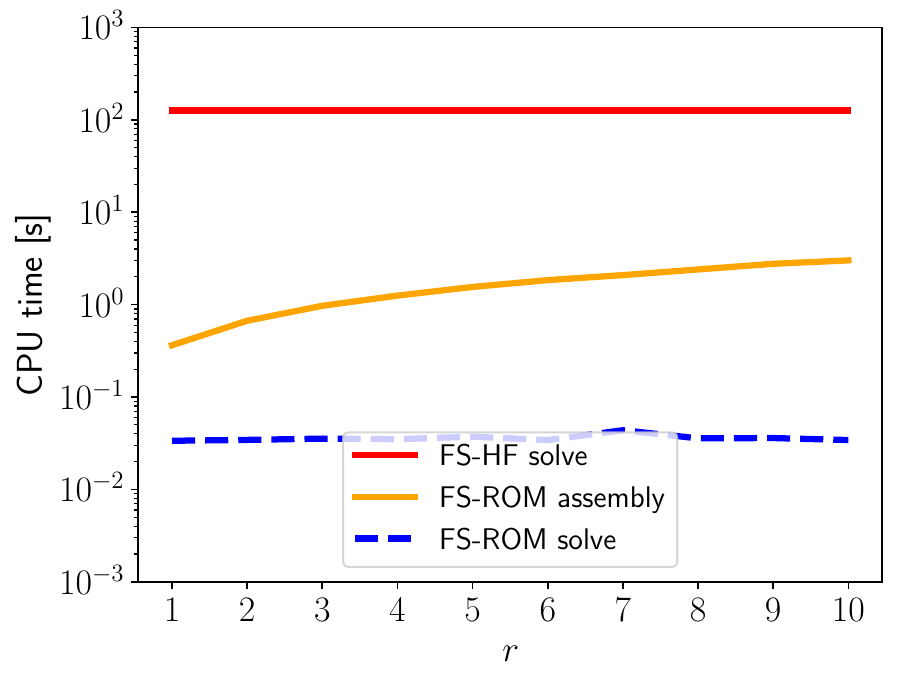}
        \caption{Iterative schemes (FS-HF and FS-ROM).}
        \label{fig:test4_CPU_time_fs}
    \end{subfigure}
    \caption{Example 1D: total CPU time.}
    \label{fig:test4_CPU_time}
\end{figure}

\begin{figure}
    \centering
    \includegraphics[width=0.48\textwidth]{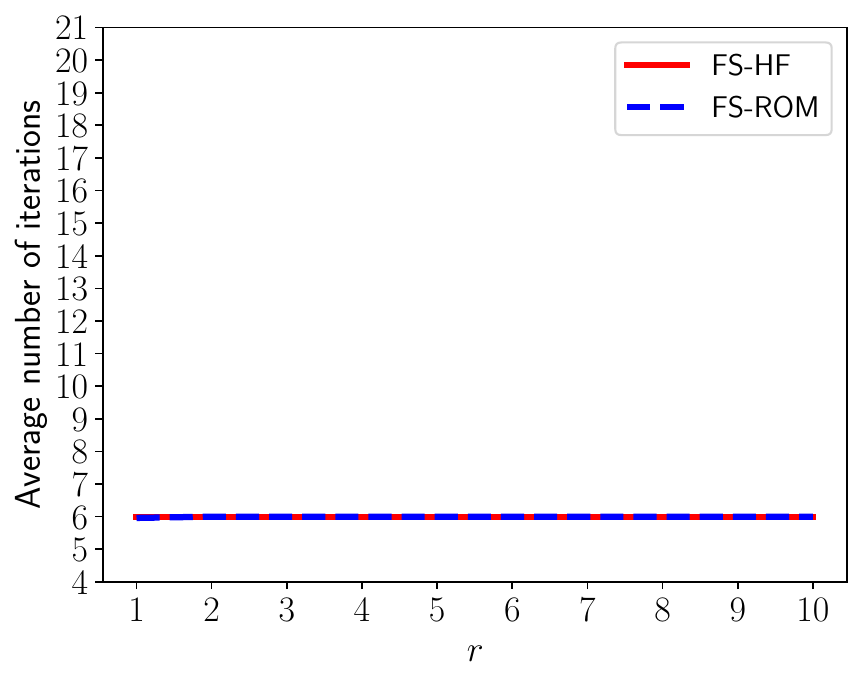}
    \caption{Example 1D: average number of fixed-stress iterations.}
    \label{fig:test4_iterations}
\end{figure}

For the last test case in Example 1, we explore the possibility of employing different time step sizes $\dt$ between the HF and the ROM solvers.
Indeed, since the evaluation of the ROM still relies on time stepping through the interval $(0, T]$,
adopting a larger time step presents a way to further enhance the speed of online evaluation as it entails a reduction in the number of time instances. Here, we set $\dt = 0.001$ for FS-HF as in Example 1B, but FS-ROM with $\dt = 0.01$, i.e.\ ten times larger than the one in Example 1B. We consider $T=1$.

Figure \ref{fig:test4_relative_error_wrt_exact_fs} displays the relative errors of the displacement, pressure, and temperature (in $\|\cdot\|_1$ and $\| \cdot \|$ norms) with respect to the analytical solutions in \eqref{eq:analytical solution} for the FS-HF and FS-ROM schemes.
Comparing Figures \ref{fig:test2_relative_error_wrt_exact_fs} and \ref{fig:test4_relative_error_wrt_exact_fs}, we conclude that employing the FS-ROM scheme with a larger time step introduces a larger time discretization error in the approximation of the pressure $p$ and the temperature $\theta$.
However, in Figure \ref{fig:test4_relative_error_wrt_exact_fs}, the errors in $\pr$ and $\thetar$ still remain in the same order of magnitude as those for $\ph$ and $\thetah$, thus offering an approximation that could be accurate enough in several practical scenarios.
Furthermore, the displacement $\bu$, which exhibits minimal variation over time in this setup, can be precisely approximated even by a ROM that employs a larger time step for time integration. In fact, there are negligible distinctions observed between $\bur$ and $\buh$.

As expected, choosing a larger time step decreases the online computational cost. Comparing Figure \ref{fig:test4_CPU_time} to Figure \ref{fig:test2_CPU_time}, we realize that the CPU time is decreased by a factor of ten since the number of time instances is divided by ten. Finally, choosing a larger time step does not have a detrimental effect on the average number of FS-ROM iterations, which are still equal to the ones of the FS-HF scheme, see Figure \ref{fig:test4_iterations}.

\FloatBarrier

\subsection{Example 2. A domain with parametric
heterogeneities}
\label{sec:appl}

The purpose of this example is to assess the capabilities of the proposed FS-ROM in a parametric setting with a more realistic application problem.
The domain $\Omega=(0,1)^2$ is divided into two subdomains, $\Omega_1$ and $\Omega_2$, with $\Omega_1$ being characterized by higher permeability and conductivity than $\Omega_2$.


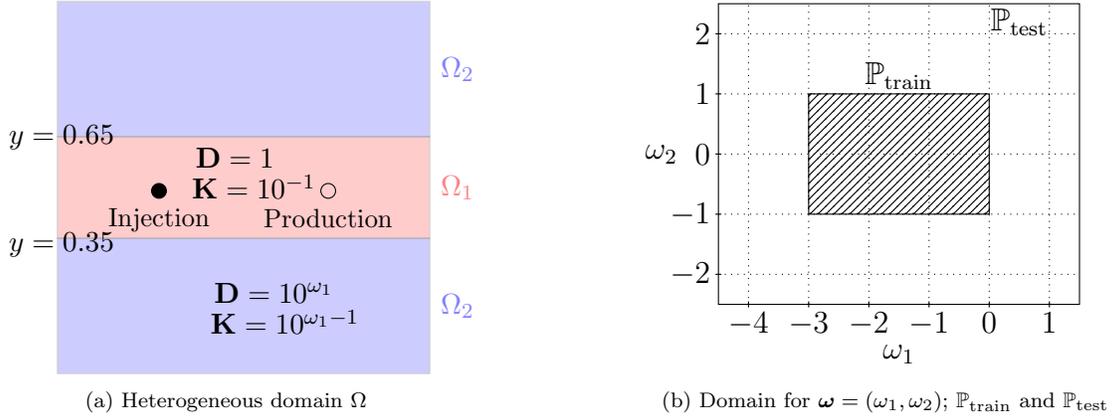
\begin{figure}[!h]
\centering
\begin{subfigure}[b]{0.4\textwidth}
\begin{tikzpicture}[scale=1.8]


\draw[fill=blue, opacity=0.2] (0,0) rectangle (2.75,1.375-0.375);
\draw[fill=red,  opacity=0.2] (0,1.375-0.375) rectangle (2.75,1.375+0.375);
\draw[fill=blue, opacity=0.2] (0,1.375+0.375) rectangle (2.75,2.75);
\draw[] (1.31,1.75) node[below] {\small $\bD = 1$};
\draw[] (1.45,1.55) node[below] {\small $\bK = 10^{-1}$};
\draw[] (1.6,0.75) node[below] {\small $\bD = 10^{\omega_1}$};
\draw[] (1.68,0.55) node[below] {\small $\bK = 10^{\omega_1-1}$};
\draw[black,fill=black] (0.75,1.35) circle (.3ex);

\draw[black] (2.,1.35) circle (.3ex);

\draw[] (0.75,1.3) node[below] {\footnotesize Injection};
\draw[] (2.,1.3) node[below] {\footnotesize Production};
\draw[opacity=0.5] (2.75,0.5) node[right] {\color{blue}\small $\Omega_2$};
\draw[opacity=0.5] (2.75,2.25) node[right] {\color{blue}\small $\Omega_2$};
\draw[opacity=0.5] (2.75,1.375) node[right] {\color{red}\small $\Omega_1$};

\draw[] (0.5,0.95) node[left] {\small
$y= 0.35$};
\draw[] (0.5,1.75) node[left] {\small
$y=0.65$};
\end{tikzpicture}
\caption{Heterogeneous domain $\Omega$}
\label{fig:ex3_setup_a}
\end{subfigure}
\hfill
\begin{subfigure}[b]{0.44\textwidth}
\begin{tikzpicture}[scale=0.8]


\draw[] (-4.5,-2.5) rectangle (1.5,2.5);
\draw[pattern=north east lines] (-3.,-1.) rectangle (0,1.);

\foreach \x in {-4,-3,...,1.} {%
    \draw ($(\x,-2.5) + (0,-0.05)$) -- ($(\x,-2.5) + (0,0.05)$)
        node [below] {$\x$};

    \draw [dotted] (\x,-2.5) -- (\x,2.5);
}

\foreach \y in {-2,-1,...,2} {%
    \draw ($(-4.5,\y) + (-0.05,0)$) -- ($(-4.5,\y) + (0.05,0)$)
        node [left] {$\y$};
    \draw [dotted] (-4.5,\y) -- (1.5,\y);
}

\draw[] (-1.5,0.9) node[above] {$\mathbb{P}_{\text{train}}$};

\draw[] (0.5,1.8) node[above] {$\mathbb{P}_{\text{test}}$};

\draw[] (-1.5,-3) node[below] {$\omega_1$};
\draw[] (-5.,0) node[left] {$\omega_2$};
\end{tikzpicture}
\caption{Domain for $\bomega = (\omega_1, \omega_2)$;  $\mathbb{P}_\text{train}$ and $\mathbb{P}_\text{test}$}
\label{fig:ex3_setup_b}
\end{subfigure}
\caption{Example 2. (a) setup of the computational domain $\Omega$, (b) domain for the pair $\bomega = (\omega_1, \omega_2)$. We note that $\mathbb{P}_\text{test}$ is larger than $\mathbb{P}_\text{train}$ to consider extrapolation of ROM. }
\label{fig:ex3_setup}
\end{figure}
The heterogeneity of the domain is due to the differences in the values of $\bD$ and $\bK$ in $\Omega_1$ and $\Omega_2$, and it will be represented by introducing a parameter $\omega_1$.
Specifically, $\bD = 1$ and $\bK = 10^{-1}$ in $\Omega_1$, while $\bD = 10^{\omega_1}$ and $\bK = 10^{\omega_1-1}$ in $\Omega_2$. See Figure \ref{fig:ex3_setup_a} for illustration.
We introduce another parameter $\omega_2$ for the order of magnitude of the elastic properties of the material in $\Omega$, where $\lambda = \mu = 10^{\omega_2}$.
The remaining physical properties are independent of  $\omega_1$ or $\omega_2$ and  fixed as $\alpha = 1$,
 $\alpha_T = \num{e-4}$,  $\alpha_M = \num{e-6}$,
$c_0 = \num{e-2}$, $C_d = \num{1} $,  and
$\theta_0 = 1$.
The injection and production of the pressure and temperature are realized by the functions defined by
$$
g(x,y) = \eta(x,y) =10^{-2} \left( e^{(-1000 (x-x_1)^2 - 1000(y-y_1)^2)} - e^{(-1000 (x-x_2)^2 - 1000(y-y_2)^2)}\right),
$$
where $(x_1,y_1) = (0.25, 0.5)$ and $(x_2,y_2) = (0.75, 0.5)$ \cite{YiLee23}. The body force $\mathbf{f}$ and the initial conditions for the displacement, pressure, and temperature are all set to zero.
As for the boundary conditions, we impose homogeneous Dirichlet boundary conditions for the displacement and homogeneous Neumann boundary conditions for the pressure and temperature on the entire boundary, $\partial \Omega$. Additionally, we choose a mesh size  $h=0.01$, {linear finite elements}, and a time step size  $\Delta t  = 0.1$ with the final time $T = 2$, and the stabilization parameter is again set to $L=1$.

\begin{figure}
    \centering
    \begin{subfigure}[b]{0.48\textwidth}
        \centering
        \includegraphics[width=\textwidth]{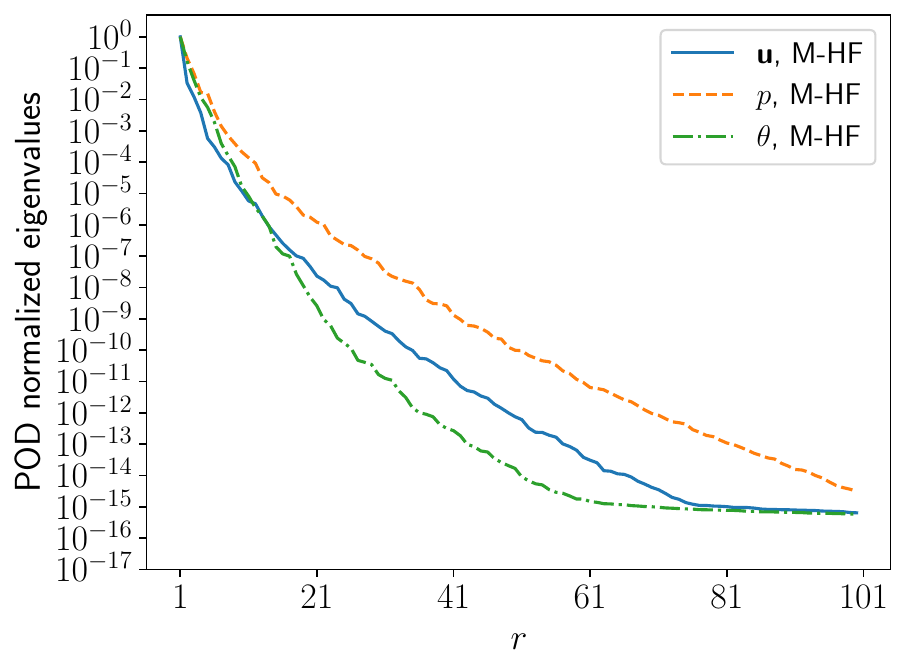}
        \caption{Monolithic schemes (M-HF to train M-ROM).}
        \label{fig:application_POD_eigs_cg}
    \end{subfigure}
    \hfill
    \begin{subfigure}[b]{0.48\textwidth}
        \centering
        \includegraphics[width=\textwidth]{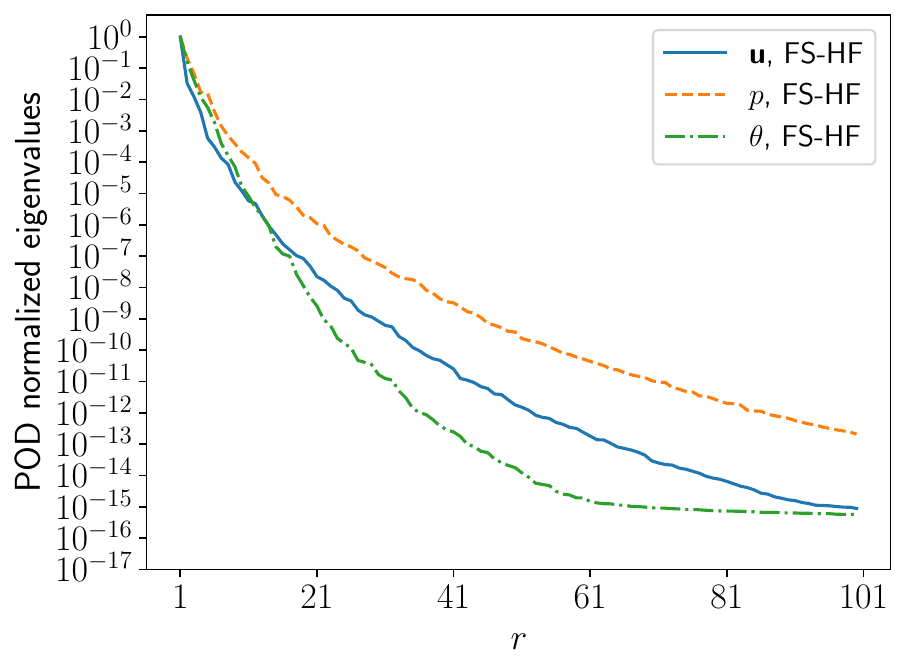}
        \caption{Iterative schemes (FS-HF to train FS-ROM).}
        \label{fig:application_POD_eigs_fs}
    \end{subfigure}
    \caption{Example 2. POD normalized eigenvalues $\nu_k / \nu_0$.}
    \label{fig:application_POD_eigs}
\end{figure}

Let $\bomega := (\omega_1, \omega_2)$.
During the training of the ROM, the parameter pair $\bomega = (\omega_1, \omega_2)$ is varied in a given training set  $\mathbb{P}_{\text{train}} = [-3, 0] \times [-1, 1]$.
Then, to assess the extrapolation capabilities of
M-ROM and FS-ROM, the parameter $\bomega$ is allowed to have a larger variation during the evaluation of the ROM. More specifically, $\bomega$ is chosen from the test parameter space $\mathbb{P}_{\text{test}} = [-4, 1] \times [-2, 2]$, which is larger than $\mathbb{P}_{\text{train}}$. See Figure \ref{fig:ex3_setup_b} for illustration.

On the other hand, we let iterations continue until the relative increment in the $\normo{\cdot}$-norm becomes less than $\epsilon = 10^{-3}$ or it reaches a maximum of $20$ iterations, whichever is met first. We note that the value of $\epsilon$ is considerably larger than the one employed in the Examples 1A-1D. Even with this larger $\epsilon$, however, the FS-HF scheme converges for the most challenging parametric configuration (i.e., $\bomega = (-3, -1))$ within the maximum allowed number of iterations.
Since FS-HF will need to be queried for several values of $\bomega$, by requiring a tolerance $\epsilon$ as tight as in the previous validation (and, thus, increasing the maximum allowed number of iterations) we would make the training process extremely expensive.

We then define a discrete training set by sampling 25 equidistant points on a $5 \times 5$ uniform grid on $\mathbb{P}_{\text{train}}$ and proceed to train both M-ROM and FS-ROM.
The inclusion of the parameter $\bomega$ in this example necessitates a slight modification to the POD method presented in Section \ref{sec: POD}. This alteration is introduced to facilitate the applicability of the ROM for parameter values that may deviate from those employed during the training phase. Consequently, the construction of the reduced basis spaces $\Vr$, $\Wr$, and $\Thetar$ involves compressing the parametric and temporal evolution of the M-HF and FS-HF solutions at the same time. Notably, since each parametric evaluation entails 21 time steps, this compression encompasses a sequence of $25 \times 21 = 525$ M-HF and FS-HF solutions.

Figure \ref{fig:application_POD_eigs} illustrates the POD normalized eigenvalues resulting from the training of the M-ROM and FS-ROM schemes.
Two main differences compared to Examples 1A-1D are observed here.
First, the decay of the POD eigenvalues is considerably slower than those observed in Examples 1A-1D due to the higher complexity of this application. Still, the $60$-th, $80$-th and $100$-th POD normalized eigenvalue are less than $10^{-12}$ for $\bu$, $p$ and $\theta$, respectively. We thus set the maximum value of $r = 100$. Though $r$ is larger than the one used previously, $100$ ROM degrees of freedom are still a very large reduction compared to the $20402$, $10201$, and $10201$ HF degrees of freedom required by $\Vh$, $\Wh$ and $\Thetah$, respectively.
Secondly, comparing Figure \ref{fig:application_POD_eigs_cg} and Figure \ref{fig:application_POD_eigs_fs} we note that the eigenvalue decay obtained by applying POD from M-HF data and FS-HF data, respectively, is slightly different, with the pressure presenting the most noticeable difference. This is due to the relatively large value of $\epsilon$ employed in the FS-HF, due to which FS-HF solutions are not as accurate as M-HF ones.

\begin{figure}
    \centering
    \begin{subfigure}[b]{\textwidth}
        \centering
        \rotatebox{90}{\hspace{1.0cm}M-ROM}
        \includegraphics[width=0.235\textwidth, trim={0 1.25cm 4cm 0},clip]{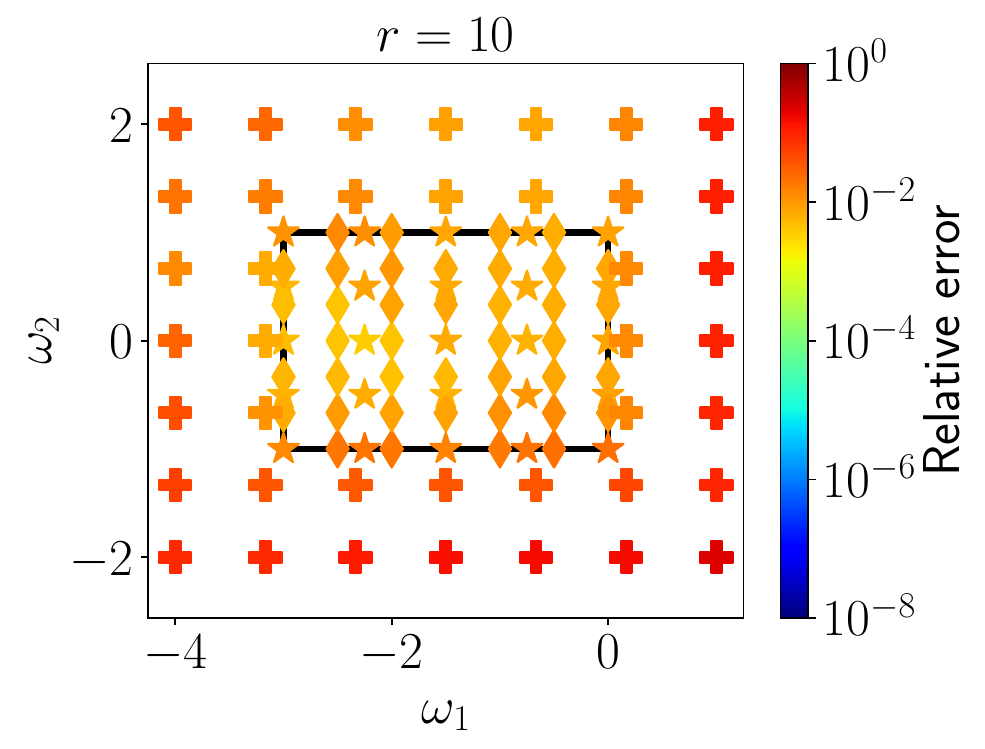}
        \includegraphics[width=0.212\textwidth, trim={1.25cm 1.25cm 4cm 0},clip]{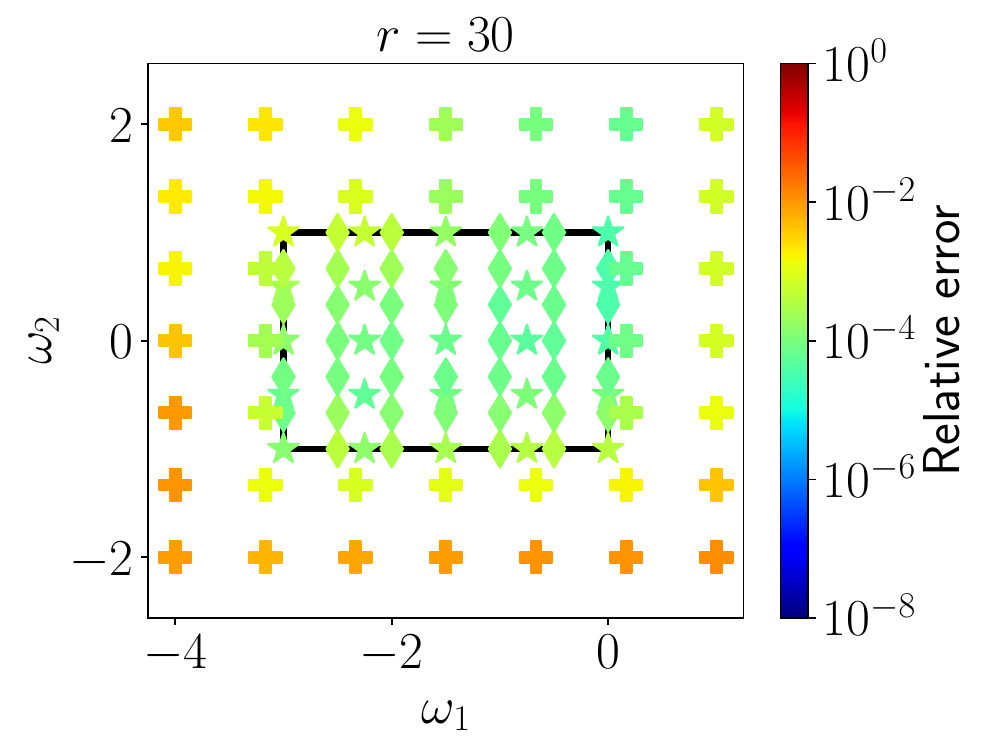}
        \includegraphics[width=0.212\textwidth, trim={1.25cm 1.25cm 4cm 0},clip]{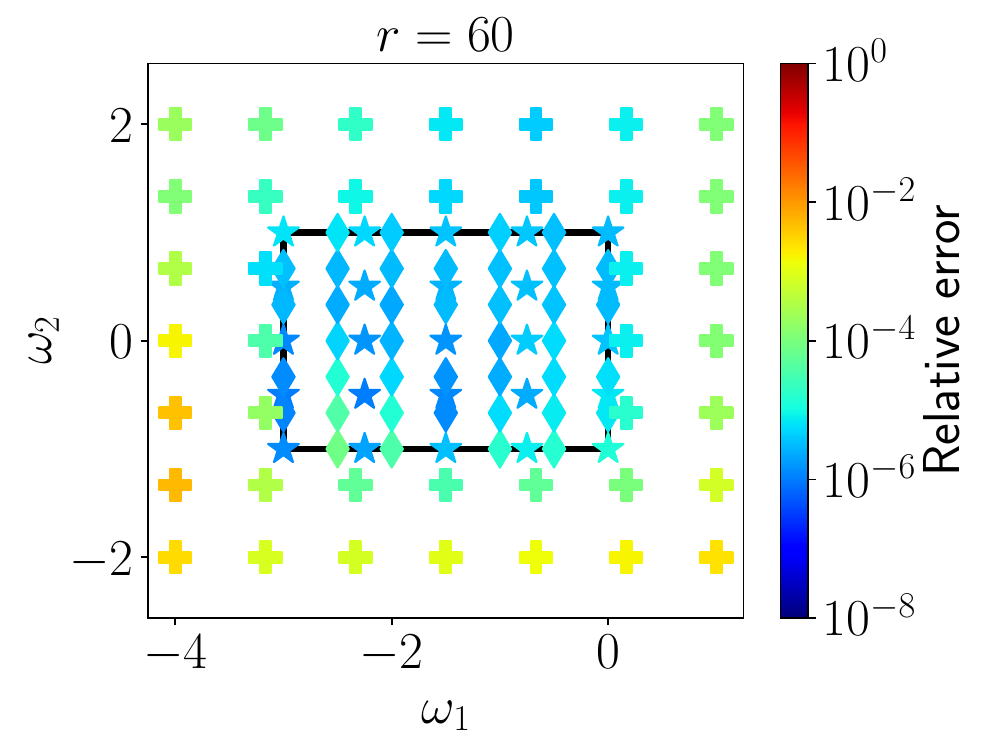}
        \includegraphics[width=0.287\textwidth, trim={1.25cm 1.25cm 0 0},clip]{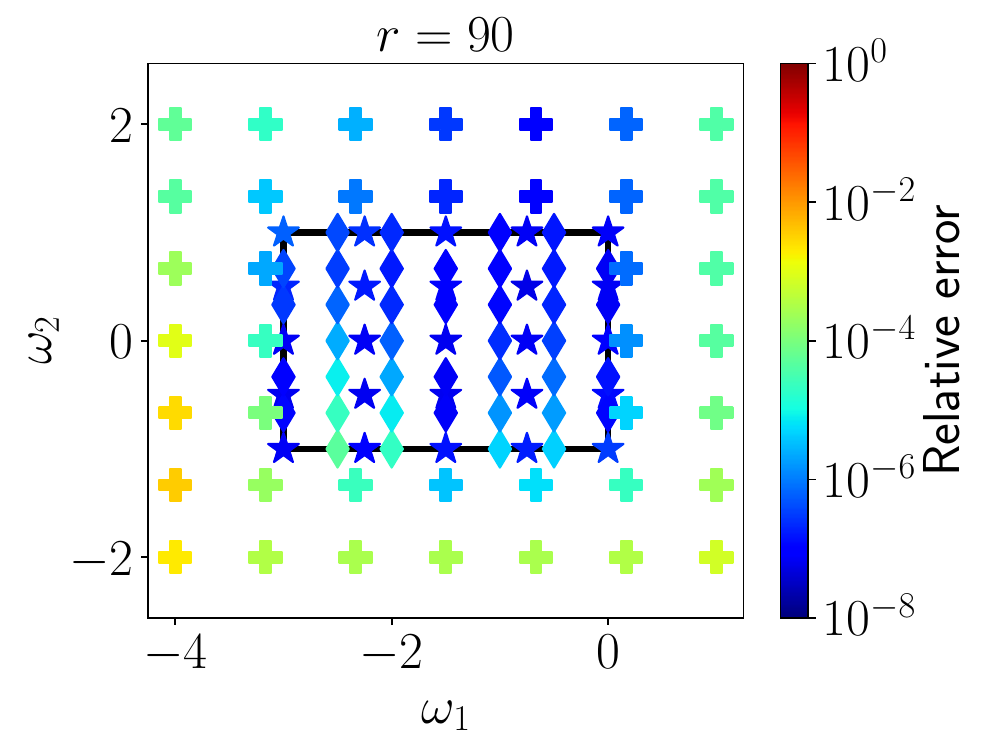}
        \\
        \rotatebox{90}{\hspace{1.1cm}FS-ROM}
        \includegraphics[width=0.235\textwidth, trim={0 0.25cm 4cm 0.15cm},clip]{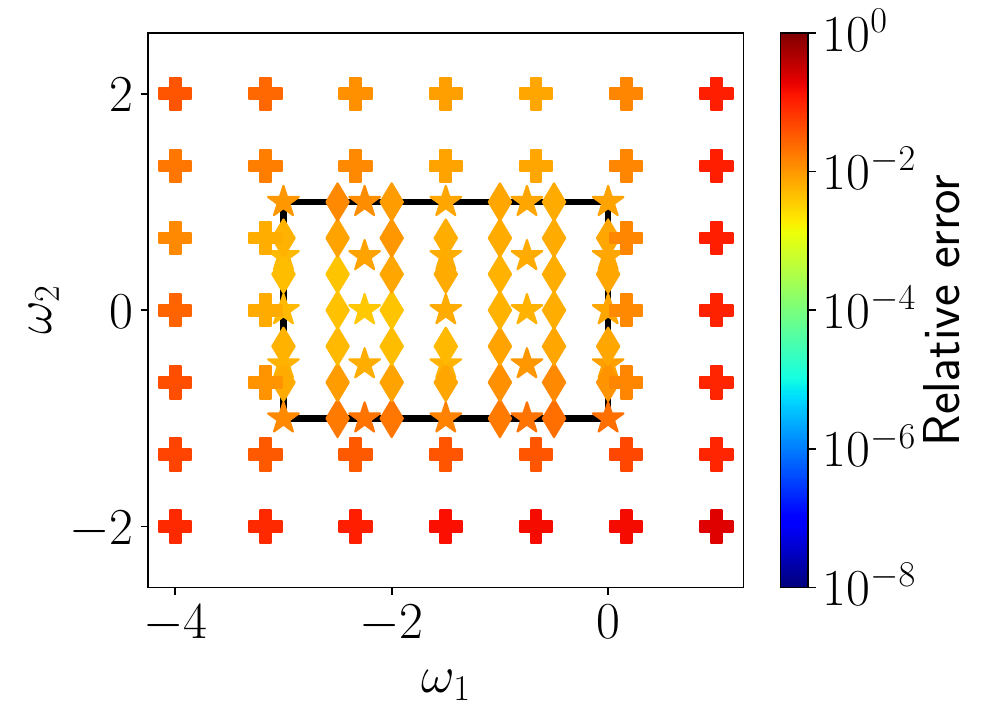}
        \includegraphics[width=0.212\textwidth, trim={1.25cm 0.25cm 4cm 0.15cm},clip]{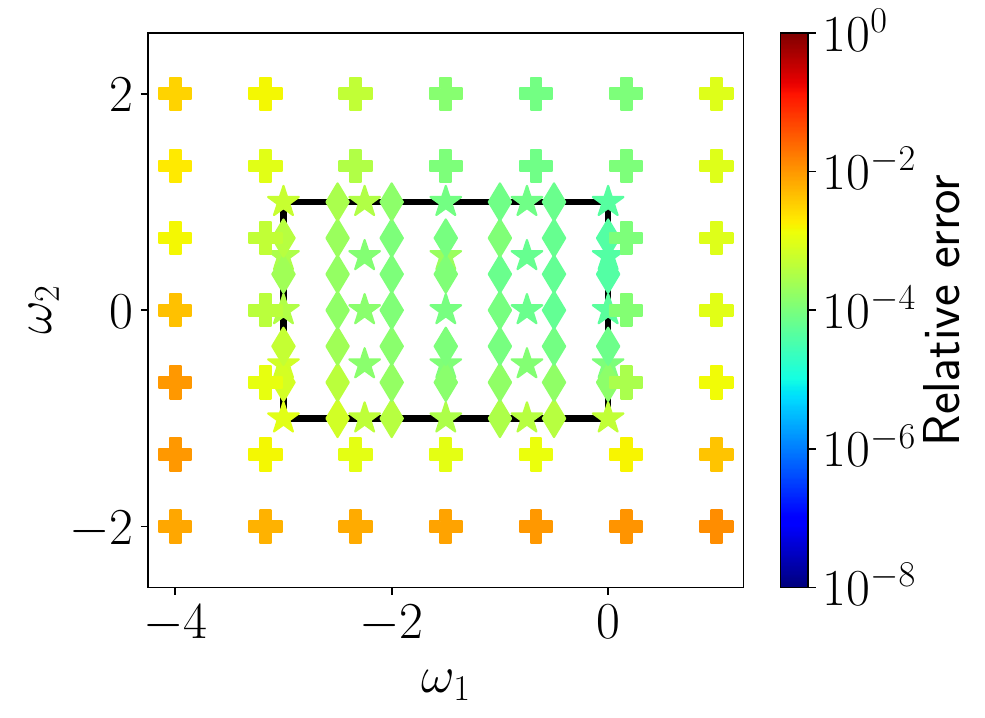}
        \includegraphics[width=0.212\textwidth, trim={1.25cm 0.25cm 4cm 0.15cm},clip]{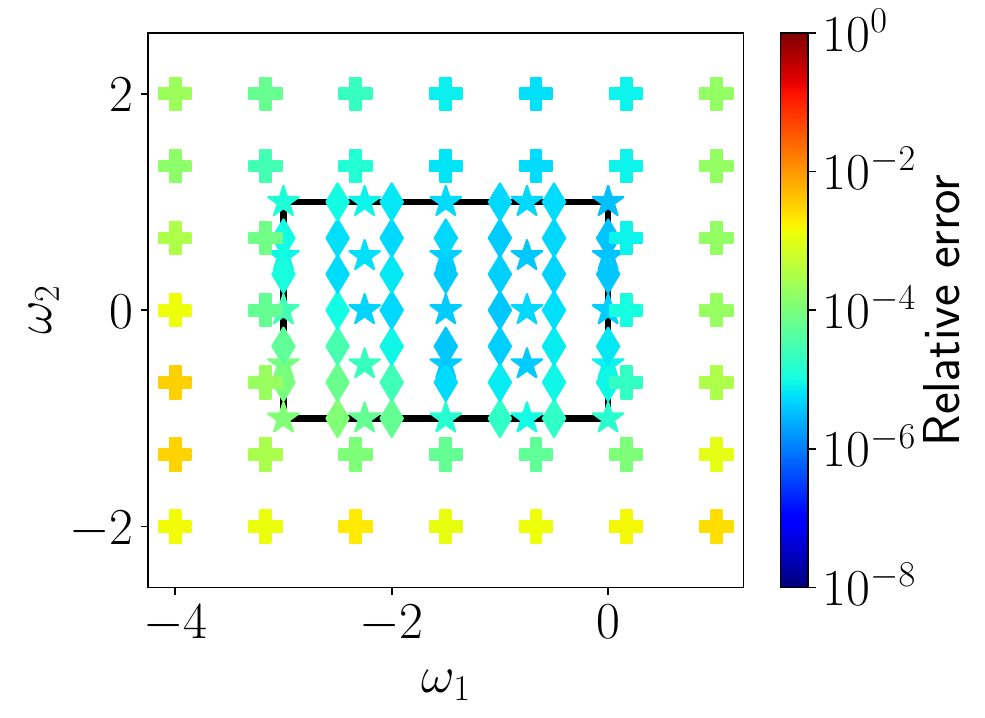}
        \includegraphics[width=0.287\textwidth, trim={1.25cm 0.25cm 0 0.15cm},clip]{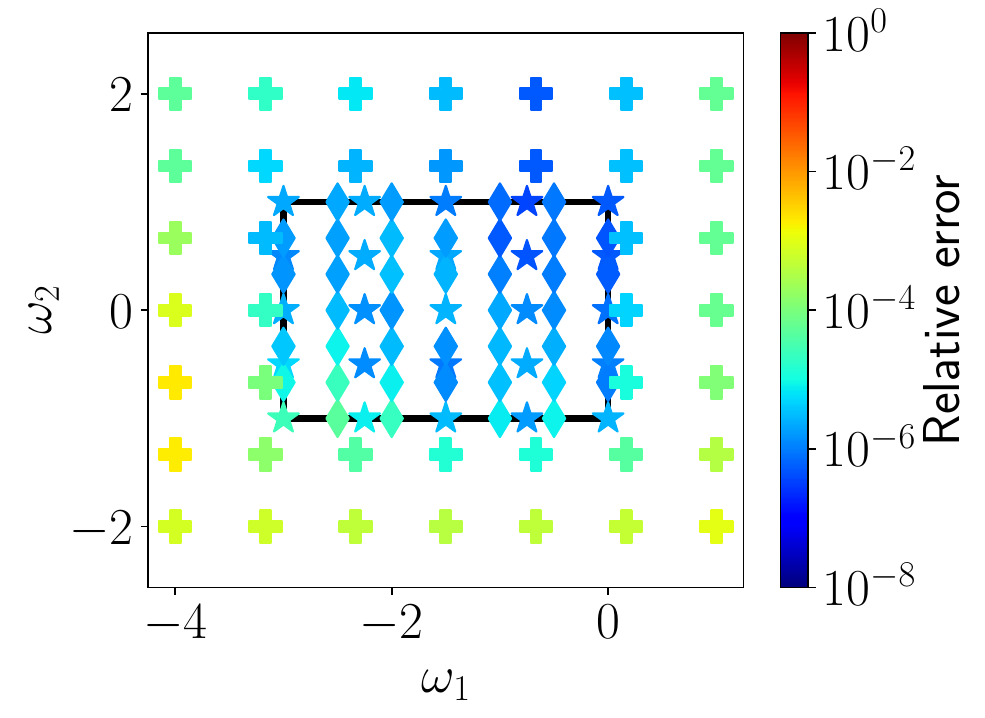}
        \caption{$\normo{\bur - \buh} / \normo{\buh}$.}
        \label{fig:application_u_H1_relative_errors}
    \end{subfigure}\\
    \begin{subfigure}[b]{\textwidth}
        \centering
        \rotatebox{90}{\hspace{1.0cm}M-ROM}
        \includegraphics[width=0.235\textwidth, trim={0 1.25cm 4cm 0},clip]{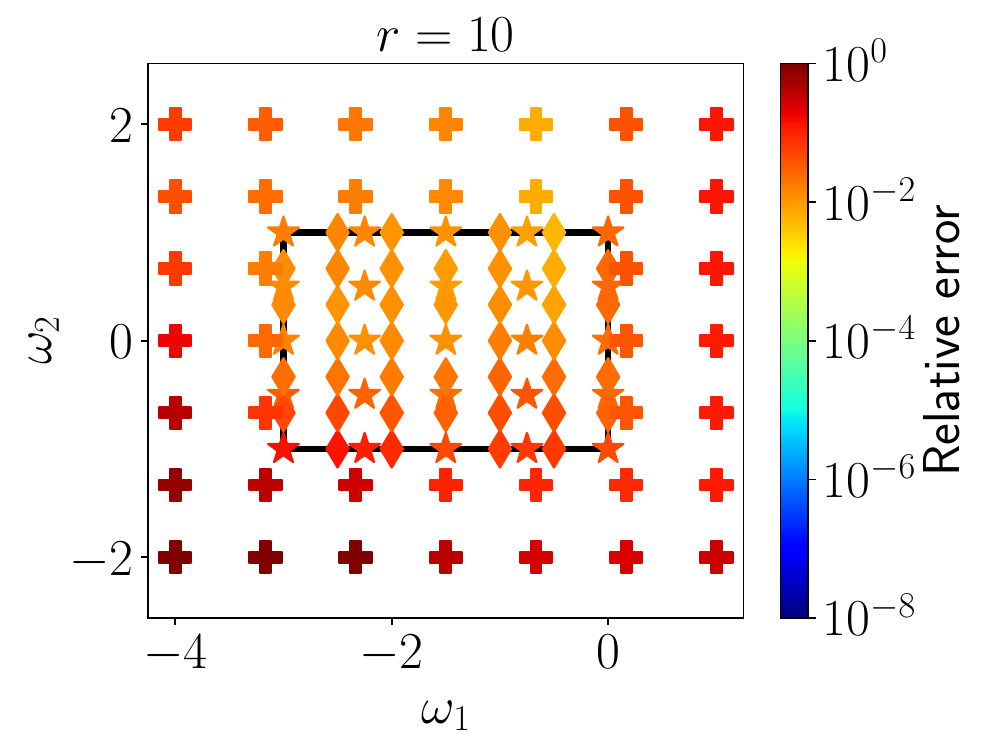}
        \includegraphics[width=0.212\textwidth, trim={1.25cm 1.25cm 4cm 0},clip]{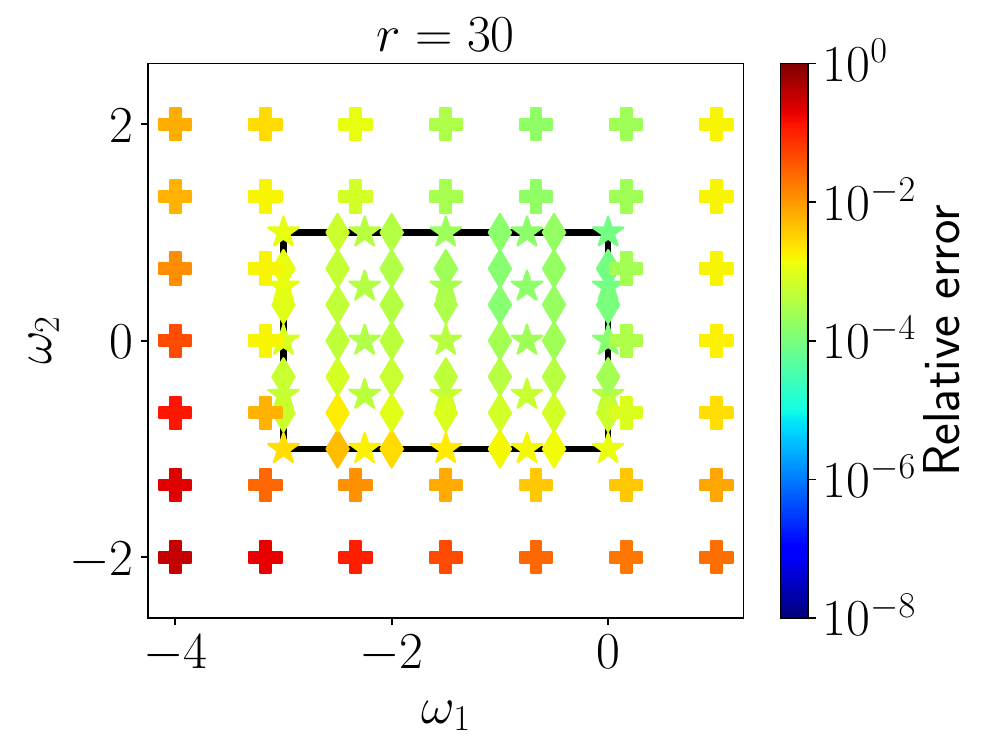}
        \includegraphics[width=0.212\textwidth, trim={1.25cm 1.25cm 4cm 0},clip]{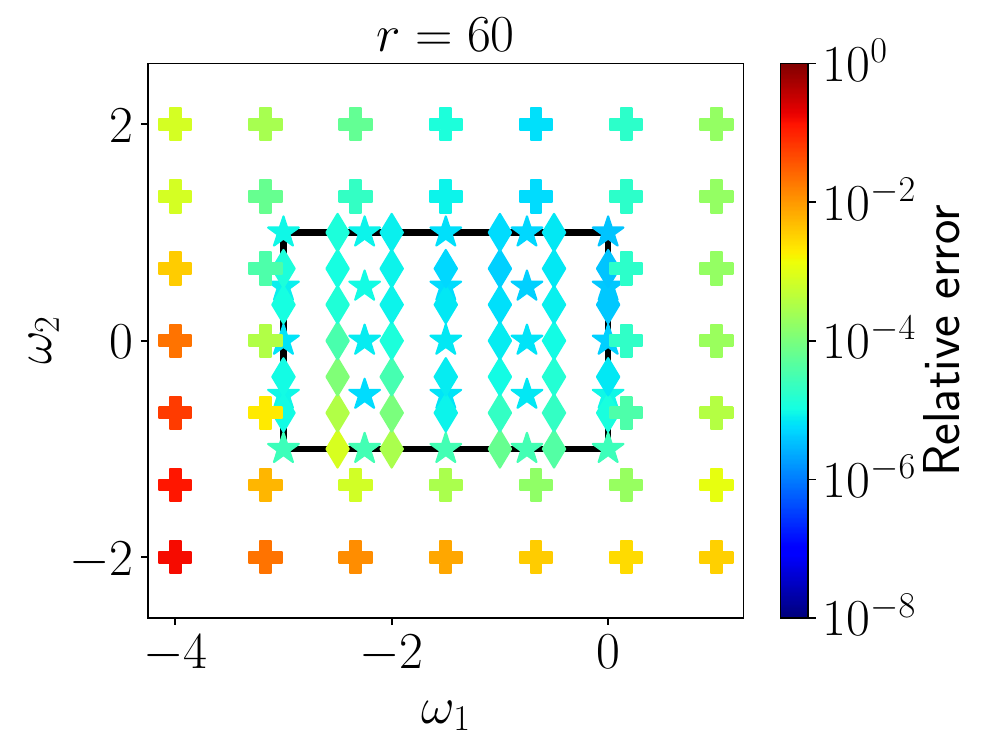}
        \includegraphics[width=0.287\textwidth, trim={1.25cm 1.25cm 0 0},clip]{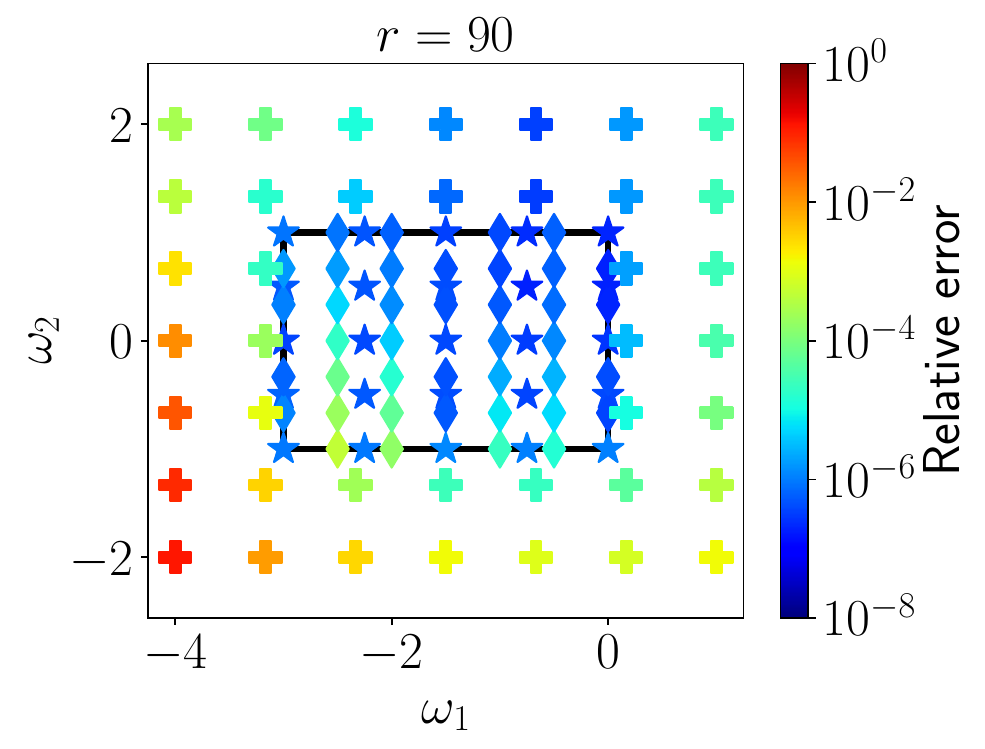}
        \\
        \rotatebox{90}{\hspace{1.1cm}FS-ROM}
        \includegraphics[width=0.235\textwidth, trim={0 0.25cm 4cm 0.15cm},clip]{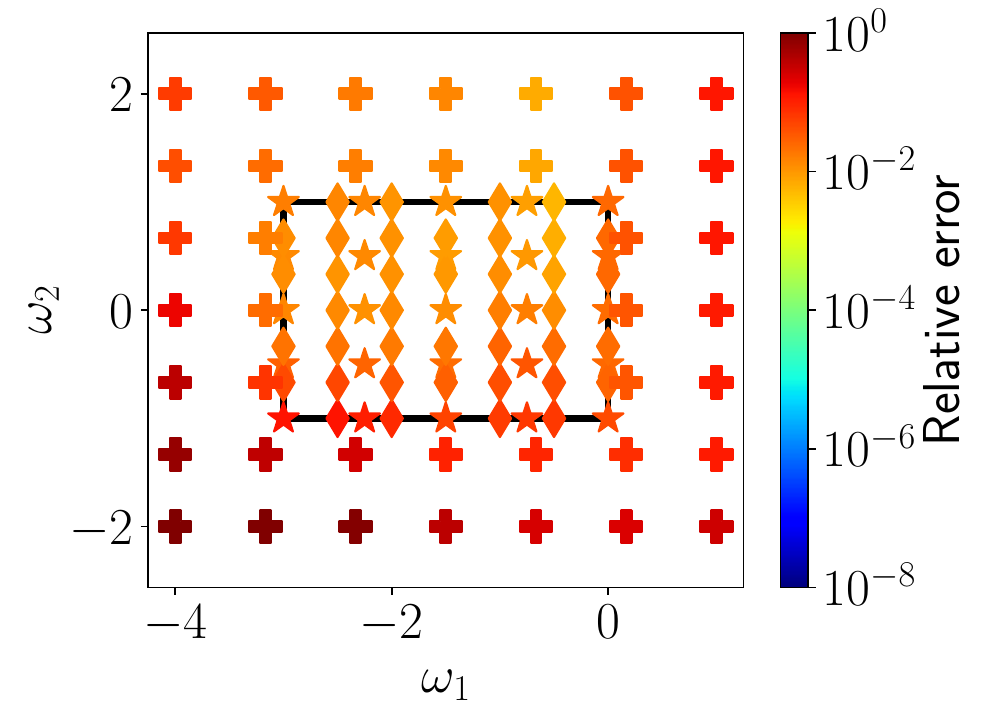}
        \includegraphics[width=0.212\textwidth, trim={1.25cm 0.25cm 4cm 0.15cm},clip]{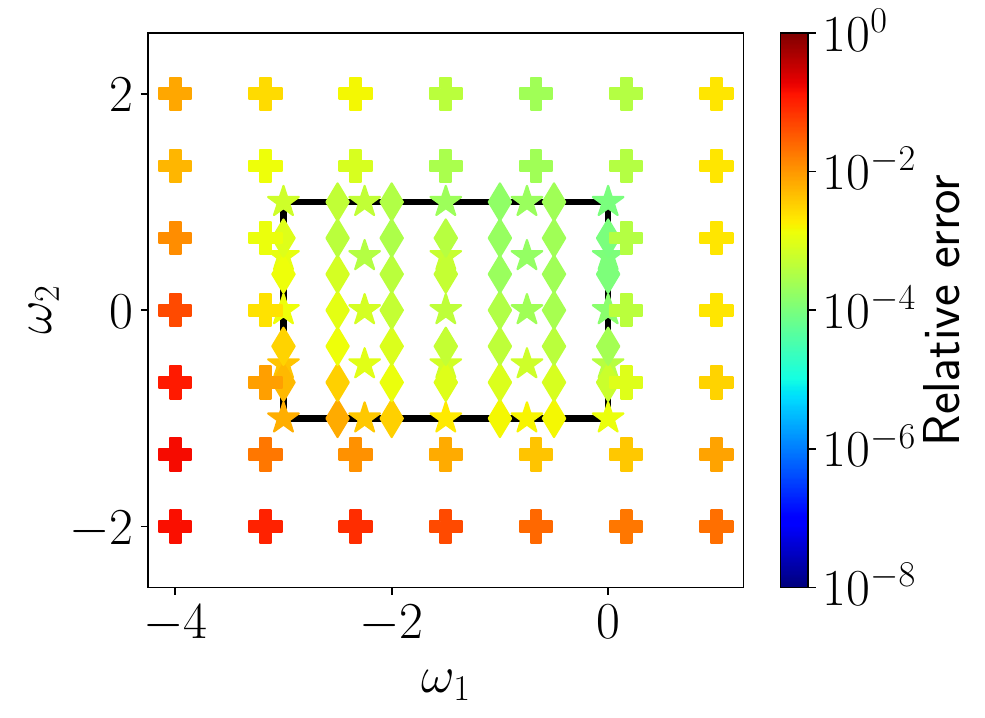}
        \includegraphics[width=0.212\textwidth, trim={1.25cm 0.25cm 4cm 0.15cm},clip]{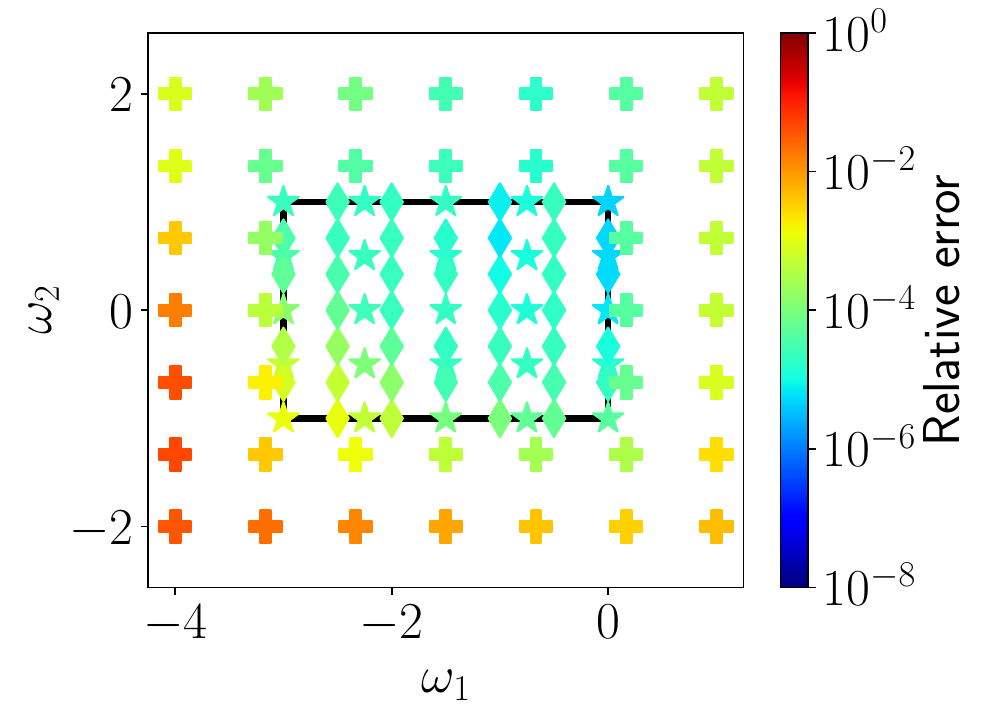}
        \includegraphics[width=0.287\textwidth, trim={1.25cm 0.25cm 0 0.15cm},clip]{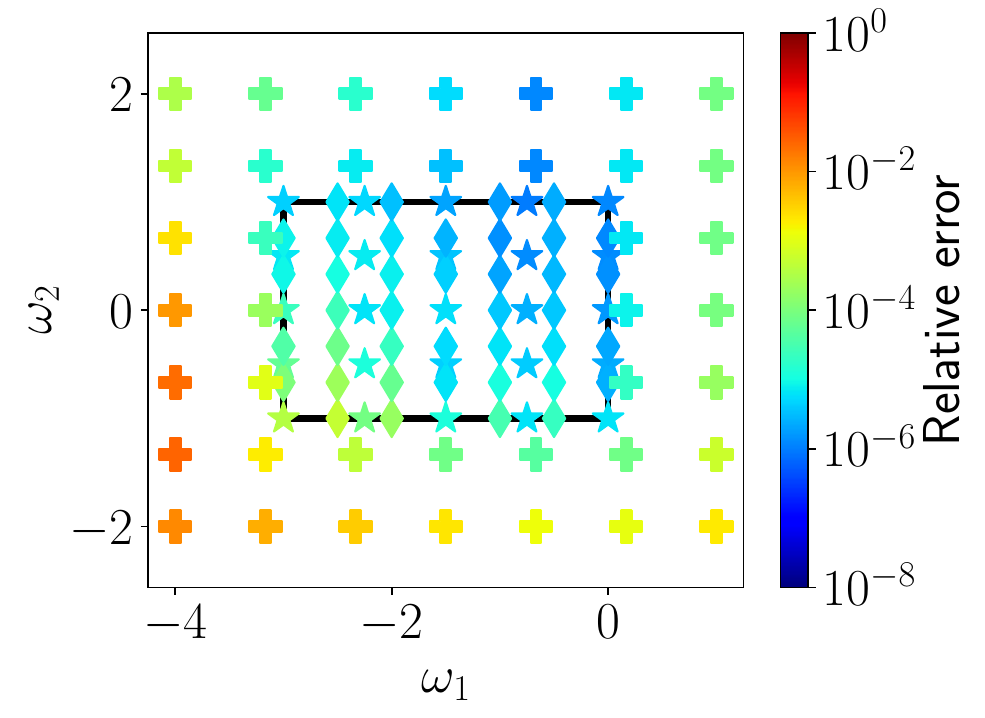}
        \caption{$\normo{\pr - \ph} / \normo{\ph}$.}
        \label{fig:application_p_H1_relative_errors}
    \end{subfigure}\\
        \begin{subfigure}[b]{\textwidth}
        \centering
        \rotatebox{90}{\hspace{1.0cm}M-ROM}
        \includegraphics[width=0.235\textwidth, trim={0 1.25cm 4cm 0},clip]{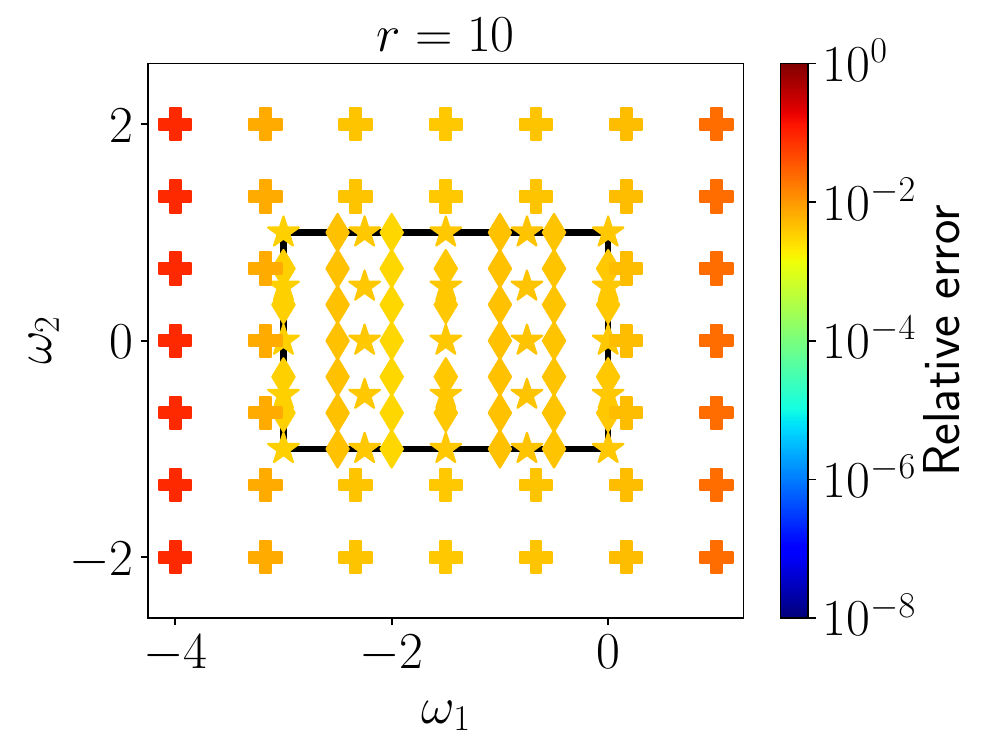}
        \includegraphics[width=0.212\textwidth, trim={1.25cm 1.25cm 4cm 0},clip]{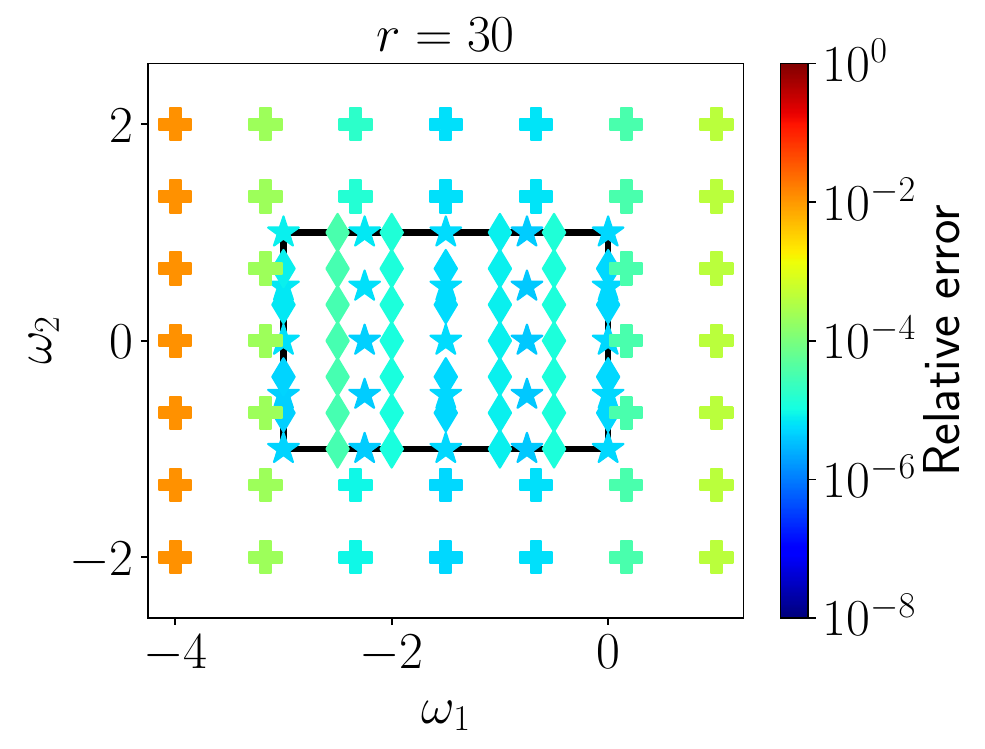}
        \includegraphics[width=0.212\textwidth, trim={1.25cm 1.25cm 4cm 0},clip]{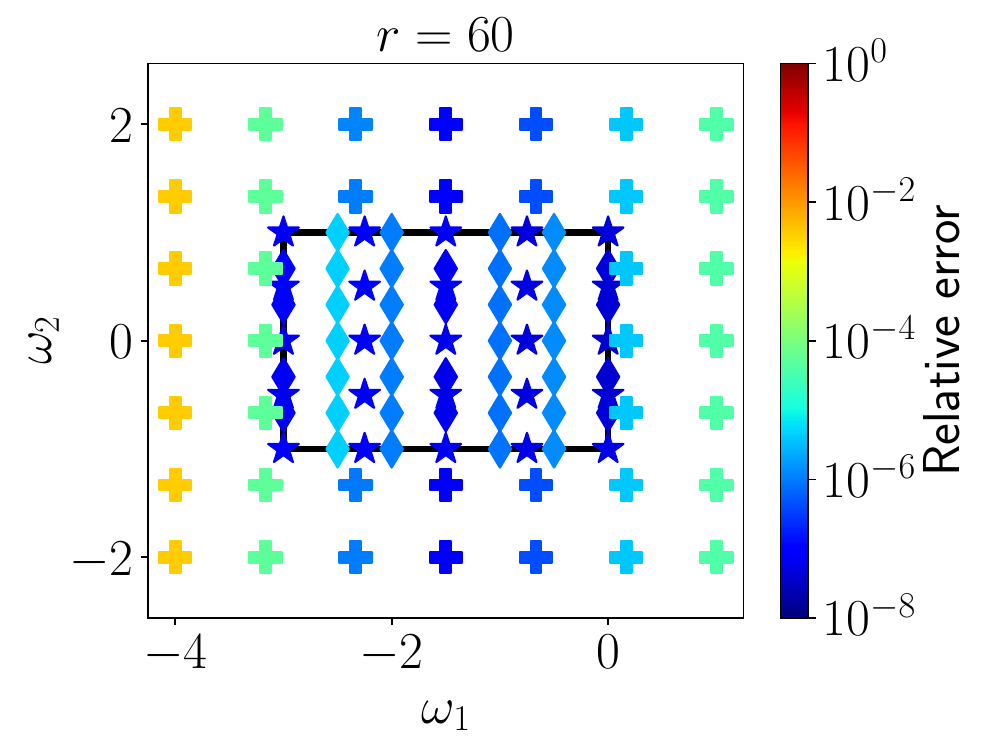}
        \includegraphics[width=0.287\textwidth, trim={1.25cm 1.25cm 0 0},clip]{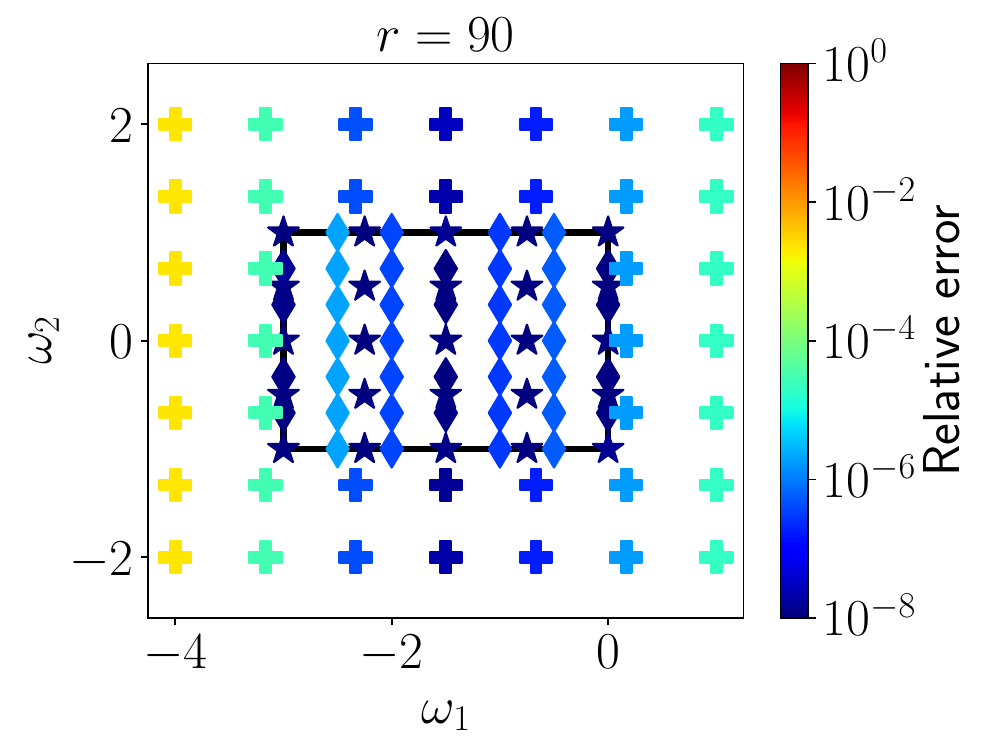}
        \\
        \rotatebox{90}{\hspace{1.1cm}FS-ROM}
        \includegraphics[width=0.235\textwidth, trim={0 0.25cm 4cm 0.15cm},clip]{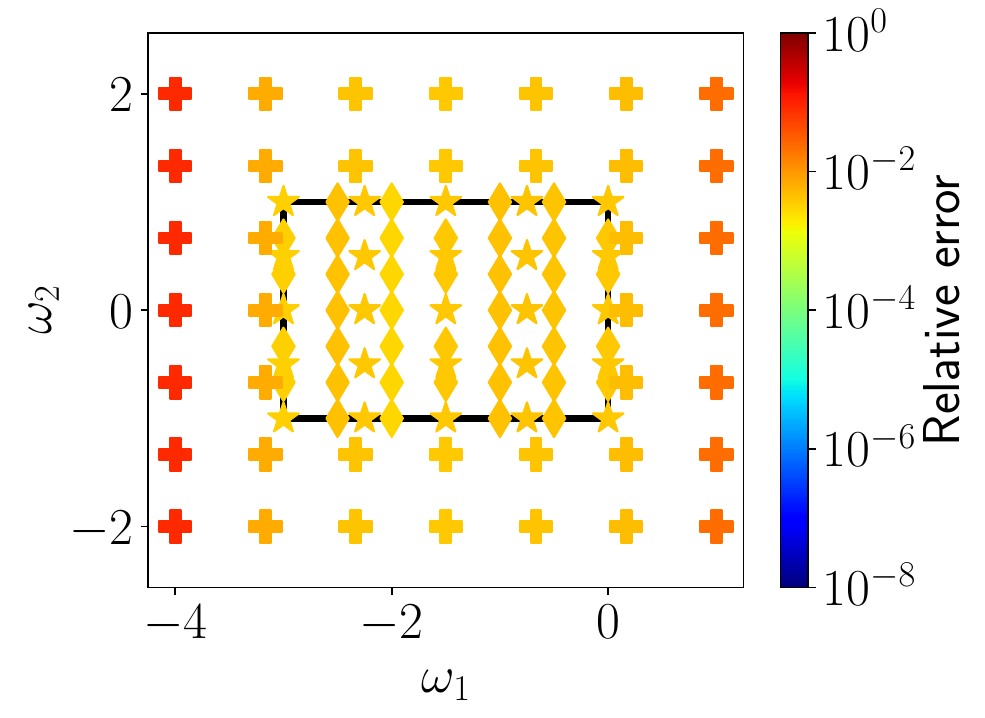}
        \includegraphics[width=0.212\textwidth, trim={1.25cm 0.25cm 4cm 0.15cm},clip]{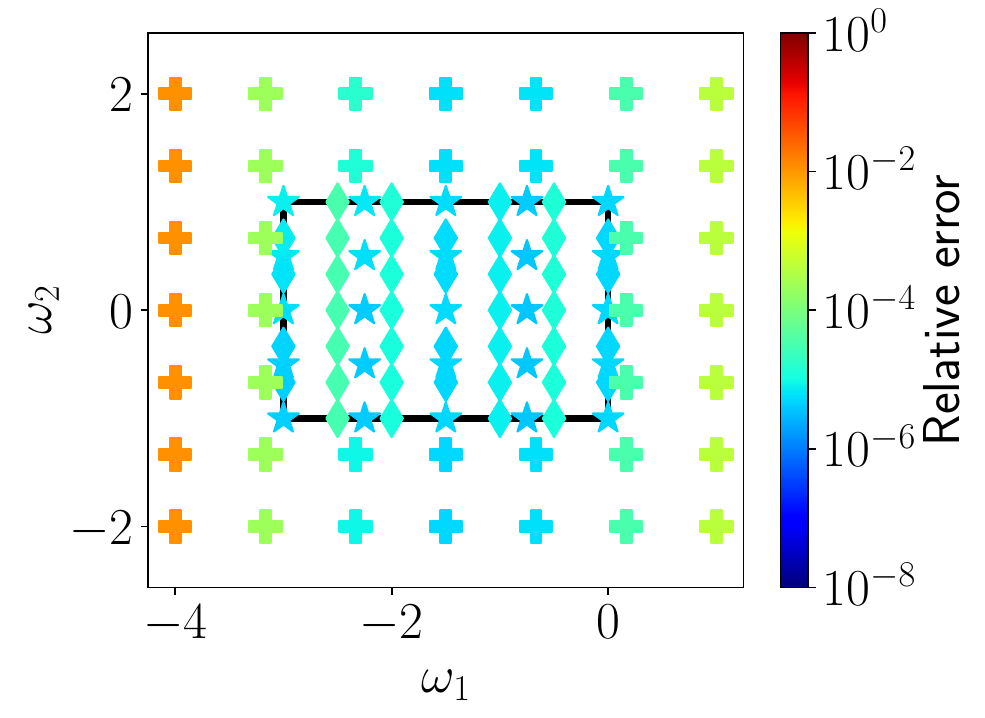}
        \includegraphics[width=0.212\textwidth, trim={1.25cm 0.25cm 4cm 0.15cm},clip]{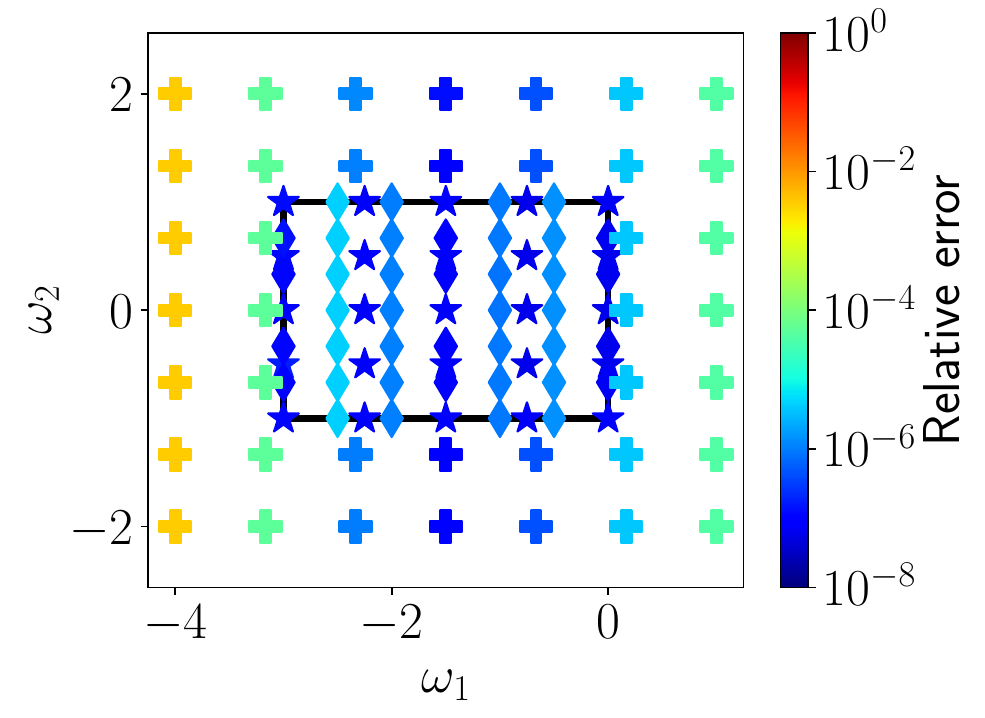}
        \includegraphics[width=0.287\textwidth, trim={1.25cm 0.25cm 0 0.15cm},clip]{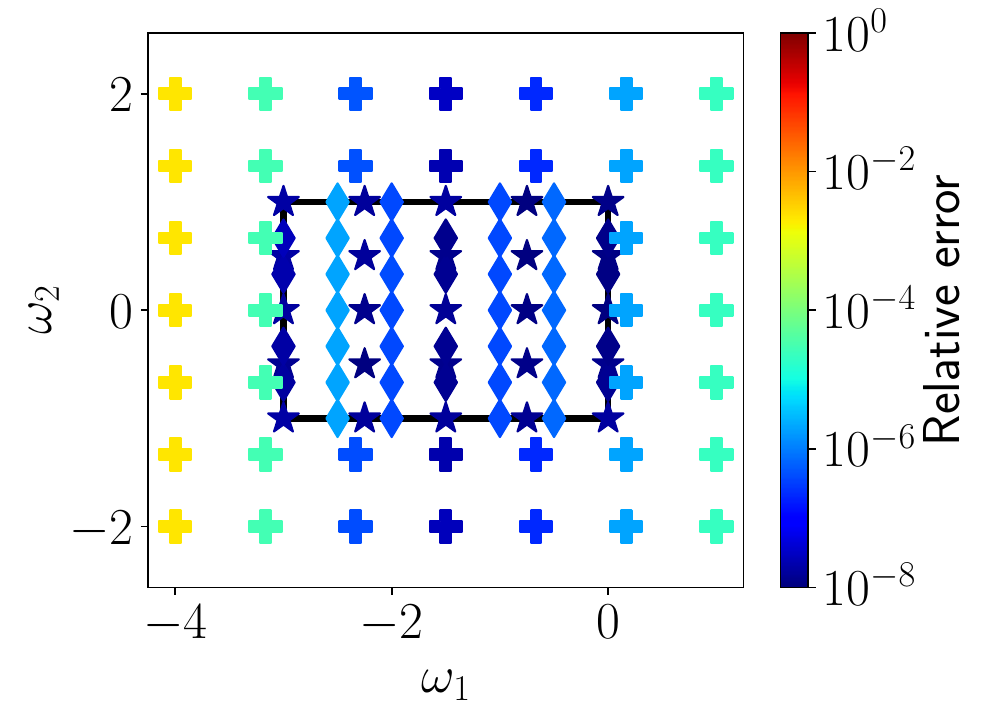}
        \caption{$\normo{\thetar - \thetah} / \normo{\thetah}$.}
        \label{fig:application_theta_H1_relative_errors}
    \end{subfigure}
    \caption{Example 2. Illustration of $\normo{\cdot}$-norm relative errors of ROM to HF for (a) $\bu$, (b) $p$ and (c) $\theta$ for the given $\mathbb{P}_{\text{test}}$ ($\mbox{\fontsize{10pt}{\baselineskip}\selectfont\FiveStar}$ for case i), $\blacklozenge$ for case ii), \ding{58} for case iii)). The inner black rectangle denotes $\mathbb{P}_{\text{train}}$ as shown in Figure \ref{fig:ex3_setup_b}.}
\label{fig:application}
\end{figure}

Next, Figure \ref{fig:application} demonstrates the assessment of the ROM through $\normo{\cdot}$-norm relative errors. This evaluation is conducted on three distinct test cases, each representing progressively more demanding conditions for the ROM.
\begin{itemize}
    \item Case i) The first testing set is composed of the same 25 parameter values used while training the ROM. These cases  are represented by star-shaped markers ($\mbox{\fontsize{10pt}{\baselineskip}\selectfont\FiveStar}$) in Figure \ref{fig:application}.
    \item Case ii) The second testing set is obtained by taking a finer $7 \times 7$ uniform grid (i.e., 49 parameters) on $\mathbb{P}_{\text{train}}$ and discarding any parameter values already present in the training set (e.g., a parameter $\bomega = (-3, -1)$), resulting in 40 parameter values. This set is used to evaluate the accuracy of the ROM on parameters that were not seen during training but still belong to the training range. The results are depicted in Figure \ref{fig:application} by diamond-shaped markers ($\blacklozenge$).
    \item Case iii) Finally, the third testing set amounts to sampling a $7 \times 7$ uniform grid on $\mathbb{P}_{\text{test}}$ and discarding parameter values lying in $\mathbb{P}_{\text{train}}$, for a total of 40 parameter configurations in $\mathbb{P}_{\text{test}} \setminus \mathbb{P}_{\text{train}}$. These values are employed to assess the accuracy of the ROM while extrapolating outside of the training range and are reported with plus-shaped markers (\ding{58}) in Figure \ref{fig:application}.
\end{itemize}

The evaluation of M-ROM and FS-ROM follows through as described in Section~\ref{sec: rom}, with a minor modification regarding the assembly of the ROM operators in Table~\ref{table:rom_operators}.
Here, we consider only
$[\mat{A}_r^{\theta \theta}]_{\beta\gamma} = (\bD \nabla \varphi_\gamma^{\theta} , \nabla \varphi_\beta^{\theta})$ for the sake of exposition. Since now $\bD = \bD(\bomega)$, the resulting $\mat{A}_r^{\theta \theta}(\bomega)$ cannot be preassembled as reported in Table~\ref{table:rom_operators} without knowing the value of $\bomega$ a priori. To retain the computational efficiency, the ROM makes use of the so-called \textit{affine parametric dependence} or \textit{parameter separability} \cite{hesthaven2016certified,quarteroni2015reduced} and rewrites $[\mat{A}_r^{\theta \theta}]_{\beta\gamma}$ as
\begin{align*}
[\mat{A}_r^{\theta \theta}]_{\beta\gamma}(\bomega) &= (\bD(\bomega) \nabla \varphi_\gamma^{\theta} , \varphi_\beta^{\theta}) =
(\bD(\bomega) \nabla \varphi_\gamma^{\theta} , \nabla \varphi_\beta^{\theta})_{\Omega_1} +
(\bD(\bomega) \nabla \varphi_\gamma^{\theta} , \nabla \varphi_\beta^{\theta})_{\Omega_2}\\
& =
1 \ (\nabla \varphi_\gamma^{\theta} , \nabla \varphi_\beta^{\theta})_{\Omega_1} +
10^{\omega_1} \ (\nabla \varphi_\gamma^{\theta} , \nabla \varphi_\beta^{\theta})_{\Omega_2} :=
[\mat{A}_r^{\theta \theta, (1)}]_{\beta\gamma} + 10^{\omega_1} [\mat{A}_r^{\theta \theta, (2)}]_{\beta\gamma}.
\end{align*}
The matrices $\mat{A}_r^{\theta \theta, (1)}$ and $\mat{A}_r^{\theta \theta, (2)}$ are now independent of $\bomega$ and can be preassembled at the end of the offline stage.

\begin{figure}
    \begin{subfigure}[t]{\textwidth}
        \centering
        \rotatebox{90}{\hspace{1.1cm}FS-ROM}
        \includegraphics[width=0.23\textwidth, trim={0 0.25cm 5cm 0},clip]{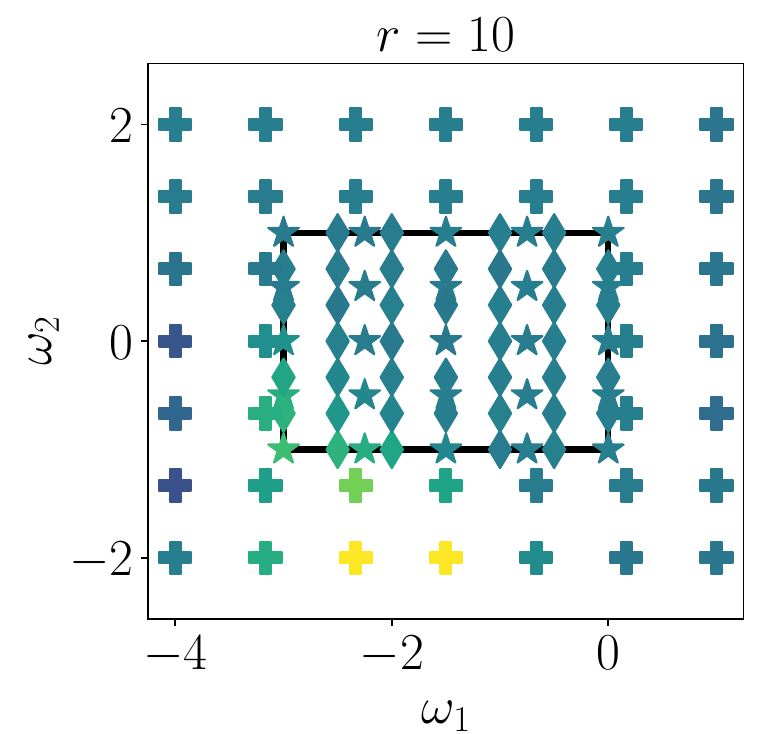}
        \includegraphics[width=0.208\textwidth, trim={1.25cm 0.25cm 5cm 0},clip]{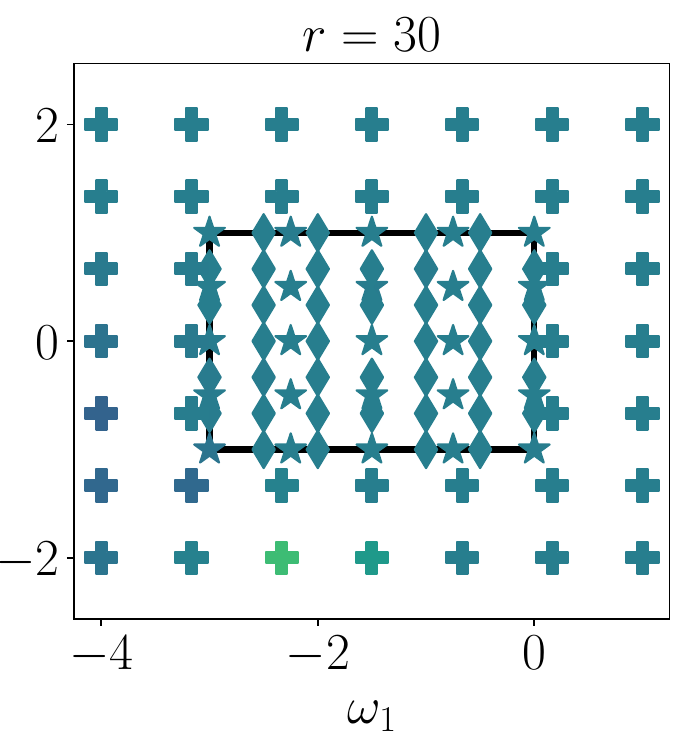}
        \includegraphics[width=0.208\textwidth, trim={1.25cm 0.25cm 5cm 0},clip]{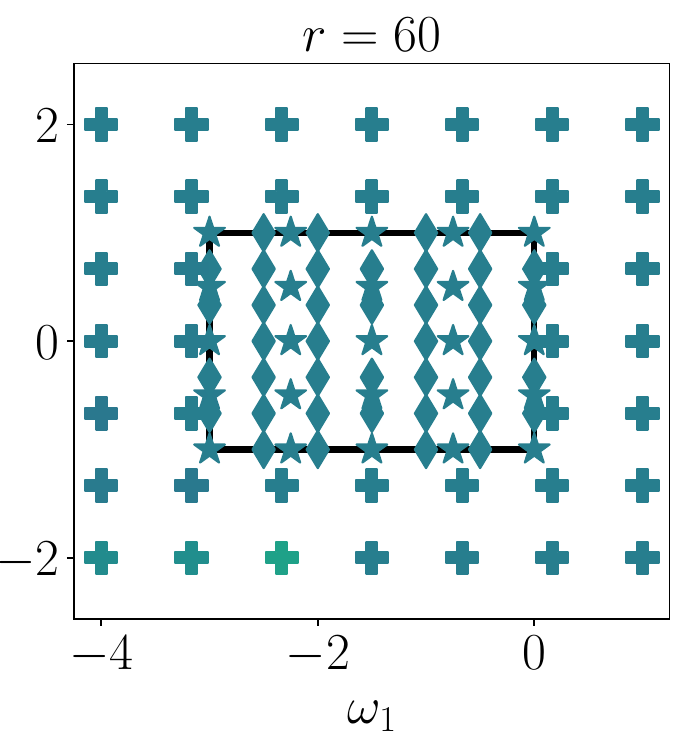}
        \includegraphics[width=0.295\textwidth, trim={1.25cm 0.25cm 0 0},clip]{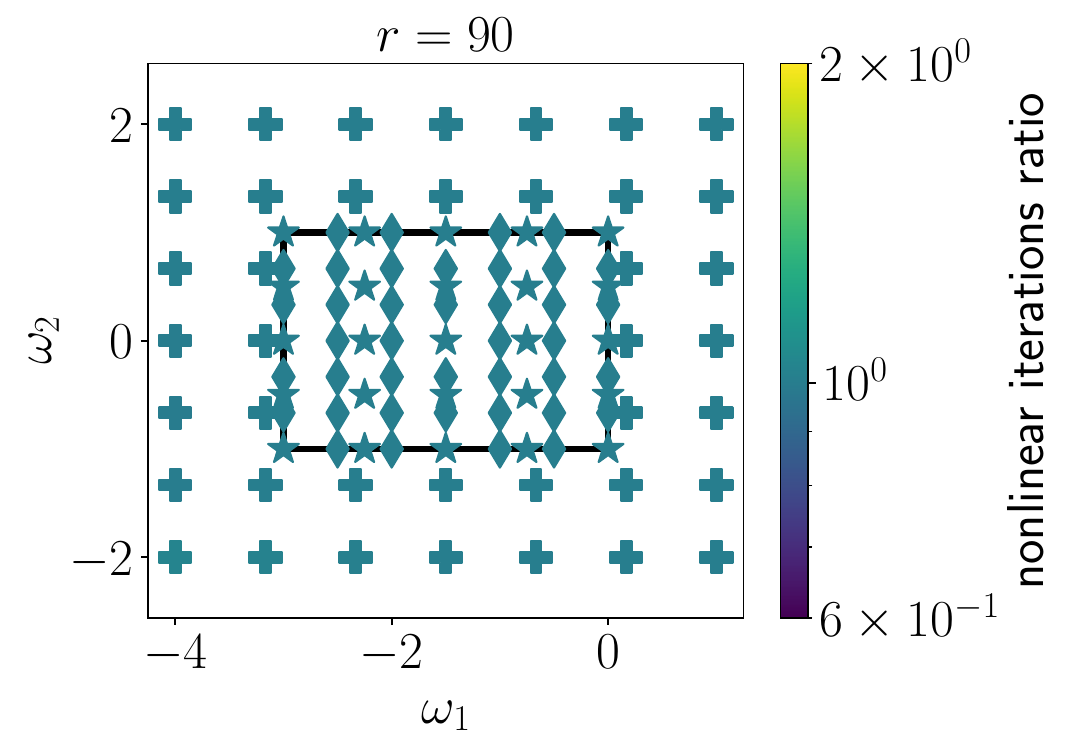}
        \caption{Ratio of the total number of FS-ROM iterations to that of FS-HF iterations.}
        \label{fig:application_iterations}
    \end{subfigure}\\

    \vspace{0.05in}
    \begin{subfigure}[t]{\textwidth}
        \centering
        \rotatebox{90}{\hspace{1.0cm}M-ROM}
        \includegraphics[width=0.237\textwidth, trim={0 1.25cm 3.5cm 0},clip]{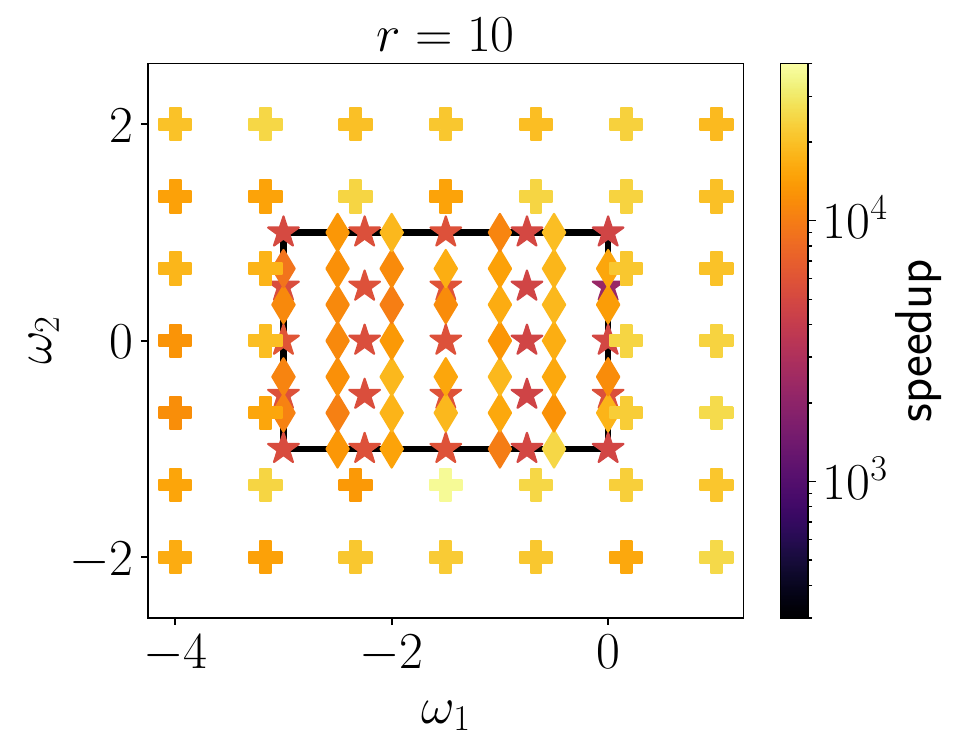}
        \includegraphics[width=0.214\textwidth, trim={1.25cm 1.25cm 3.5cm 0},clip]{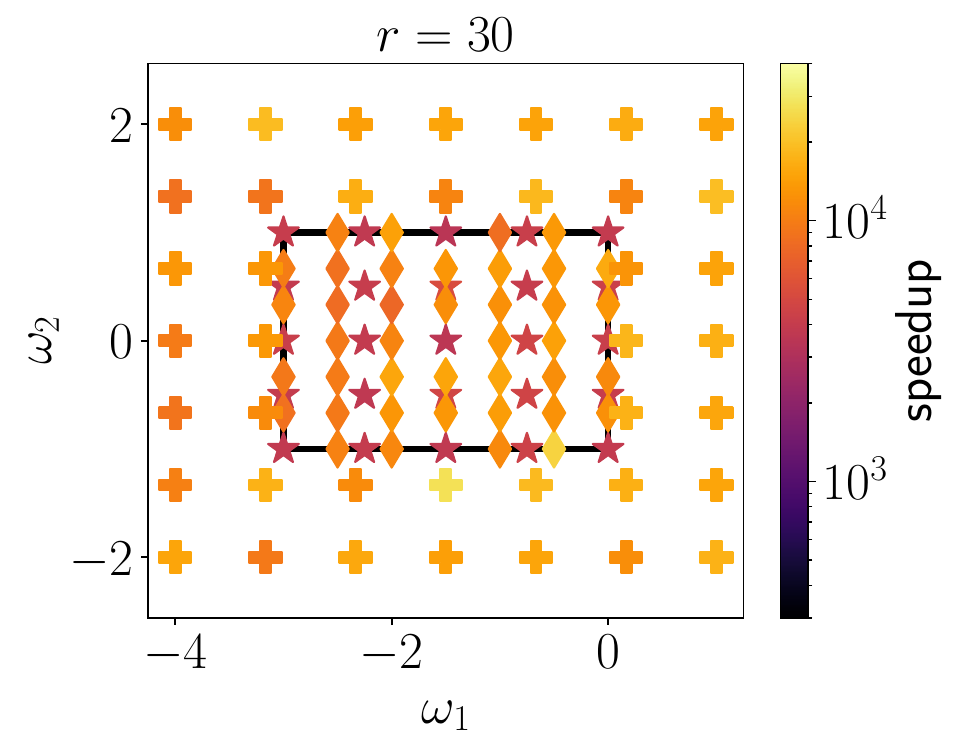}
        \includegraphics[width=0.214\textwidth, trim={1.25cm 1.25cm 3.5cm 0},clip]{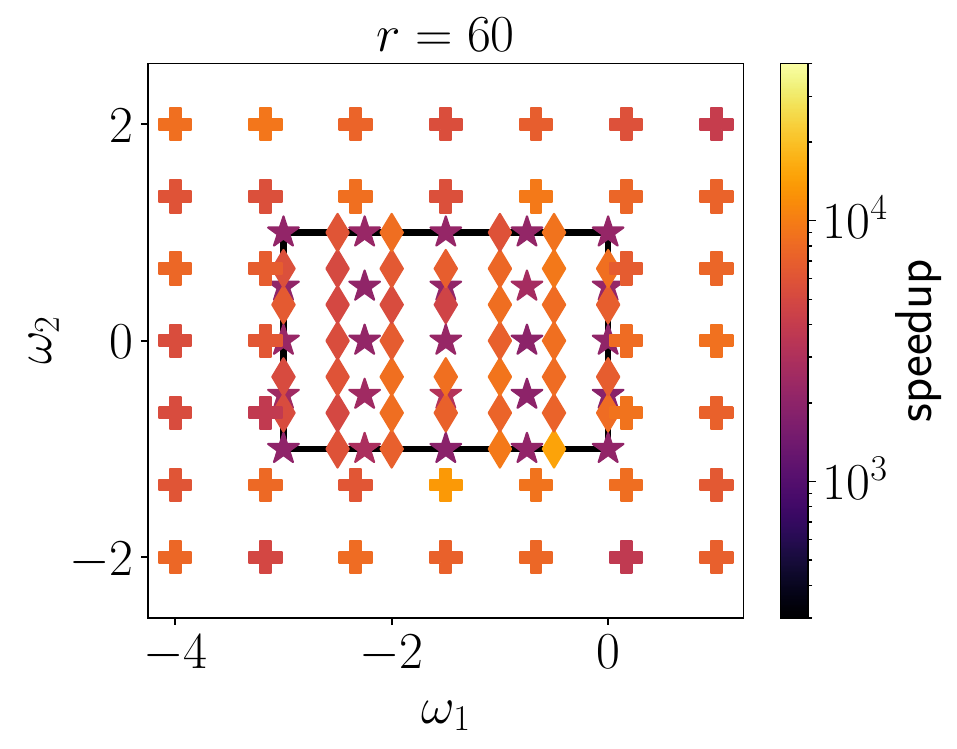}
        \includegraphics[width=0.28\textwidth, trim={1.25cm 1.25cm 0 0},clip]{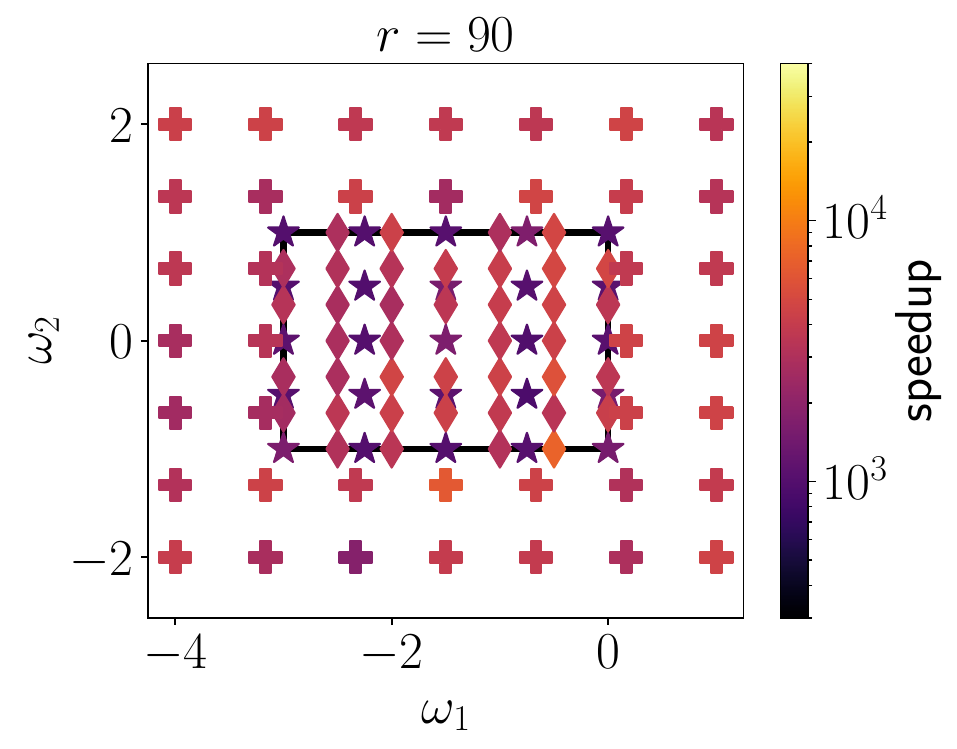}
        \\
        \rotatebox{90}{\hspace{1.1cm}FS-ROM}
        \includegraphics[width=0.237\textwidth, trim={0 0.25cm 3.5cm 0.15cm},clip]{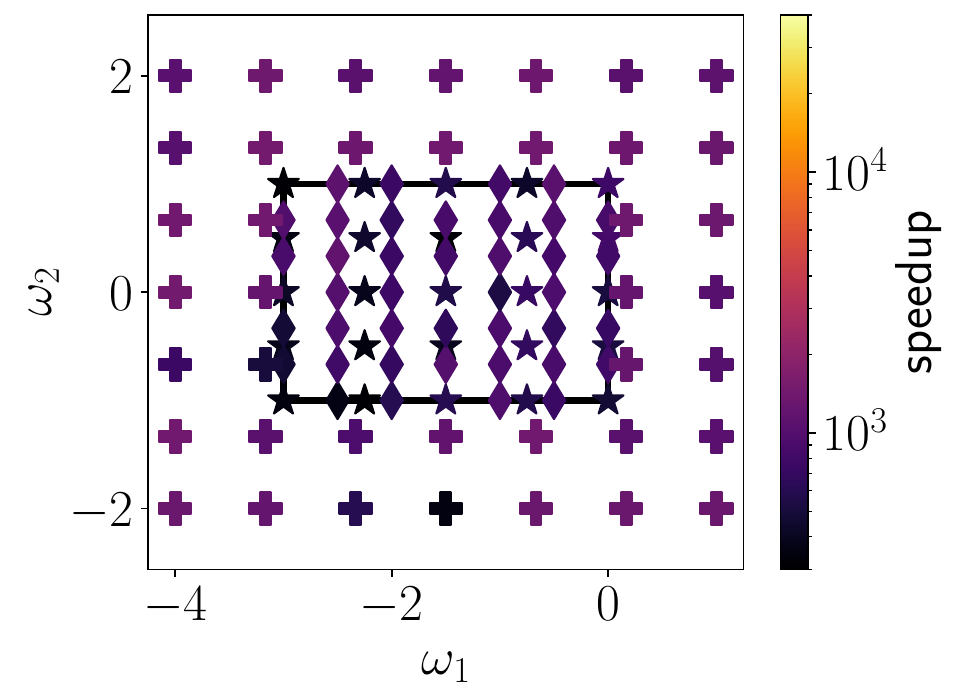}
        \includegraphics[width=0.214\textwidth, trim={1.25cm 0.25cm 3.5cm 0.15cm},clip]{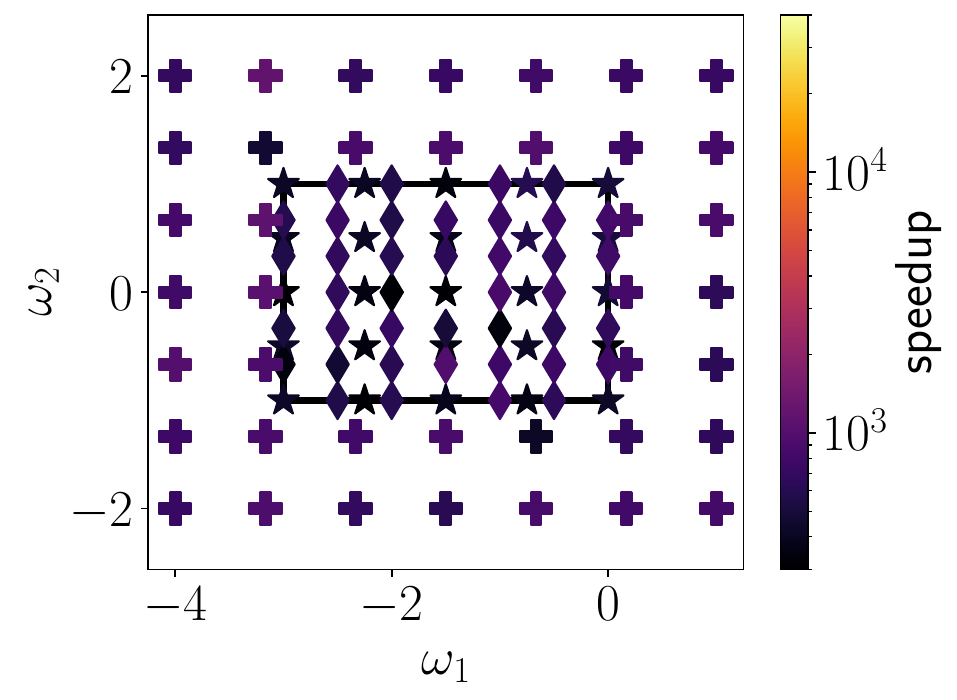}
        \includegraphics[width=0.214\textwidth, trim={1.25cm 0.25cm 3.5cm 0.15cm},clip]{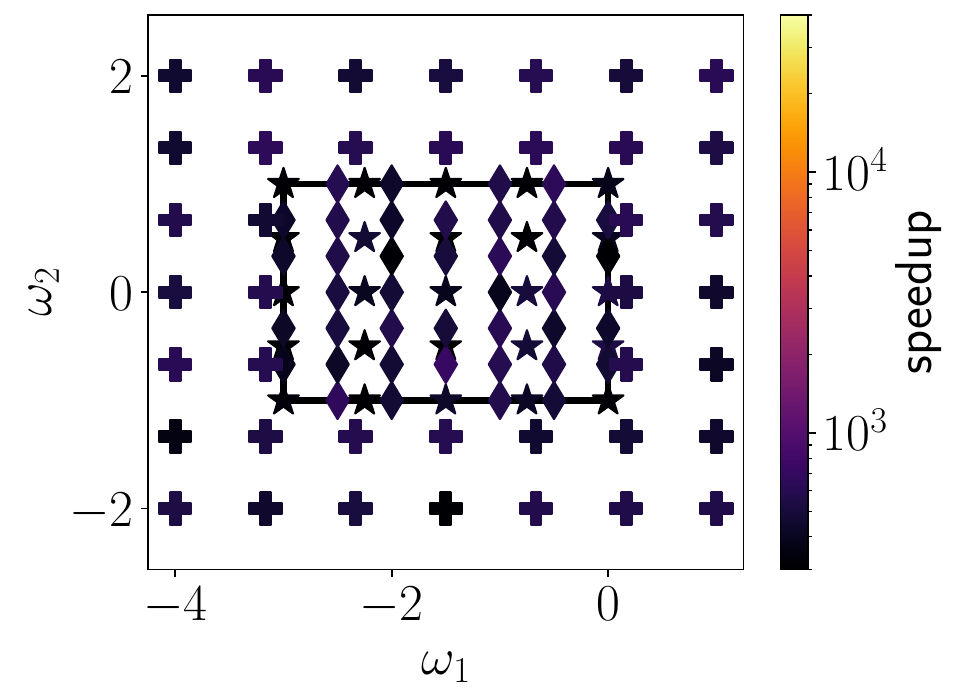}
        \includegraphics[width=0.28\textwidth, trim={1.25cm 0.25cm 0 0.15cm},clip]{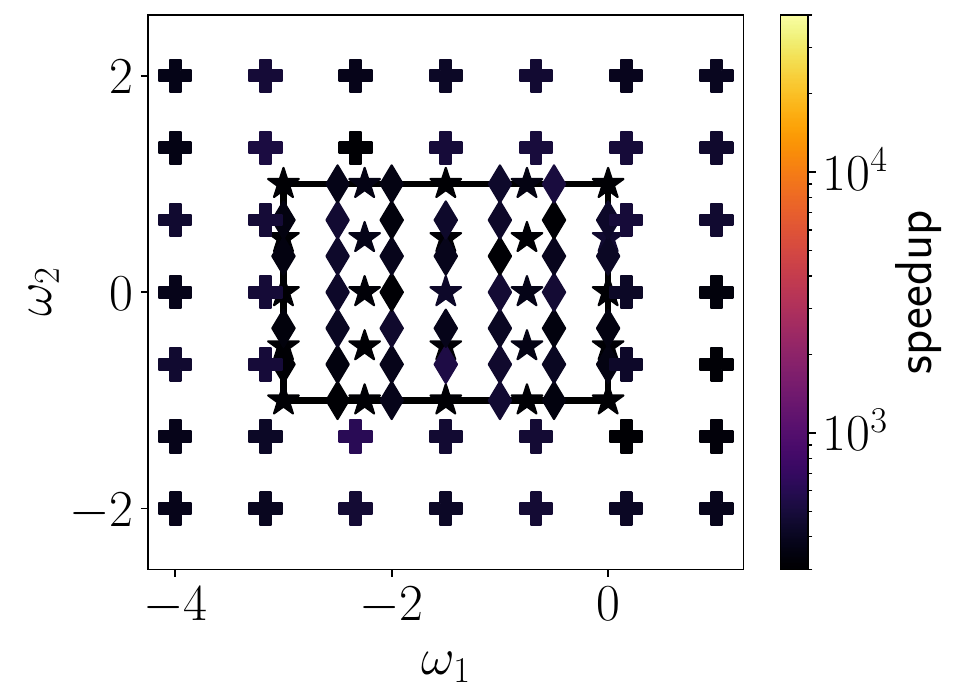}
        \caption{Speedup of the FS-ROM over FS-HF.}
        \label{fig:application_speedup}
    \end{subfigure}
    \caption{Example 2. Illustration of (a) ratio of the total iteration numbers of FS-ROM to FS-HF and (b) speedup of the ROM over HF for the given $\mathbb{P}_{\text{test}}$ ($\mbox{\fontsize{10pt}{\baselineskip}\selectfont\FiveStar}$ for case i), $\blacklozenge$ for case ii), \ding{58} for case iii)). The inner black rectangle denotes $\mathbb{P}_{\text{train}}$ as shown in Figure \ref{fig:ex3_setup_b}.}
\end{figure}

Figure \ref{fig:application} is divided into five panels, representing five different quantities of interest.
The first three panels, Figures \ref{fig:application_u_H1_relative_errors}--\ref{fig:application_theta_H1_relative_errors}, depict the $\normo{\cdot}$-norm relative errors for $\bur$, $\pr$, and $\thetar$ with respect to $\buh$, $\ph$ and $\thetah$, respectively.
In each panel, the top row displays the results obtained from the M-ROM, while the bottom row illustrates those from the FS-ROM. Each row contains four images, corresponding to increasing values of $r = 10$, $r = 30$, $r = 60$ and $r = 90$.
Markers are positioned at the provided $\bomega = (\omega_1, \omega_2)$ coordinates, signifying the errors for the three distinct cases: case i) denoted by $\mbox{\fontsize{10pt}{\baselineskip}\selectfont\FiveStar}$, case ii) by $\textstyle \blacklozenge$, and case iii)  by \ding{58}.

In the first column of Figures \ref{fig:application_u_H1_relative_errors}, \ref{fig:application_p_H1_relative_errors} and \ref{fig:application_theta_H1_relative_errors},
it is evident that both the M-ROM and HF-ROM with $r = 10$ exhibit inaccuracies. The relative errors are approximately $10^{-2}$ within the training range (inside the inner rectangular region) and can escalate to as much as $1$ outside of this range.
The accuracy progressively increases when increasing $r$ to $30, 60$ and $90$. In particular, when $r = 90$ we observe that both the M-ROM and HF-ROM become more accurate within the training range, including both seen (i.e., in the training set) and unseen data (i.e., not in the training set).
We note that, at $r = 90$, the temperature ($\theta$) exhibits the lowest error within the training range, reaching a maximum of $10^{-6}$ for both the M-ROM and FS-ROM. This outcome aligns with the observation that the POD eigenvalues associated with temperature were the ones that exhibited the most rapid decrease.
On the other hand, the displacement $\bu$ and pressure $p$ show differences of up to an order of magnitude between the M-ROM and FS-ROM for parameter values within the training range. Nevertheless, the resulting relative errors remain at most $10^{-4}$, a value that falls below the tolerance $\epsilon$ established by the FS-HF solver.

We also used both M-ROM and FS-ROM for extrapolation outside of the training range, which resulted in larger relative errors, especially in the bottom left of $\mathbb{P}_{\text{test}}$. The highest relative error is observed at $\bomega = (-4, -2)$ and is of the order of $10^{-2}$ for FS-ROM and $10^{-1}$ for M-ROM for the pressure approximation, i.e. the variable that was characterized by the slowest POD decay.

Finally, Figures~\ref{fig:application_iterations} and \ref{fig:application_speedup} discuss the number of iterations  and the computational efficiency of the ROM.
Specifically, Figure~\ref{fig:application_iterations} illustrates the ratio of the total number of FS-ROM iterations to FS-HF iterations. Apart from the smallest values, $r = 10$ and $r = 30$, in the extrapolation region, FS-ROM converges in a comparable number of iterations to FS-HF across the entire parameter range.
Also, Figure \ref{fig:application_speedup} depicts the speedup offered by the ROM, representing the ratio of CPU time required for a HF solve compared to a ROM solve. The magnitude of the speedup for the M-ROM diminishes from $10^4$ (for $r = 10$) to $10^3$ (for $r = 90$) as the reduced basis size $r$ increases. On the other hand, FS-ROM exhibits a smaller speedup, approximately on the order of $10^2$, which remains relatively consistent across varying $r$ values. In summary, any ROM yields query times that are at least 100 times faster than the corresponding HF scheme.

\section{Conclusions}

This paper introduces a novel approach involving a fixed-stress iterative  for the solution of {linear thermo-poroelasticity} problems in conjunction with reduced order modeling techniques. The approach is validated by means of several numerical examples, which also illustrate its computational capabilities. The benefits of this methodology are two-fold. Firstly, the utilization of fixed-stress iterations aids in the management of complex multi-physics coupling scenarios. Secondly, the incorporation of reduced order modeling significantly boosts computational efficiency. In future works we plan to demonstrate the versatility of the proposed techniques by employing them in conjunction with various numerical discretization techniques, including mixed finite elements, as well as discontinuous or enriched Galerkin finite element methods.

\section*{Acknowledgments}
Francesco Ballarin acknowledges the PRIN 2022 PNRR project ``ROMEU: Reduced Order Models for Environmental and Urban flows'' (CUP J53D23015960001), and the INdAM-GNCS projects ``Metodi numerici per lo studio di strutture geometriche parametriche complesse'' (CUP E53C22001930001) and ``Metodi di riduzione computazionale per le scienze applicate: focus su sistemi complessi'' (CUP E55F22000270001).
The work of Sanghyun Lee was supported by the US National Science Foundation under Grant DMS-2208402.
The work of Son-Young Yi was supported by the US National Science Foundation under Grant DMS-2208426.
\bibliographystyle{elsarticle-num}
\bibliography{references.bib}

\end{document}